\documentclass[10pt]{article}

\usepackage{graphicx}
\usepackage{amsmath,amsfonts}
\usepackage{sidecap}
\usepackage{subfig}
\usepackage{mathrsfs}

\usepackage{mathrsfs}
\usepackage{pifont} 
\usepackage{stmaryrd}
\usepackage{setspace}

\usepackage[normalem]{ulem}
\usepackage{hyphenat}

\usepackage{todonotes}

\usepackage{wasysym}

\newcommand{\dt}[1]{\textcolor{black}{#1}}

\newcommand{\revdt}[1]{\textcolor{black}{#1}} 

\def\N{\mathbb{N}}

\def\R{\mathbb{R}}

\newcommand{\eequ}{\end{equation}}
\newcommand{\bequ}{\begin{equation}}
\newcommand{\eequd}{\end{eqnarray*}}
\newcommand{\bequd}{\begin{eqnarray*}}

\def\Bila{\mathbf{B}}

\begin{document}

\title{A multiscale mathematical model of tumour invasive growth
}

\author{Lu Peng\thanks{Division of Mathematics,
              University of Dundee, Dundee, DD1 4HN, United Kingdom;
              Current Address:
              Beijing Computational Science Research Centre, 
              Division of Applied and Computational Mathematics,
             10 Dongbeiwang West Road,
             Haidian District,
              Beijing 100193, E-mail: penglu@csrc.ac.cn}
        \and Dumitru Trucu\thanks{Division of Mathematics,
              University of Dundee, Dundee, DD1 4HN, United Kingdom, E-mail: trucu@maths.dundee.ac.uk}
        \and Ping Lin\thanks{Division of Mathematics,
              University of Dundee, Dundee, DD1 4HN, United Kingdom, E-mail: plin@maths.dundee.ac.uk}
        \and Alastair Thompson\thanks{MD Anderson Cancer Center,
            The University of Texas, 1515 Holcombe Blvd, Houston, TX 77030, USA, E-mail: AThompson1@mdanderson.org}
        \and Mark A. J. Chaplain\thanks{School of Mathematics and Statistics, Mathematical Institute,
            University of St Andrews, St Andrews KY16 9SS, United Kingdom, E-mail: majc@st-andrews.ac.uk}}

\maketitle

\begin{abstract}

Known as one of the hallmarks of cancer \cite{Hanahan_et_al_2000}, cancer cell invasion of human body tissue is a complicated spatio-temporal multiscale process which enables a localised solid tumour to transform into a systemic, metastatic and fatal disease. This process explores and takes advantage of the reciprocal relation that solid tumours establish with the extracellular matrix (ECM) components and other multiple distinct cell types from the surrounding microenvironment. 
Through the secretion of various proteolytic enzymes such as matrix metalloproteinases (MMP) or the urokinase plasminogen activator (uPA), the cancer cell population alters the configuration of the surrounding ECM composition and overcomes the physical barriers to ultimately achieve local cancer spread into the surrounding tissue. 

The active interplay between the tissue-scale tumour dynamics and the molecular mechanics of the involved proteolytic enzymes at the cell-scale underlines the biologically multiscale character of invasion, and raises the challenge of modelling this process with an appropriate multiscale approach. 
In this paper, we present a new two-scale moving boundary model of cancer invasion that explores the tissue scale tumour dynamics in conjunction with the molecular dynamics of the urokinase plasminogen activation system. Building on the multiscale moving boundary method proposed in \cite{Dumitru_et_al_2013}, the modelling that we propose here allows us to study the changes in tissue scale tumour morphology caused by the cell-scale uPA micro-dynamics occurring along the invasive edge of the tumour. Our computational simulation results demonstrate a range of heterogeneous dynamics which are qualitatively similar to the invasive growth patterns observed in a number of different types of cancer, such as the tumour infiltrative growth patterns discussed in \cite{Ito_et_al_2012}.\\


\end{abstract}

\section{Introduction}
\label{intro}
Cancer is a complicated disease that involves many cross-related processes occurring over several spatial scales, ranging from genes to cells to tissues. The abilities of cancer cells to activate invasion and metastasis, to sustain proliferative signalling, to evade growth suppressors, to enable replicative immortality, to induce angiogenesis, and to resist cell death, have been initially identified as the six main hallmarks of cancer \cite{Hanahan_et_al_2000}. A growing knowledge about cancer over the last decade has shed more light on the whole picture of the disease and another four hallmarks were added, namely: the ability of cancer to avoid immune destruction, to deregulate cellular energetics, to develop tumour-promoting inflammations, alongside genome instability and mutations \cite{Hanahan_et_al_2011}.
 
Highlighted as one of the hallmarks of cancer, cancer cell invasion is a landmark event that transforms a locally growing tumour into a systemic, metastatic and fatal disease. The past four decades have witnessed great scientific efforts focussed on gaining a better understanding of the processes involved during cancer invasion, which is of highly importance in designing early detection strategies and attempting effective therapies. 

A malignant tumour includes a complex heterotypic community of cells (such as cancer cells, immuno-inflamatory cells, stromal cell, fibroblasts, endothelial cells, macrophages) that are mixed with ECM. This community is enhanced by vastly complex signalling pathways underpinning intense molecular processes that mediates the crosstalk between the various cell populations composing the tumours, such as the interaction between the cancer cells and the peritumoural stroma cells occurring during cancer invasion \cite{Hanahan_et_al_2011,Egeblad_et_al_2010,Qian_Pollard_2010,Joyce_Pollard_2009,Kalluri_Zeisberg_2006}. As one of the main factors that affect the way the cancer cell migrate and invade, the peritumoural ECM not only plays the role of a scaffold for the tissues and physical barriers during cell migration but also conveys the signalling pathway processes, enabling the cells to communicate. These give rise to specific conditions within the tumour microenvironment that locally regulate cell migration, proliferation and differentiation. Specifically, the secretion of proteolytic enzymes such as urokinase plasminogen activator (uPA) and matrix metalloproteinases (MMPs) by the tumour cells and interactions of these enzymes with the ECM components lead to proteolytic degradation and remodelling of the ECM and represent a key step in the cancer invasion process.

In order to decipher the mechanisms behind the complicated processes involved in cancer invasion, clinical investigations and experimental observations carried out over the past few decades have started being increasingly accompanied by mathematical modelling \cite{Adam_1986,Anderson_et_al_2000,Byrne_Chaplain_1996,Chaplain_et_al_2006,Gatenby_Gawlinski_1996,Greenspan_1976,Perumpanani_et_al_1996,Perumpanani_et_al_1998,Webb_et_al_1999}. These initial modelling attempts used reaction-diffusion systems to describe the interaction between malignant and normal cells and focused on several important invasion aspects, leading to the first qualitative mathematical modelling approaches and paving the way towards a better understanding of the contribution of proteolytic enzymes in cancer invasion. This was further explored with modelling focused the role of proteolytic activities of some specific hydrolytic enzymes such as urokinase plasminogen activator (uPA) and matrix metalloproteinases (MMPs) in tumour invasive behaviour \cite{Andasari_et_al_2011,Anderson_2005,Byrne_et_al_2001,Chaplain_Lolas_2005,Deakin_Chaplain_2013,Perumpanani_et_al_1996,Perumpanani_et_al_1998}. For instance, while \cite{Chaplain_Lolas_2005} proposed a system of reaction-diffusion-taxis partial differential equations to explore the role of the uPA system (including uPA, uPA-inhibitors, plasmin and host tissue) in cancer invasion, the model introduced in \cite{Deakin_Chaplain_2013} evaluates the role of two types of simultaneously expressed MMPs on cancer growth and spread i.e., membrane-bound matrix MT1-MMP and soluble MMP-2.  

Alongside the chemotactic and haptotactic movement assumed in all the models mentioned above, cell-cell and cell-matrix adhesion were also recognised as playing a crucial role in the growth and development of carcinomas \cite{Byrne_Chaplain_1996}. Advances in addressing the importance of cell-adhesion in the cancer invasion process were obtained via several models based on systems of non-local integro-differential equations \cite{Armstrong_et_al_2006,Chaplain_et_al_2011,Domschke_et_al_2014,Gerisch_Chaplain_2008} that account for dynamic interactions within an appropriately small sensing radius $R$ between potentially mutating cancer cell populations, ECM, and the involved matrix degradation enzymes. 

Since the local tissue invasion of a malignant tumour could be regarded also as a free moving boundary problem, several appropriate numerical techniques were considered in the computational modelling of cancer invasion. Particularly, the level-set method was intensely used to study solid tumour growth in homogeneous microenvironments \cite{Frieboes_et_al_2006,macklin05,macklin06,macklin07,Zheng_et_al_2005}. A new ghost cell/level set method (based of a nonlinear nutrient equation coupled with a pressure equation with geometry-dependent jump boundary conditions) was developed in \cite{macklin08} and applied to models of tumour invasive growth in complex, heterogeneous tissues. This model was later extended into an improved model of tumour invasion including the process of tumour-induced angiogenesis \cite{Macklin_et_al_2009}. Finally, alongside the level-set method, multiphase models based on the theory of mixtures were also developed and used to investigate tumour growth and spread \cite{Byrne_Preziosi_2003,Chaplain_et_al_2006bb,Frieboes_et_al_2010,Preziosi_Tosin_2009,Wise_et_al_2008,Wise_et_al_2011}. In these approaches, the tumour was regarded as a system consisting of different phases (e.g. cellular phase, liquid phase, etc) and the development of a solid tumour was modelled by exploring the mass and momentum balances alongside the constitutive laws that distinguish the phases in the system. 

However, while recognizing the multiscale nature of cancer growth and spread, over the past two decades or so most computational and mathematical modelling has focused mainly on one scale, either at the tissue, cell or molecular scale, with the first attempts towards linking these scales being revised in \cite{Deisboeck2011}. Recently, following a series of important developments within the general multiscale framework, and based on strong insights from atomistic-to-continuum methods \cite{Lin2007}, homogenization techniques \cite{Allaire1992,Trucu2012}, and heterogeneous multiscale finite element methodology \cite{Abdulle2005,EEngquist2003,Ren2005}, a genuinely new multiscale moving boundary model for cancer invasion that links the tissue, cellular and subcellular scales was proposed in \cite{Dumitru_et_al_2013}. \revdt{In this new framework, PDE modelling at the tissue scale for the cell population dynamics and PDE modelling at the cell scale for the molecular mechanics of the proteolytic enzymes population are linked together in a two-scale model through top-down and bottom-up permanent links. The top-down link provides the source for the micro-scale dynamics, which is induced in a non-local manner by the macro-dynamics. On the other hand, the bottom-up link enables the micro-dynamics to provide the macro-dynamics with a law for the macro-scale boundary movement, whose direction and displacement magnitude is determined at microscale. This is fundamentally different from previous modelling perspectives such as the one-scale modelling presented in \cite{Anderson_et_al_2000}, or the one proposed in \cite{Ramis-Conde_et_al_2008} where the authors considered continuous modelling at the microscale but an individual-based model at the macroscale and where no top-down links were assumed.}

\revdt{In this paper we propose a novel multiscale mathematical model of cancer invasion} that explores the tissue scale cancer progression in conjunction with the cell-scale dynamics of the urokinase plasminogen activation system. Building on the multiscale moving boundary method proposed in \cite{Dumitru_et_al_2013}, the new modelling that we propose here allows us to study the changes in tissue scale (macro) tumour morphology caused by the cell-scale uPA micro-dynamics occurring in a cell-scale (micro) neighbourhood  of the invasive edge of the tumour. Assuming the uPA model developed in \cite{Andasari_et_al_2011,Chaplain_Lolas_2005} for the tissue scale dynamics, we derive governing laws for the degrading enzymes cell-scale dynamics arising in the close proximity of the tumour interface. While this results in prescribing an appropriate form for the top-down link between the macro- and micro- dynamics, by exploring the spatial interaction between the uPA micro-dynamics and the surrounding ECM from the peritumoural region we obtain a bottom-up link between the micro- and macro- dynamics. This allows us to describe the evolution of tumour invasive edge morphology and enables computational predictions for the changes occurring in the macroscopic pattern of cancer during the local invasion process.

\section{Biological Background: Components of the uPA System and their Functions}
Proteolytic degradation and remodelling of the extracellular matrix is essential for cancer cell invasion. It enables cancer cells to proliferate and migrate through surrounding tissue. In this context, one of the first steps of invasion is the production and secretion of proteolytic enzymes, i.e., urokinase plasminogen activator (uPA) and matrix metalloproteinases (MMPs) by cancer cells. These enzymes interact with the dynamics of the ECM macromolecules and pave the way for cancer invasion. Specifically, the uPA enzymatic system mainly consists of the urokinase receptor (uPAR), urokinase plasminogen activator (uPA), the matrix-like protein vitronectin (VN), plasminogen activator inhibitor type1 (PAI-1), and the degrading enzyme plasmin.

\paragraph{\it{Urokinase plasminogen activator (uPA).}}
uPA is an extracellular serine protease produced by cells. Two major functional domains of the uPA molecule are the protease domain and the growth factor domain. The protease part activates plasminogen and turns it into plasmin, which is able to digest basement membrane and extracellular matrix proteins. The growth factor domain has no protease activity but can bind a specific high affinity cell-surface receptor, uPAR. Finally, uPA has a zymogen form, pro-uPA, which can be activated by plasmin and binds to uPAR.
\paragraph{\it{Urokinase plasminogen activator receptor (uPAR).}}
uPAR is a high affinity cell-surface receptor of uPA (and of its zymogen form pro-uPA), which via the binding process localises the uPA and pro-uPA to the cell surface. Importantly, uPAR contains another binding site for the ECM component called vitronectin (VN), and since VN and uPA binding sites are distinct, uPAR can simultaneously bind both ligands, allowing coordinated regulation of proteolysis, cell adhesion, and signalling.

uPAR expression during ECM remodelling is well-controlled under normal conditions, for example, in gestational tissues during embryo implantation and placental development and in keratinocytes during epidermal would healing. uPAR is also expressed in many human cancers. It indicates poor prognosis and in some cases is predictive of invasion and metastasis. Importantly, uPAR expression in tumours can occur in tumour cells and/or tumour-associated stromal cells, such as fibroblasts and macrophages. Moreover, there is a certain crosstalk between these two binding processes, as the ligand-binding of uPA to uPAR enhances the VN binding by uPAR \cite{Harvey_Chris_2010}.
\paragraph{\it{Vitronectin (VN).}}
VN is an abundant versatile glycoprotein found in serum and the ECM and promotes cell adhesion and spreading. Vitronectin binds strongly to glass surfaces, as the name indicates ({\em vitro} = glass), and it has binding sites for several ligands, including heparin, urokinase plasminogen activator receptor (uPAR), plasminogen activator inhibitor type-1(PAI-1), and integrins, such as $\alpha_v\beta_3$. When vitronectin binds to uPAR, it is thought to bring PAI-1 closer to uPA, thereby promoting inhibition and clearance of uPA from the receptor.
\paragraph{\it{Urokinase plasminogen activator Inhibitor-1(PAI-1)}}
One of the inhibitors of urokinase plasminogen activator, PAI-1, belongs to the serpin (serine protein inhibitors) family and it is believed to be the most abundant, fast-acting inhibitor of uPA \textit{in vivo}. It can specifically bind to soluble and membrane-bound uPA to inhibit plasminogen activation. When PAI-1 binds to the uPA/uPAR complex, it triggers the internalisation of the uPA/uPAR/PAI-1 complex by receptor-mediated endocytosis, meaning that the complex will be dissociated and PAI-1 and uPA will be digested, but the receptor will be recycled to the cell surface. This process helps with the clearance of PAI-1 from the vicinity of the cell surface. Additionally, as a major binding protein of VN, PAI-1 competes with uPAR for binding to VN.
\paragraph{\it{Plasmin.}}
Plasmin is a widespread enzyme that cleaves many extracellular matrix proteins, such as fibronectin, laminin, vitronectin and thrombospondin. In addition, plasmin can also activate many matrix metalloproteinases (MMPs), enhancing even more the degradation of extracellular matrix. It can also influence the composition of the extracellular environment by affecting the activity of cytokines and growth factors, for example, \revdt{decreasing the activation of TGF-$\beta 1$ \cite{Venkatraman2012}}.


\section{Mathematical Modelling}
\label{sec:Math model}
In this section, we will detail the three main components of the uPA multiscale model that we propose for cancer invasion, namely: the macroscopic dynamics, the microscopic dynamics and regulation of the tumour boundary relocation. Thus, the description will cover the model at both the macro level (tissue-scale) and the micro level (cell-scale), and will explore the link between these two biological scales.
\subsection{The macroscopic dynamics}\label{subsec:macroscopic}
\revdt{For the tumour macroscopic dynamics, we adopt here the modelling hypotheses formulated in \cite{Chaplain_Lolas_2005}.} We denote the cancer cell density by $c$, the extracellular matrix density by $v$ (without making the distinction between ECM and its component VN), the urokinase plasminogen activator (uPA) concentration by $u$, the plasminogen activator inhibitor (PAI-1) concentration by $p$ and the plasmin concentration by $m$. Further, since we assume a fixed average number of receptors uPAR located on each cancer cell surface, there is no explicit modelling of uPAR. Therefore, the concentration of uPAR is considered to be proportional to the cancer cell density. Another important assumption is that the supply of plasminogen is unlimited in this model. Finally, the macroscopic model is obtained by accounting for the biological considerations described in previous section in conjunction with the following presumptions: 

\paragraph{\it{The cancer cell dynamics.}} It is assumed that cancer cell migration is mainly governed by diffusion, chemotaxis due to uPA, and PAI-1 and haptotaxis due to VN and other ECM components. Additionally, a logistic growth law is used to model cancer cell proliferation. Thus, the mathematical equation for cancer cell density that is considered here is as follows:
\begin{equation}
\frac{\partial c}{\partial t} = \underbrace{D_c \Delta c}_{\text{diffusion}} - \nabla \cdot [ \underbrace{\chi_u c \nabla u}_{\text{uPA-chemo}} + \underbrace{\chi_p c \nabla p}_{\text{PAI-1-chemo}} + \underbrace{\chi_v c \nabla v}_{\text{VN-hapo}} ] + \underbrace{\mu_1 c(1-\frac{c}{c_0})}_{\text{profieration} },
\end{equation}
where $D_c$ is the diffusion coefficient of cancer cells, $\chi_u$ and $\chi_p$ are the chemotaxis coefficients relevant to uPA and PAI-1 respectively, $\chi_v$ is the VN-mediated haptotaxis rate, $\mu_1$ is the cancer cell proliferation rate, and $c_0$ is the maximum carrying capacity for cancer cells.

\paragraph{\it{The ECM/VN dynamics.}}
As ECM is not assumed to move, we rule out any migration terms in the governing law. While VN (which is an important ECM component) is degraded in contact with enzymes $m$, the binding of PAI-1 to uPA inhibits the activation of plasminogen, leading to the protection of VN and other ECM molecular constituents and indirectly contributing to their production. Simultaneously, the binding of PAI-1 to VN results in less binding to cell-surface receptors such as uPAR, and so, through the regulation of cell-matrix-associated signal transduction pathways, this inhibits the production of VN. Therefore, assuming a logistic ECM remodelling, the governing equation for ECM is given by 
\begin{equation}
\frac{\partial v}{\partial t}= -\underbrace{\delta v m}_{\text{degradation}} +\underbrace{\phi_{21} u p}_{\text{uPA/PAI-1}} -\underbrace{\phi_{22} v p}_{\text{PAI-1/VN} }+ \underbrace{\mu_2 v(1-\frac{v}{v_0})}_{\text{remodelling}}.
\end{equation}
where $\delta$ is the rate of ECM degradation by plasmin, $\phi_{21}$ is the binding rate of PAI-1 to uPA, $\phi_{22}$ is the binding rate of PAI-1 to VN, and $\mu_2$ is the matrix remodelling rate.

\paragraph{\it{The uPA dynamics.}} The mathematical modelling of the uPA concentration dynamics accounts for the following aspects. While being produced by the cancer cells and removed from the system due to its binding with PAI-1 and uPAR, per unit time the uPA exercises a local diffusion. Therefore, this can be formalised mathematically as follows:

\begin{equation}
\frac{\partial u}{\partial t} = \underbrace{D_u \Delta u}_{\text{diffusion} }- \underbrace{\phi_{31} p u}_{\text{uPA/PAI-1}} -\underbrace{\phi_{33} c u}_{\text{uPA/uPAR}} +\underbrace{\alpha_{31} c}_{\text{production}},
\end{equation}
where $D_u$ is the diffusion coefficient, $\phi_{31}$ and $\phi_{33}$ are binding rates of uPA/PAI-1 and uPA/uPAR accordingly, and $\alpha_{31}$ is the production rate of uPA by the cancer cells.

\paragraph{\it{The PAI-1 dynamics.}} Similarly, the equation for PAI-1 simply includes a diffusion term with coefficient $D_p$, removal caused by binding to uPA and VN with binding rates $\phi_{41}$ and $\phi_{42}$ respectively, and production as a result of plasmin formation at a rate $\alpha_{41}$. Thus, these considerations lead us to the following governing equation:
\begin{equation}
\frac{\partial p}{\partial t} = \underbrace{D_p \Delta p}_{\text{diffusion}} - \underbrace{\phi_{41}p u}_{\text{uPA/PAI-1}} -\underbrace{\phi_{42} p v}_{\text{PAI-1/VN}} +\underbrace{\alpha_{41} m}_{\text{production}}.
\end{equation}

\paragraph{\it{The plasmin dynamics.}} The evolution of plasmin concentration is modelled as follows. While assuming that per unit time this exercises a local diffusion, we consider that the binding of uPA to uPAR provides an opportunity for pericellular proteolytic activity through plasminogen activation leading to plasmin formation. Moreover, the binding of PAI-1 to VN indirectly enhances the binding of uPA to uPAR, therefore bringing additional contribution to plasmin formation. Thus, these assumptions give us the following evolution law:

\begin{equation}
\frac{\partial m}{\partial t}  = \underbrace{D_m \Delta m}_{\text{diffusion}} + \underbrace{\phi_{52} p v}_{\text{PAI-1/VN}}+\underbrace{\phi_{53} c u}_{\text{uPA/uPAR}}-\underbrace{\phi_{54} m}_{\text{degradation}}, 
\end{equation}
where $D_m$ is the diffusion coefficient, $\phi_{52}$ and $\phi_{53}$ are the binding rates of PAI-1/VN and uPAR/uPA accordingly, and $\phi_{54}$ is the rate of decay of plasmin.

To summarise the macro-dynamics, the dimensionless mathematical model of the uPA system adopted here is the one that was initially proposed in \cite{Chaplain_Lolas_2005}, namely:
\begin{align}
\frac{\partial c}{\partial t} &= \underbrace{D_c \Delta c}_{\text{diffusion}} - \nabla \cdot [ \underbrace{\chi_u c \nabla u}_{\text{uPA-chemo}} + \underbrace{\chi_p c \nabla p}_{\text{PAI-1-chemo}} + \underbrace{\chi_v c \nabla v}_{\text{VN-hapo}} ] + \underbrace{\mu_1 c(1-c)}_{\text{profieration} } \; ,\label{equ:Cancer}\\
\frac{\partial v}{\partial t}&= -\underbrace{\delta v m}_{\text{degradation}} +\underbrace{\phi_{21} u p}_{\text{uPA/PAI-1}} -\underbrace{\phi_{22} v p}_{\text{PAI-1/VN} }+ \underbrace{\mu_2 v(1-v)}_{\text{remodelling}} \;, \label{equ:ECM} \\
\frac{\partial u}{\partial t} &= \underbrace{D_u \Delta u}_{\text{diffusion} }- \underbrace{\phi_{31} p u}_{\text{uPA/PAI-1}} -\underbrace{\phi_{33} c u}_{\text{uPA/uPAR}} +\underbrace{\alpha_{31} c}_{\text{production}}, \label{equ:uPA_macro}\\
\frac{\partial p}{\partial t} &= \underbrace{D_p \Delta p}_{\text{diffusion}} - \underbrace{\phi_{41}p u}_{\text{uPA/PAI-1}} -\underbrace{\phi_{42} p v}_{\text{PAI-1/VN}} +\underbrace{\alpha_{41} m}_{\text{production}} , \label{equ:PAI-1_macro} \\
\frac{\partial m}{\partial t} & = \underbrace{D_m \Delta m}_{\text{diffusion}} + \underbrace{\phi_{52} p v}_{\text{PAI-1/VN}}+\underbrace{\phi_{53} c u}_{\text{uPA/uPAR}}-\underbrace{\phi_{54} m}_{\text{degradation}} . \label{equ:plasmin_macro}
\end{align}

\subsection{The microscopic dynamics}
\label{subsec:microscopic}

Turning now our attention to the micro-scale setting, in the following we will derive and propose a system of three coupled PDEs to describe the microdynamics of the plasminogen activation system taking place within a cell-scale $\epsilon$-neighbourhood $\mathcal{P}_{\epsilon}$ of the tumour invasive edge $\partial\Omega(t_{0})$ (which is introduced in Appendix \ref{sec:Appendix}). 

Assuming that PAI-1 and uPAR are uniformly expressed on the cell surface of various cell types in the tumour, this system will capture the leading edge micro-dynamics of the uPA, PAI-1 and plasmin by accounting for the following biological considerations. On one hand, the urokinase plasminogen activator (uPA) is assumed to bind to the cancer surface receptor uPAR to activate plasminogen, leading to degradation of pericellular ECM through a series of proteolytic activities. On the other hand, the membrane-bound MMPs (such as MT1-MMP) are secreted from within the tumour cell population distributed on the outer proliferating rim along the entire tumour periphery. Their region of proteolytic activities is therefore restricted around the tumour interface \cite{Deakin_Chaplain_2013,Farideh_et_al_2009}. Thus, based on these considerations, we propose a coupled governing law for the leading edge microdynamics, which is detailed as follows:

\paragraph{\it{The uPA microdynamics.}}

In each micro region $\epsilon Y$ the dynamics of the uPA molecular population is governed by a diffusion process whose source is induced from the tumour cell macro-dynamics. At each point $y\in\epsilon Y$, a source of uPA arises as a collective contribution of the tumour cells distributed within a certain neighbouring area within the tumour's outer proliferating rim. Therefore, this source is denoted by $f_{\epsilon Y}(\cdot,\cdot)\times[0,\Delta t]:\epsilon Y\rightarrow \mathbb{R}_+$ and is defined by:
\begin{equation}
f_1^{\epsilon Y} (y, \tau) =
\left\{\begin{array}{ll}
 \frac{1}{\lambda(B(y,\gamma)\cap \Omega(t_0))} \int\limits_{B(y,\gamma)\cap \Omega(t_0)} c\;(x,t_0+\tau) \;dx, & y \in \epsilon Y \cap \Omega(t_{0}),\\
 0, & \textrm{outside cancer},
\end{array}
\right.
\label{eq:source1_chap5}
\end{equation}
where $\lambda(\cdot)$ is the standard Lebesgue measure, and $\gamma$ represents the maximal thickness of the outer proliferating rim.
Thus, per unit time, under the presence of source \eqref{eq:source1_chap5} the uPA is locally diffusing and is binding to both PAI-1 and uPAR, and so its microdynamics can be formally written as: 
 
\begin{equation}
\frac{\partial u}{\partial \tau} = \underbrace{D_u \Delta u}_{\text{diffusion} }- \underbrace{\phi_{31} p u}_{\text{uPA/PAI-1}} + \;(\!\!\!\! \underbrace{\alpha_{31} }_{\text{production}} -\underbrace{\phi_{33}  u}_{\text{uPA/uPAR}}\!\!\!\!)  f_1^{\epsilon Y} (y, \tau) \label{equ:uPA}
\end{equation}

\paragraph{\it{The PAI-1 microdynamics.}}

The equation for PAI-1 accounts for diffusive motion, production due to plasmin activation, and loss due to binding with uPA and VN. Specifically, the binding between PAI-1 and VN is as a collective effect of the ECM distribution within $\epsilon Y$. Therefore, proceeding similarly to the case of the source term in \eqref{eq:source1_chap5}, we define

\begin{equation}
f_2^{\epsilon Y} (y, \tau) = \frac{1}{\lambda(B(y,2\epsilon))} \int_{B(y,2\epsilon)} v\;(x,t_0+\tau) \;dx,\;\;\;y \in \epsilon Y, 
\label{eq:source2_chap5}
\end{equation}
which finally enable us to write the following governing law for the PAI-1 microdynamics, namely:
\begin{equation}
\frac{\partial p}{\partial \tau} = \underbrace{D_p \Delta p}_{\text{diffusion}} - \underbrace{\phi_{41}p u}_{\text{uPA/PAI-1}} -\underbrace{\phi_{42} p\; f_2^{\epsilon Y} (y, \tau)}_{\text{PAI-1/VN}} +\underbrace{\alpha_{41} m}_{\text{production}}.   \label{equ:PAI-1}
\end{equation}

\paragraph{\it{The plasmin microdynamics.}} For the spatio-temporal evolution of plasmin, it is assumed that, per unit time, this exercises a local diffusion in the presence of the following source and decay circumstances. Considering that the binding of uPA to uPAR is required to provide the cell surface with a potential proteolytic activity, the plasmin source accounts on one hand on the contribution of the binding uPA/uPAR. On the other hand, as PAI-1 collectively competes with uPAR for binding to VN, the binding of PAI-1 to VN gives more opportunities to uPAR to bind with uPA,  and indirectly results in more plasmin formation. Finally, plasmin can be deactivated either by degradation or by the action of the plasmin inhibitor $\alpha_2$-antiplasmin. Thus, the equation that we obtain to describe these biological interactions is:

\begin{equation}
\frac{\partial m}{\partial \tau}  = \underbrace{D_m \Delta m}_{\text{diffusion}} + \underbrace{\phi_{52} p \; f_2^{\epsilon Y} (y, \tau)}_{\text{PAI-1/VN}}+\underbrace{\phi_{53} u \; f_1^{\epsilon Y} (y, \tau) }_{\text{uPA/uPAR}}-\underbrace{\phi_{53} m}_{\text{degradation}}  \label{equ:plasmin}
\end{equation}

In summary, the leading edge microdynamics is therefore given by the following system: 
\begin{align}
\frac{\partial u}{\partial \tau} &= \underbrace{D_u \Delta u}_{\text{diffusion} }- \underbrace{\phi_{31} p u}_{\text{uPA/PAI-1}} + \;(\!\!\!\! \underbrace{\alpha_{31} }_{\text{production}} -\underbrace{\phi_{33}  u}_{\text{uPA/uPAR}}\!\!\!\!)  f_1^{\epsilon Y} (y, \tau), \label{eq:uPA_chap5}\\ 
\frac{\partial p}{\partial \tau} &= \underbrace{D_p \Delta p}_{\text{diffusion}} - \underbrace{\phi_{41}p u}_{\text{uPA/PAI-1}} -\underbrace{\phi_{42} p\; f_2^{\epsilon Y} (y, \tau)}_{\text{PAI-1/VN}} +\underbrace{\alpha_{41} m}_{\text{production}}\!\!\!,  \label{eq:PAI-1_chap5}\\
\frac{\partial m}{\partial \tau} & = \underbrace{D_m \Delta m}_{\text{diffusion}} + \underbrace{\phi_{52} p \; f_2^{\epsilon Y} (y, \tau)}_{\text{PAI-1/VN}}+\underbrace{\phi_{53} u \; f_1^{\epsilon Y} (y, \tau) }_{\text{uPA/uPAR}}-\underbrace{\phi_{54} m}_{\text{degradation}}\!\!\!\!\!.  \label{eq:plasmin_chap5}
\end{align}

\subsection{The macroscopic tumour boundary relocation induced by the leading edge micro-dynamics }
Following the multiscale approach described in Appendix \ref{sec:Appendix}, the set of points $\{x_{\epsilon Y}^* \}$ on the boundary of tumour at the current time moves towards a set of new spatial positions  $\{\widetilde{x_{\epsilon Y}^* }\}$ to form the new boundary at the next multiscale-stage, provided that the local transitional probability $q^{*}$ is in agreement with the circumstances in the surrounding peritumoural microenvironment. When the invading strength is above a tissue threshold, the point $x_{\epsilon Y}^* $ will relocate to a new position $\widetilde{x_{\epsilon Y}^*}$ following a direction and displacement magnitude that represents the choreographic movement of all the points from the part of the invasive edge captured by the micro-domain $\epsilon Y$.   

As described in Appendix \ref{sec:The_MovingBoundaryMethod}, on any micro-domain $\epsilon Y$, provided that a sufficient amount of plasmin has been produced across the invading edge, it is the pattern of the front of the advancing spatial distribution of plasmin that characterises ECM degradation. Therefore, the movement direction and displacement magnitude of the part of the invading edge of the tumour caputred by the current micro-domain $\epsilon Y$ will be determined by the spatial distribution pattern of the advancing front of plasmin $m(\cdot, \tau_f)$ in the peritumoural region. As detailed in Appendix \ref{sec:The_MovingBoundaryMethod}, these movement characteristics are obtained by accounting the contribution of all peaks (baricentred at the spatial points $y_{l}$) at the front of advancing plasmin that are above the mean value of the entire mass of plasmin produced on $\epsilon Y \backslash  \Omega(t_0)$ and are located at the furthest away Euclidean distance from $\{x_{\epsilon Y}^* \}$. Thus, under these conditions, the moving direction $\eta_{\epsilon Y}$ and displacement magnitude $\xi_{\epsilon Y} $ derived in Appendix \ref{sec:The_MovingBoundaryMethod} in (\ref{eq:direction})-(\ref{eq:magnitude}) have the following expressions:

\bequd
\eta_{\epsilon Y} = x_{\epsilon Y}^* + \nu \underset{l \in \mathcal{I}_{\delta}}{\sum} \bigg( \int_{\mathcal{D}_l} m(y,\tau_f) dy \bigg) (y_l - x_{\epsilon Y}^*), \nu \in [0,\infty],
\eequd
\bequd
\xi_{\epsilon Y} := \underset{l \in \mathcal{I}_{\delta}}{\sum}  \frac{\int_{\mathcal{D}_l} m(y,\tau_f) dy  }{\underset{l \in \mathcal{I}_{\delta}}{\sum} \int_{\mathcal{D}_l} m(y,\tau_f) dy  } \big | \overrightarrow{x_{\epsilon Y}^*y_l} \big |.
\eequd
Finally, the transitional probability $q^*$ defined in (\ref{eq:probability}) is a quantification of the amount of plasmin in $\epsilon Y \backslash \Omega(t_0)$  relative to the total amount of plasmin concentration in $\epsilon Y$ and characterizes the invading strength. Therefore, the point $x^*_{\epsilon Y}$ will exercise the movement into the new spatial position $\widetilde{x^*_{\epsilon Y}}$ if and only if $q(x^*_{\epsilon Y}):=q^*(\epsilon Y)$ exceeds a certain tissue local threshold $\omega (\beta, \epsilon Y) \in (0,1)$ associated with the micro-domain $\epsilon Y$ under a given state of favourable conditions $\beta$.

\section{Multiscale Computational Simulation Results}
\label{sec:simulation}

The multiscale model of cancer invasion that we proposed here was numerically solved in a rectangular region $Y:= [0,4 \times [0,4]$. For all our simulations, we discretise the entire cube $Y$ uniformly, using the spatial mesh size: $\Delta x=\Delta y=\frac{\epsilon}{2} = 0.03125$. In this context, the initial conditions for the macro-dynamics are specified as follows. The cancer cell population on the initially considered tumour region $\Omega(0):=\Bila((2,2),0.5)$ is assumed to be given by a translated Gaussian centred at $(2,2)$ and mollified to $0$ after a radius of $0.5$, namely:
\bequ\label{inicond_c_2016}
c(x,0) = \frac{\bigg(\!\!\exp\big(-\frac{||x-(2,2)||_2^2}{\sqrt{\Delta x \Delta y}}\big) - \exp(-28.125)\!\!\bigg)\big(\chi_{_{\Bila((2,2),0.5-\gamma)}}\ast \psi_{_{\gamma}}\big) }{2},\quad x \in Y, 
\eequ
where $\psi_{\gamma}$ is the mollifier defined in \eqref{mollifier_def1}-\eqref{mollifier_def2} with $\gamma<<\frac{\Delta x}{3}$.
Further, the initial macroscopic conditions for the enzymatic components entering the uPA system are considered as follows:
\bequ
\begin{array}{llr}
u(x,0) & = 1-\frac{1}{2}c(x,0),& \qquad x\in Y\\[0.2cm]
p(x,0) & = \frac{1}{2}c(x,0), & x\in Y \\[0.2cm]
m(x,0) & = \frac{1}{20}c(x,0), & x\in Y 
\end{array}
\eequ

\begin{figure}
\vspace{-15pt}
\hspace{-5em}
\begin{tabular}{cc}
\subfloat{\label{square_n}\hspace{5pt}\includegraphics[scale = 0.375]{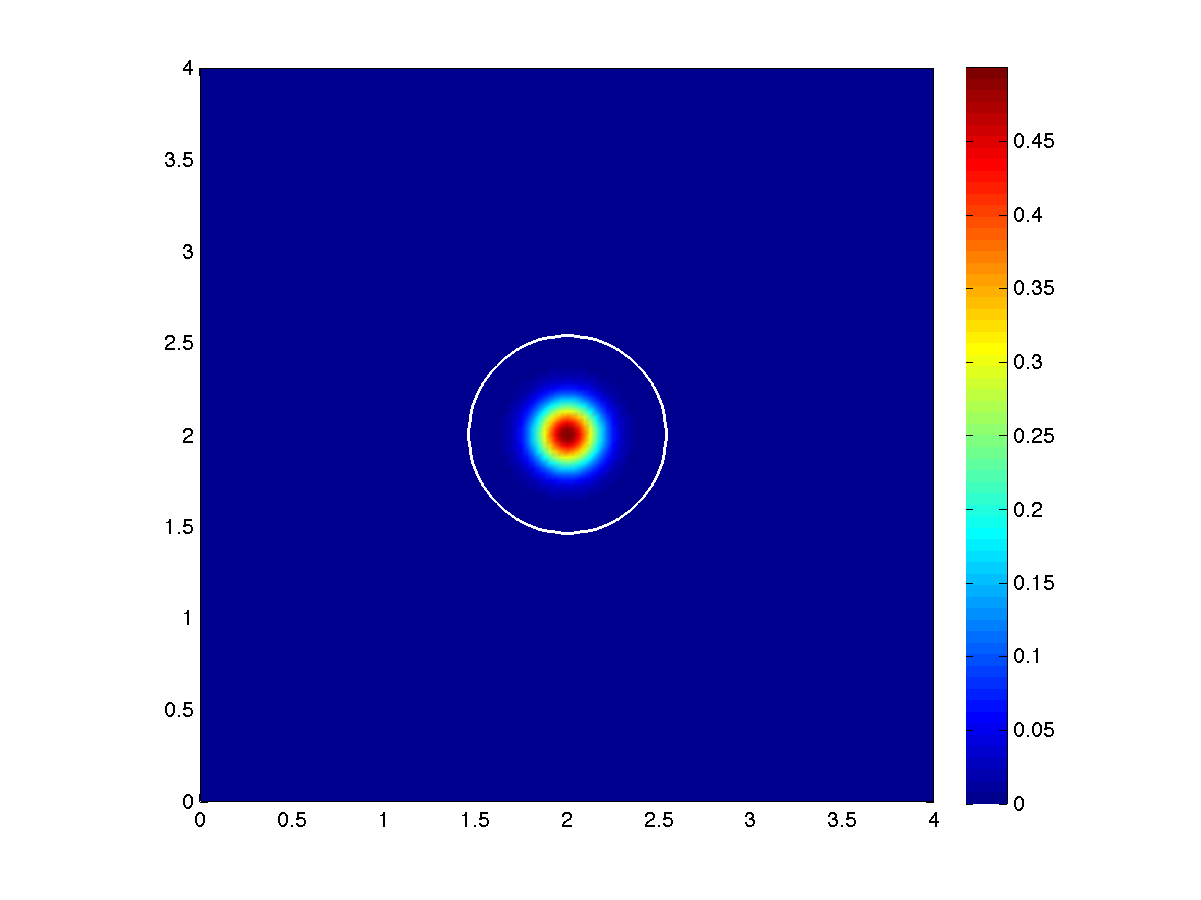}} \hspace{-35pt}
\subfloat{\label{square_v} \includegraphics[scale = 0.375]{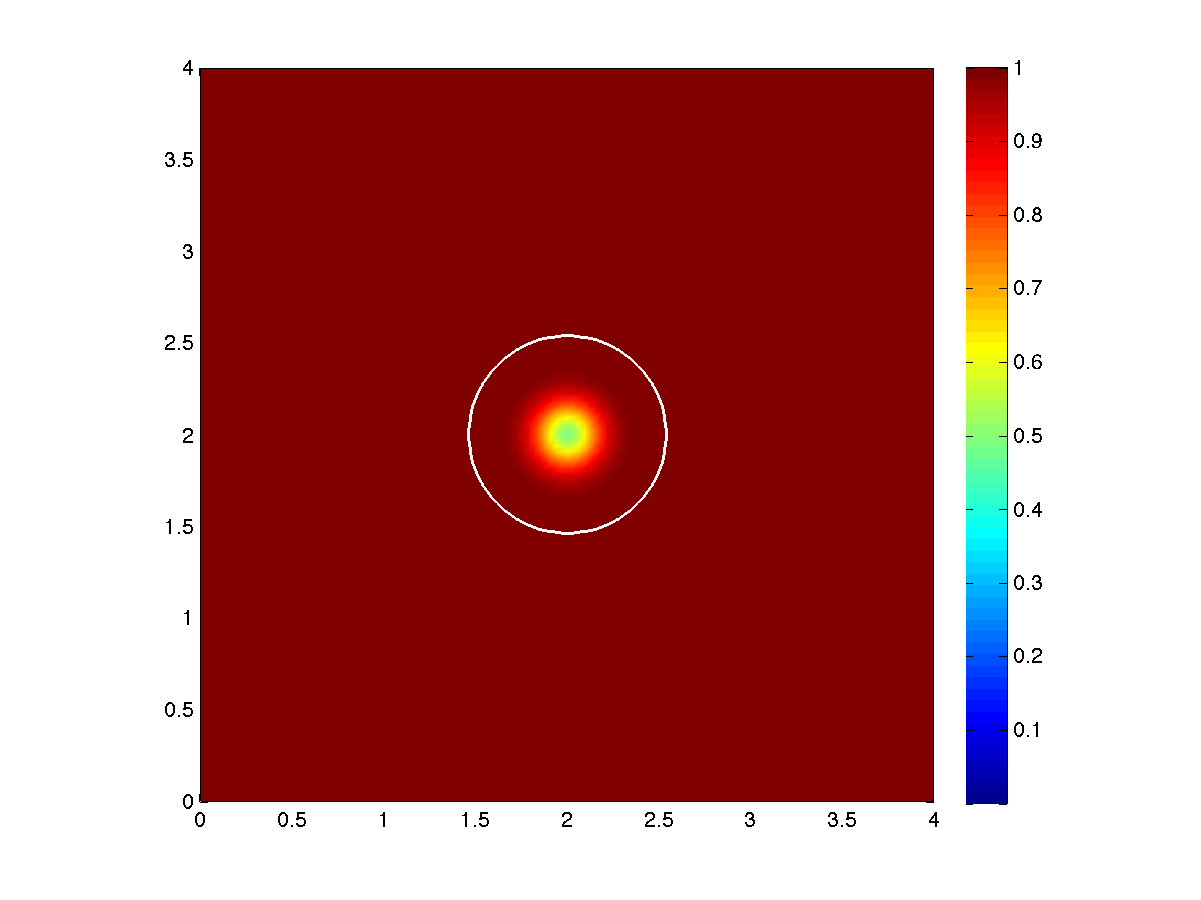}}\\[-20pt]
\subfloat{\label{square_n}\hspace{5pt}\includegraphics[scale = 0.375]{Cancer0.png}}  \hspace{-35pt}
\subfloat{\label{square_v} \includegraphics[scale = 0.375]{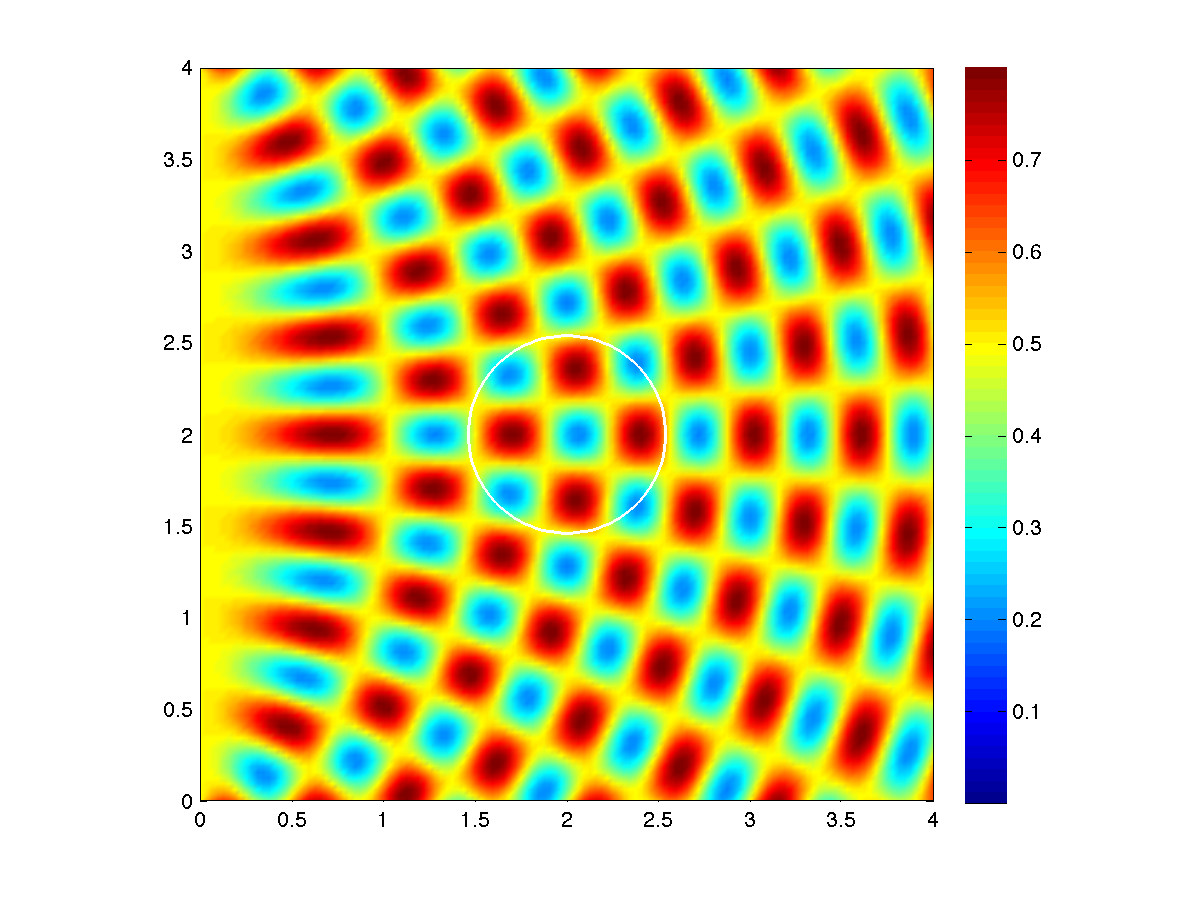}}
\end{tabular}
\vspace{-15pt}
\caption{\revdt{Initial conditions of the distributions of cancer cells (left column) and ECM (right column) and the invasive boundary of the tumour (white line) for the homogeneous case (top row) and the heterogeneous case(bottom row).}}
\label{fig:inicond2016}
\end{figure}
\noindent At the same time, in all the subsequent simulations, we consider the following homogeneous and heterogeneous initial condition for ECM, namely:
\bequ\label{inicond_v_2016}
\begin{array}{lll}
\textrm{homogeneous case:}\quad& v(x,0) = 1 - c(x,0),  &\quad x\in Y, \\[0.2cm]
\textrm{heterogeneous case:}\quad& v(x,0)  = \frac{1+ 0.3\sin{(4\pi||x ||_2)} + \sin{(4\pi || (4,0)-x ||_2)}}{2},  &\quad x\in Y.
\end{array}
\eequ
\noindent\revdt{Figure~\ref{fig:inicond2016} shows the initial conditions (\ref{inicond_c_2016}) and (\ref{inicond_v_2016}) of ECM and cancer cell distributions for both homogeneous and heterogeneous cases that are used in all the simulations presented in this paper.} 

Finally, in the absence of medical data,  for the tissue threshold $\omega (\beta, \epsilon Y)$ controlling whether or not a point on the boundary will exercised the movement according to the direction and displacement magnitude defined in Appendix \ref{sec:The_MovingBoundaryMethod}, we adopt the following functional form:
\begin{equation}
 \omega (\beta, \epsilon Y)\!\! :=\!\! \left\{ 
 \!\! \begin{array}{l l}
   \textrm{sin}\!\bigg(\!\!\frac{\pi}{2} \!\bigg(\!1\!-\!\frac{1}{\beta} \frac{v_{\omega(t_0)}(x_{\epsilon Y}^*, t_0+\Delta t)}{\underset{\xi \in \partial\Omega(t_0)}{\mathrm{sup}}  \!\!v_{\Omega(t_0)}(\xi,t_0 +\Delta t)}\!\bigg)\!\!\!\bigg) & \quad \textrm{if }  \frac{v_{\omega(t_0)}(x_{\epsilon Y}^*, t_0+\Delta t)}{\underset{\xi \in \partial\Omega(t_0)}{\mathrm{sup}} \!\!v_{\Omega(t_0)}(\xi,t_0 +\Delta t)} \leq\! \beta \\
   \textrm{sin}\!\bigg(\!\!\frac{\pi}{2(1-\beta)}\! \bigg(\!\frac{v_{\omega(t_0)}(x_{\epsilon Y}^*, t_0+\Delta t)}{\underset{\xi \in \partial\Omega(t_0)}{\mathrm{sup}}  \!\!v_{\Omega(t_0)}(\xi,t_0 +\Delta t)} \! -\! \beta \!\bigg)\!\!\!\bigg)  & \quad \textrm{if }  \frac{v_{\omega(t_0)}(x_{\epsilon Y}^*, t_0+\Delta t)}{\underset{\xi \in \partial\Omega(t_0)}{\mathrm{sup}} \!\!v_{\Omega(t_0)}(\xi,t_0 +\Delta t)} > \!\beta
  \end{array} \right.
     \label{eq:threshold}
  \end{equation}
 where $\beta \in(0,1)$ is a parameter that controls a certain ``optimal level" of ECM degradation and consider this as being the indicator of the most favourable invasion conditions at the level of tumour and tissue microenvironment. \revdt{This functional form of $\omega$ rules out any invasion if either a complete destruction or a very superficial degradation of the surrounding ECM is performed by the MDEs. This is due to the fact that while for invasion the cancer cells need considerable level of the ECM degradation, a total destruction of the surrounding ECM structure prevents them to advance further into the tissue, as they need some of the ECM components to adhere to in order to migrate. Finally, while this functional form of the tissue threshold aimed to showcase the proposed multiscale modelling framework, future work seeks to infer  $\omega (\cdot, \epsilon Y)$ from medical imaging data of the peritumoural tissue.}

The following figures show the simulation results of the evolving cancer cell and ECM spatial distributions and of the invasive tumour boundary at macro-time stage 20, 40, 60. The images are presented in two columns, with the left columns representing the cancer cell distribution and right columns showing the corresponding ECM concentrations. Furthermore, all these images include the tumour boundary. 

In the subsequent simulations, for the macroscopic part of the model, except otherwise stated, we will generally be using the following basic set of parameter values \revdt{$\mathscr{P}:$ 
\begin{align}
D_n \!&= 4.3 \!\times \!10^{-3},  &  \chi_u \!&= 3.05\! \times \!10^{-2},  & \chi_p \!&= 3.75 \!\times \!10^{-2},  & \chi_v \!&= 2.85\! \times\! 10^{-2}, \notag \\ 
 \mu_1 \!&= 0.25, & \delta \!&= 8.15,  & \phi_{21} \!&=0.75, & \phi_{22} \!&=0.55, \notag\\ 
  \mu_2 \!&=0.15, &  D_u \!&= 2.5 \times 10^{-3} , & \phi_{31}  \!&=0.75, & \phi_{33} \!&=0.3, \notag \\
  \alpha_{31} \!&=0.215, & D_p \!&= 3.5 \times 10^{-3},  & \phi_{41} \!&=0.75, & \phi_{42} \!&=0.55, \notag \\
   \alpha_{41} \!&=0.5 , & D_m \!&= 4.91 \times 10^{-3} ,& \phi_{52} \!&=0.11, & \phi_{53} \!&= 0.75, \notag \\
     \phi_{54} \!&=0.5 .  
     \label{MAY06_2011_param}
\end{align}
which are detailed in Table 1 in Appendix \ref{param_app}}. However, in order to investigate different cancer growth patterns, as described in the following simulation, we will perform numerical tests also for slightly changed values for the diffusion coefficient of cancer cells $D_c$, the ECM proliferation rate $\mu_2$, or the ECM degradation rate $\delta$. Also, we will explore a range of values for parameter $\beta$, to highlight the correlation between changes in tissue microenvironment conditions conditions and the resulting cancer invasion patterns.

In order to analyse the effect of each variables (namely, ECM initial condition, cancer cell diffusion coefficient, ECM proliferation \& degradation rates, and threshold coefficient), we will split the results into four groups as follows:

\begin{figure}[htp]
\vspace{-15pt}
\hspace{-5em}
\begin{tabular}{cc}
\subfloat{\label{circle_n} \hspace{5pt}\includegraphics[scale = 0.375]{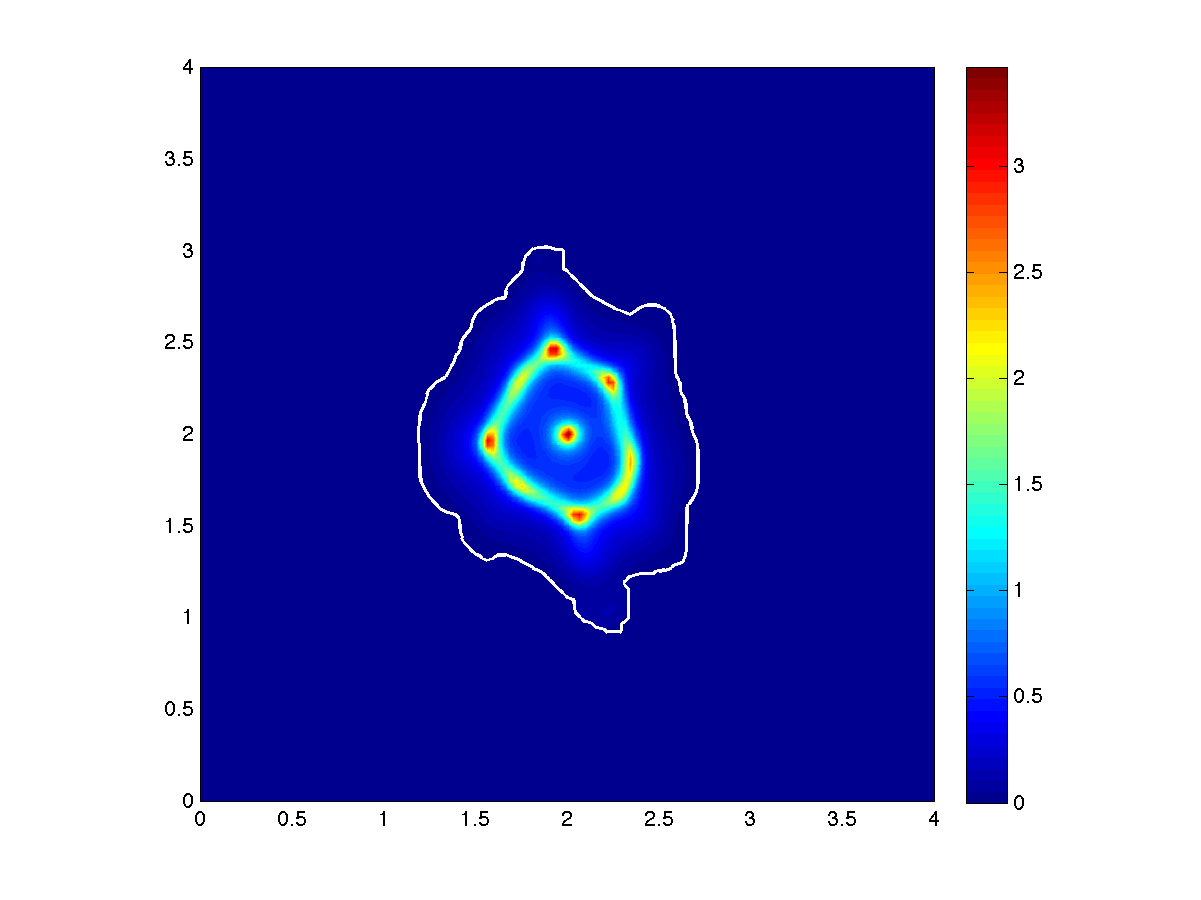}} \hspace{-35pt}
\subfloat{\label{circle_v} \includegraphics[scale = 0.375]{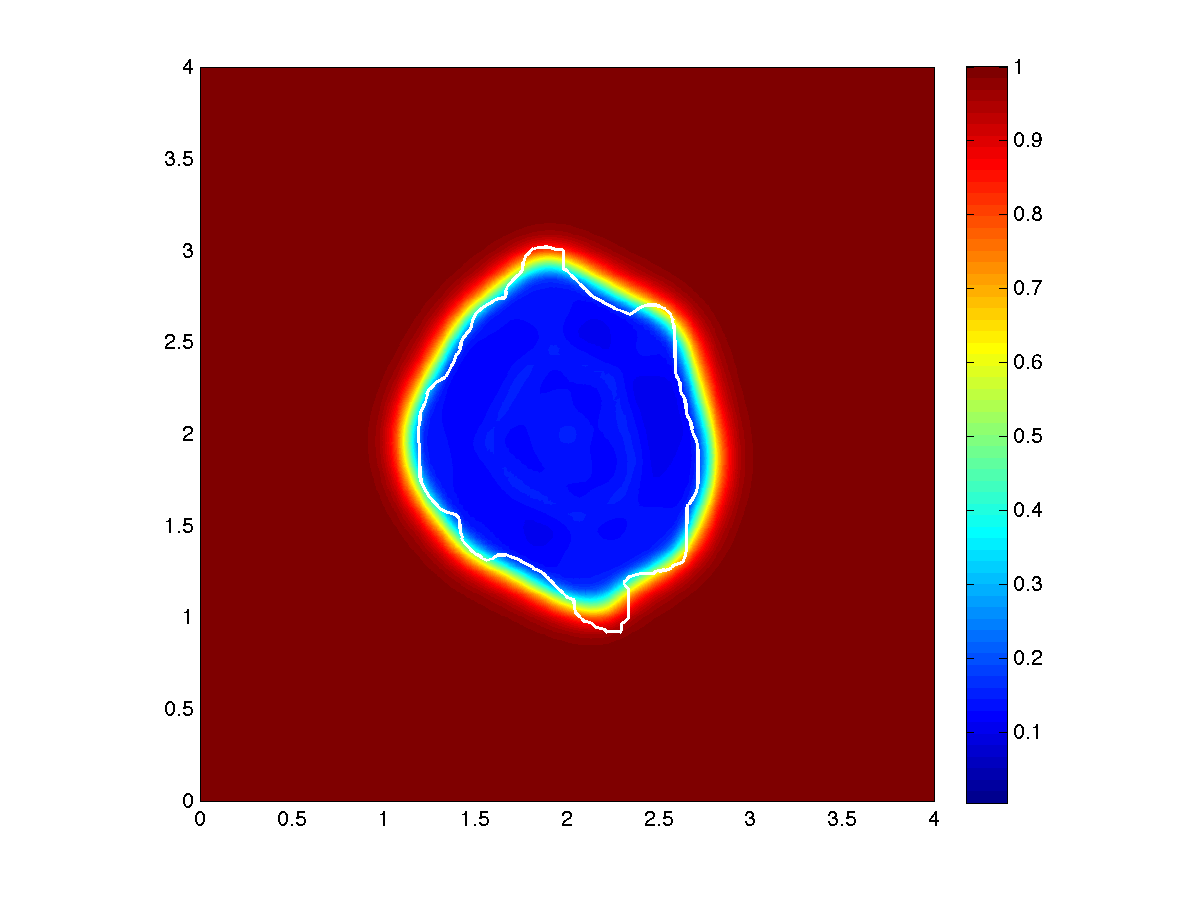}}\\[-20pt]
\subfloat{\label{circle_n} \hspace{5pt}\includegraphics[scale = 0.375]{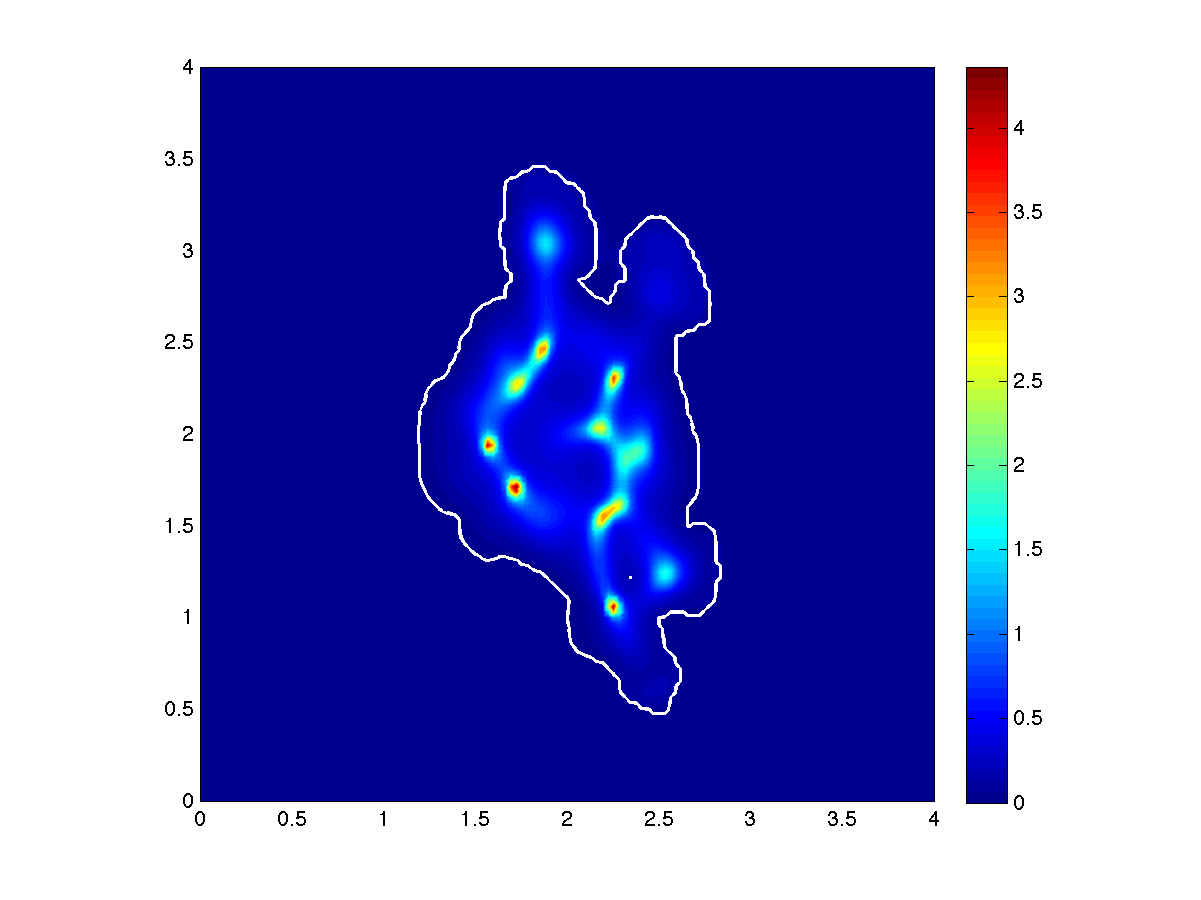}} \hspace{-35pt}
\subfloat{\label{circle_v} \includegraphics[scale = 0.375]{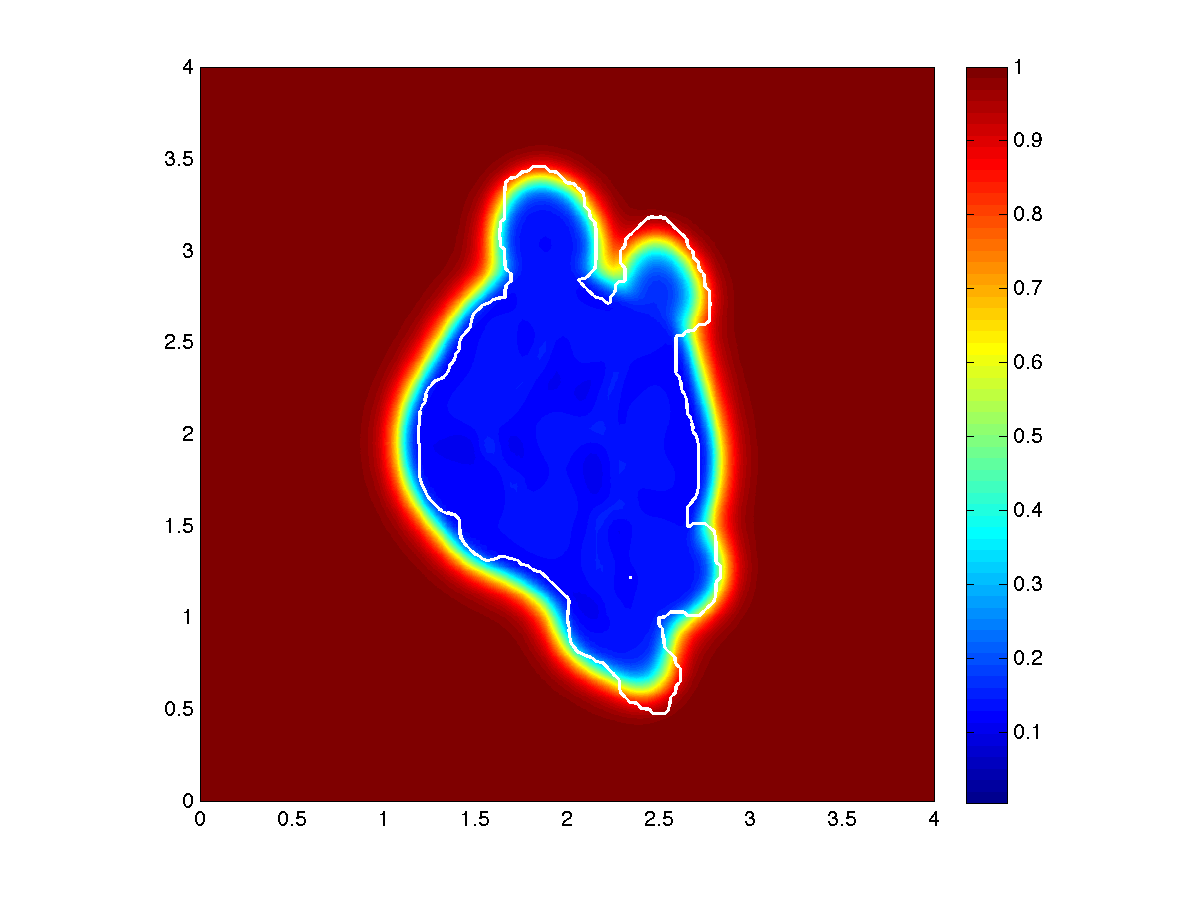}}\\[-20pt]
\subfloat{\label{circle_n} \hspace{5pt}\includegraphics[scale = 0.375]{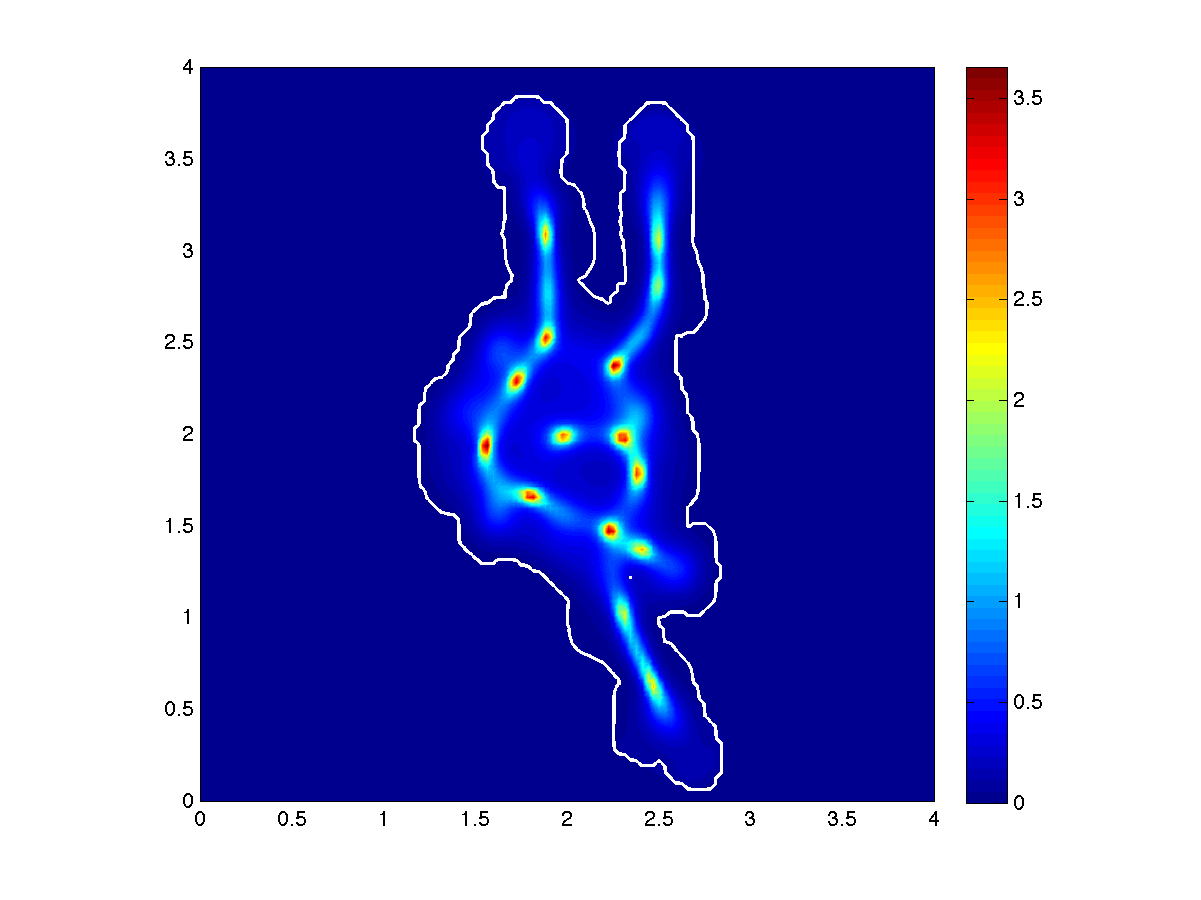}}  \hspace{-35pt}
\subfloat{\label{circle_v} \includegraphics[scale = 0.375]{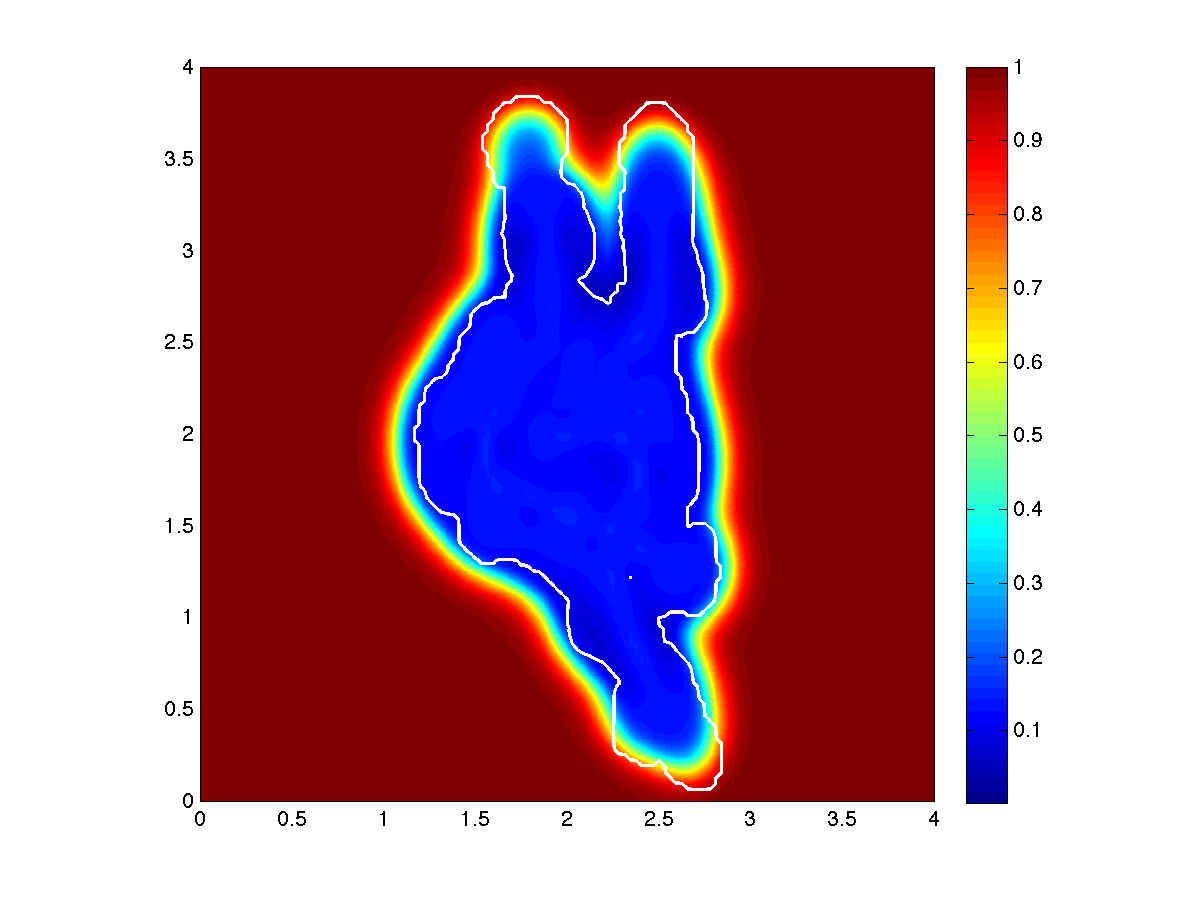}}
\end{tabular}
\vspace{-15pt}
\caption{Simulation results showing distributions of cancer cells (left column) and ECM (right column) and the invasive boundary of the tumour (white line) at various macro-micro stages: Stage 20, 40, 60. \revdt{Starting from the homogeneous initial conditions shown in Figure \ref{fig:inicond2016}, these results were obtained for $D_c~=~4.3 \times~10^{-3}$, $\beta = 0.775$, $\mu_2 = 0.01$ and $\delta = 1.5$.}}
\label{fig:homo}
\end{figure}


\begin{figure}[htp]
\vspace{-15pt}
\hspace{-5em}
\begin{tabular}{cc}
\subfloat{\label{circle_n} \hspace{5pt}\includegraphics[scale = 0.375]{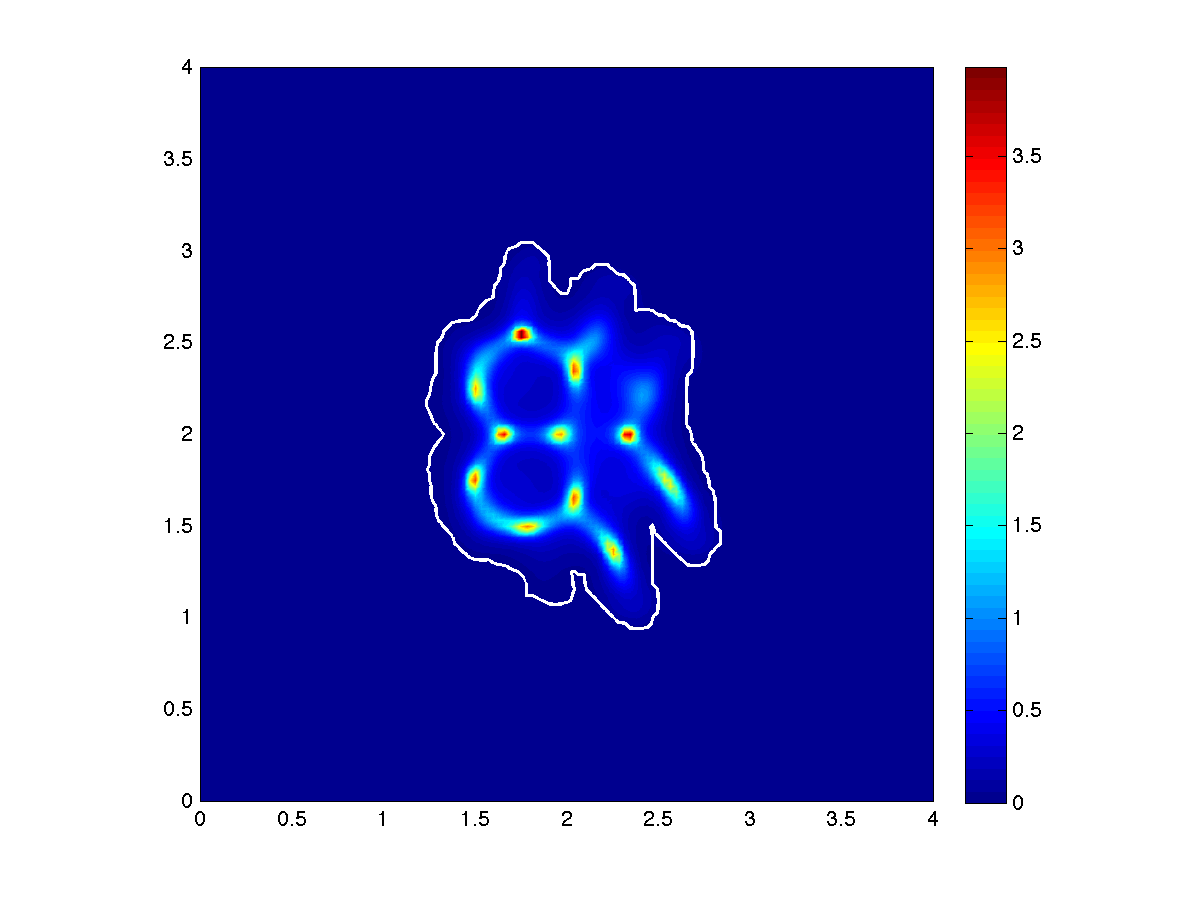}}  \hspace{-35pt}
\subfloat{\label{circle_v} \includegraphics[scale = 0.375]{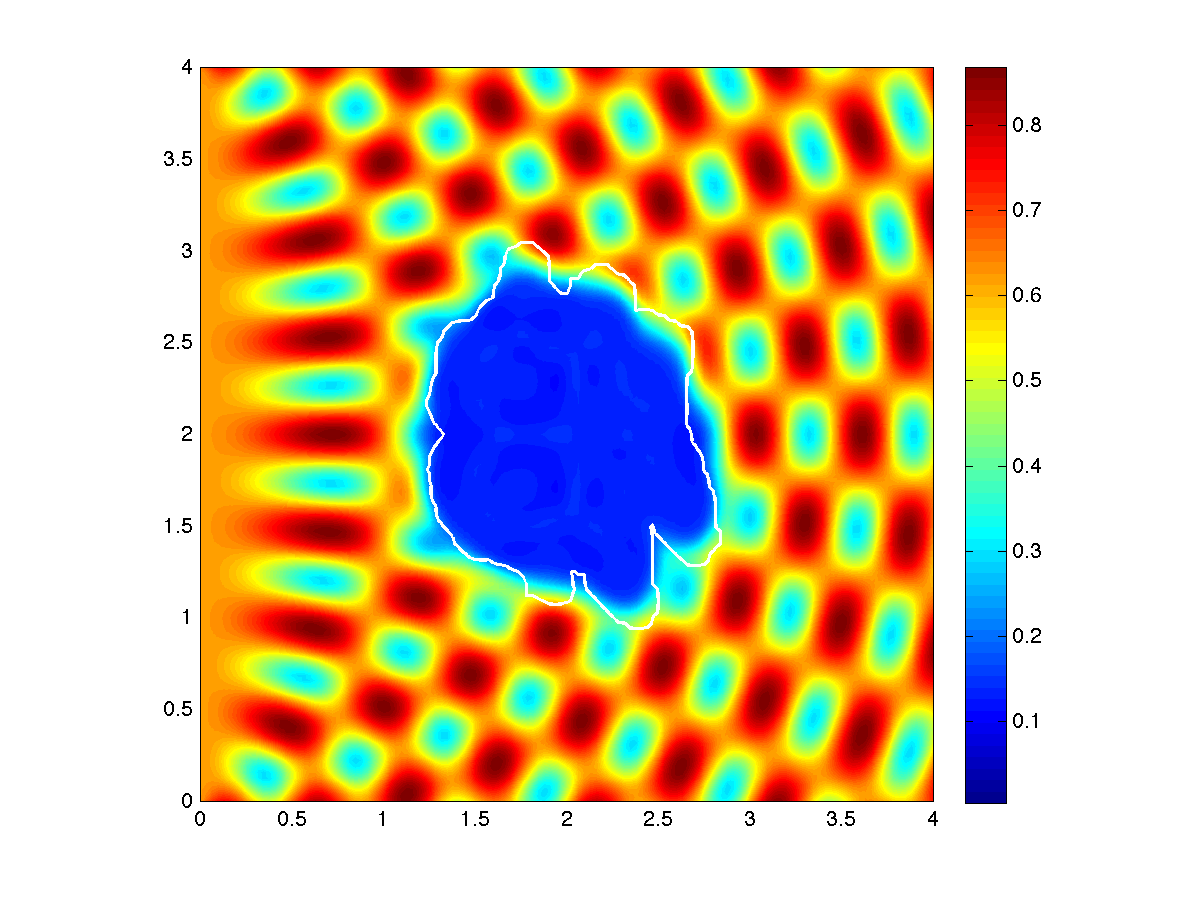}}\\[-20pt]
\subfloat{\label{circle_n} \hspace{5pt}\includegraphics[scale = 0.375]{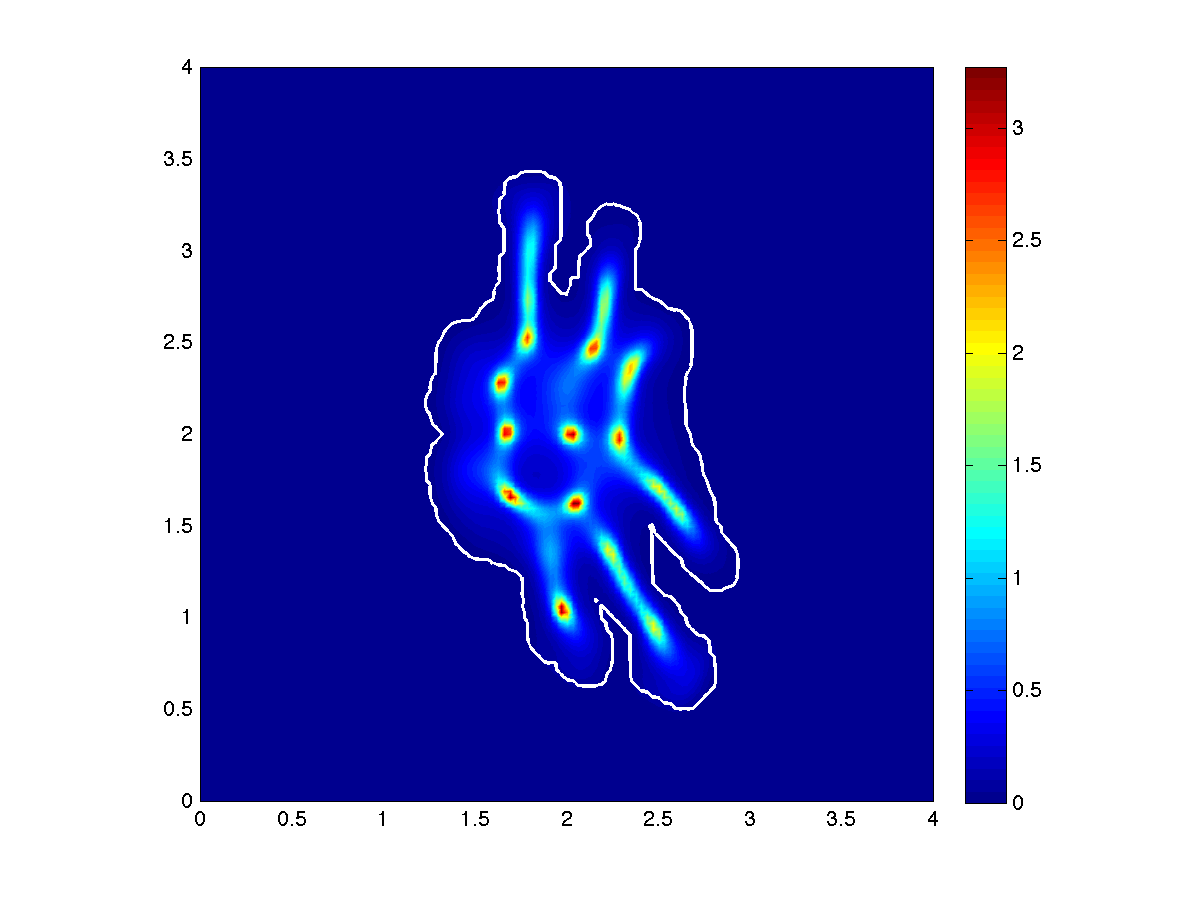}}  \hspace{-35pt}
\subfloat{\label{circle_v} \includegraphics[scale = 0.375]{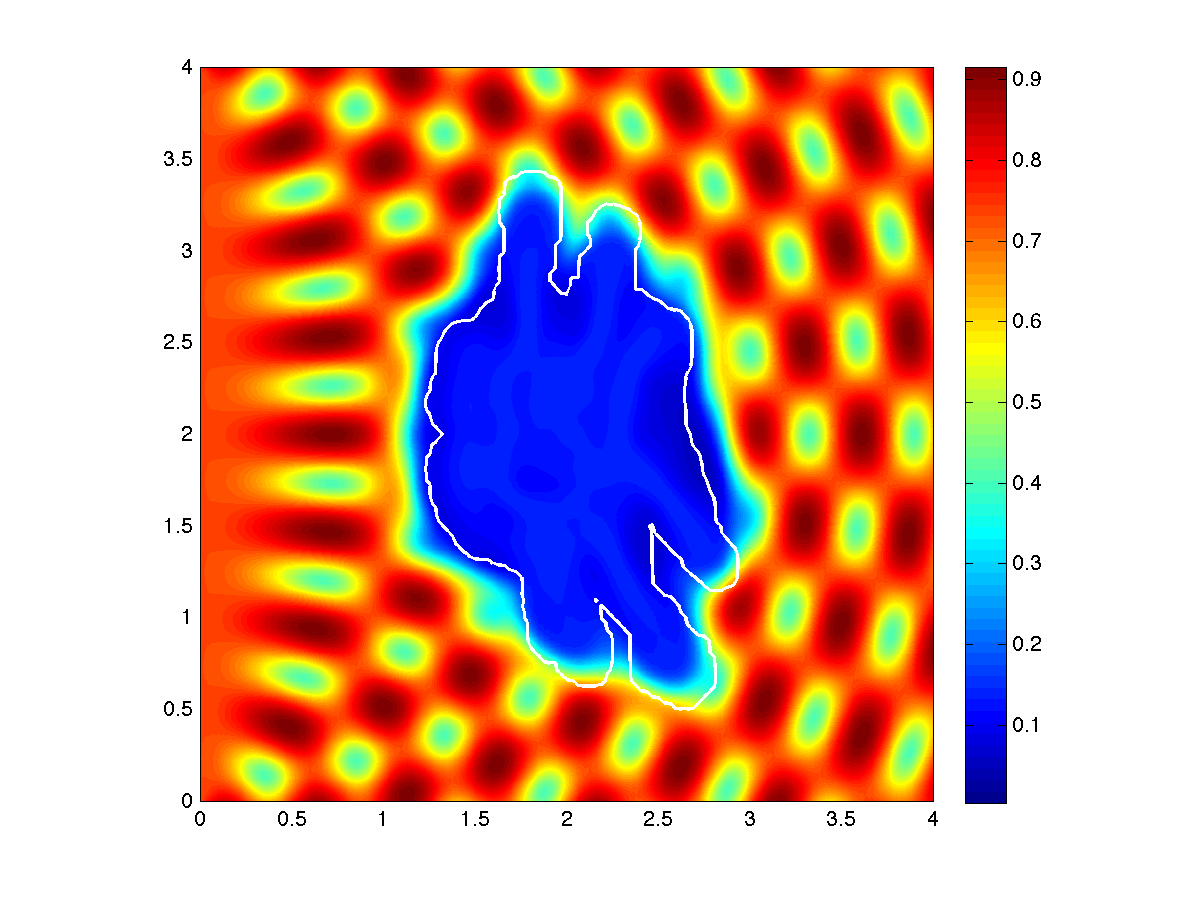}}\\[-20pt]
\subfloat{\label{circle_n} \hspace{5pt}\includegraphics[scale = 0.375]{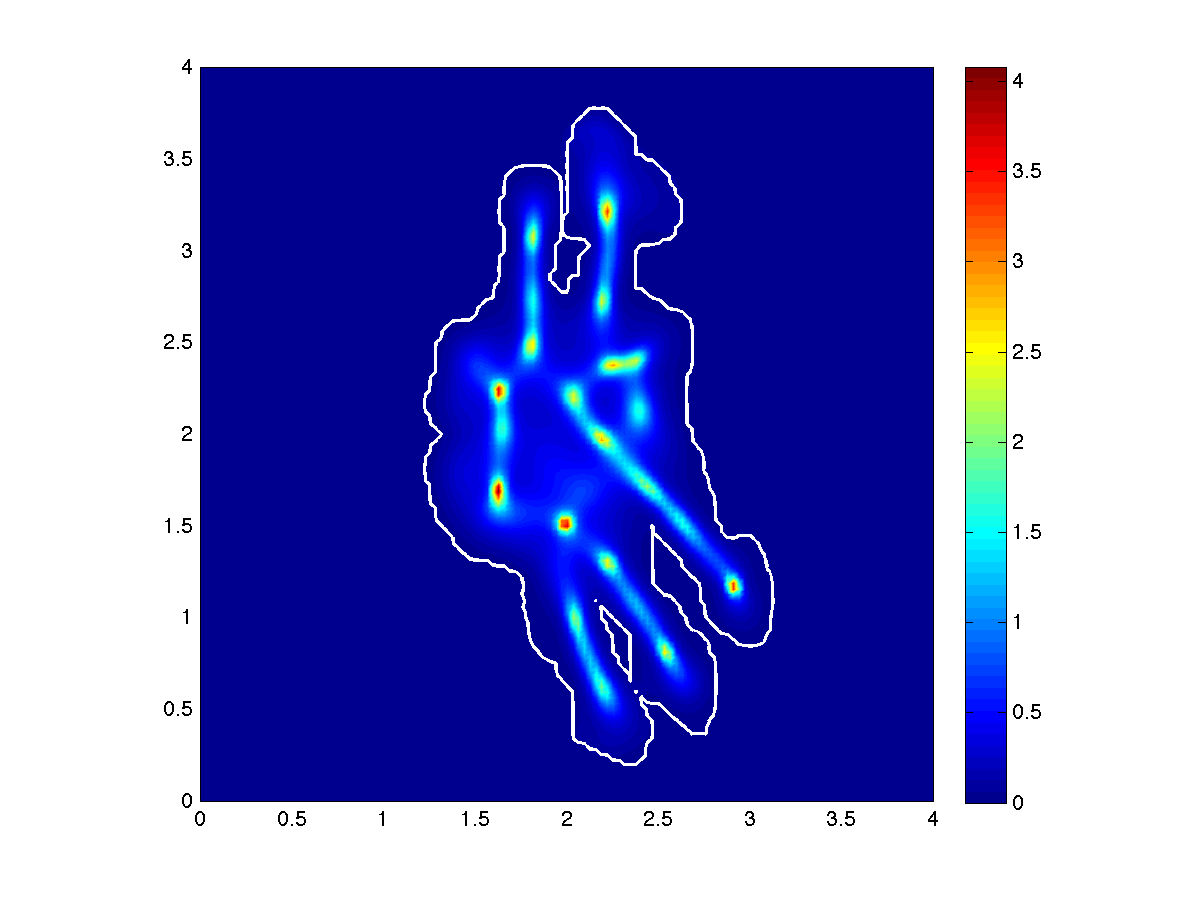}}   \hspace{-35pt}
\subfloat{\label{circle_v} \includegraphics[scale = 0.375]{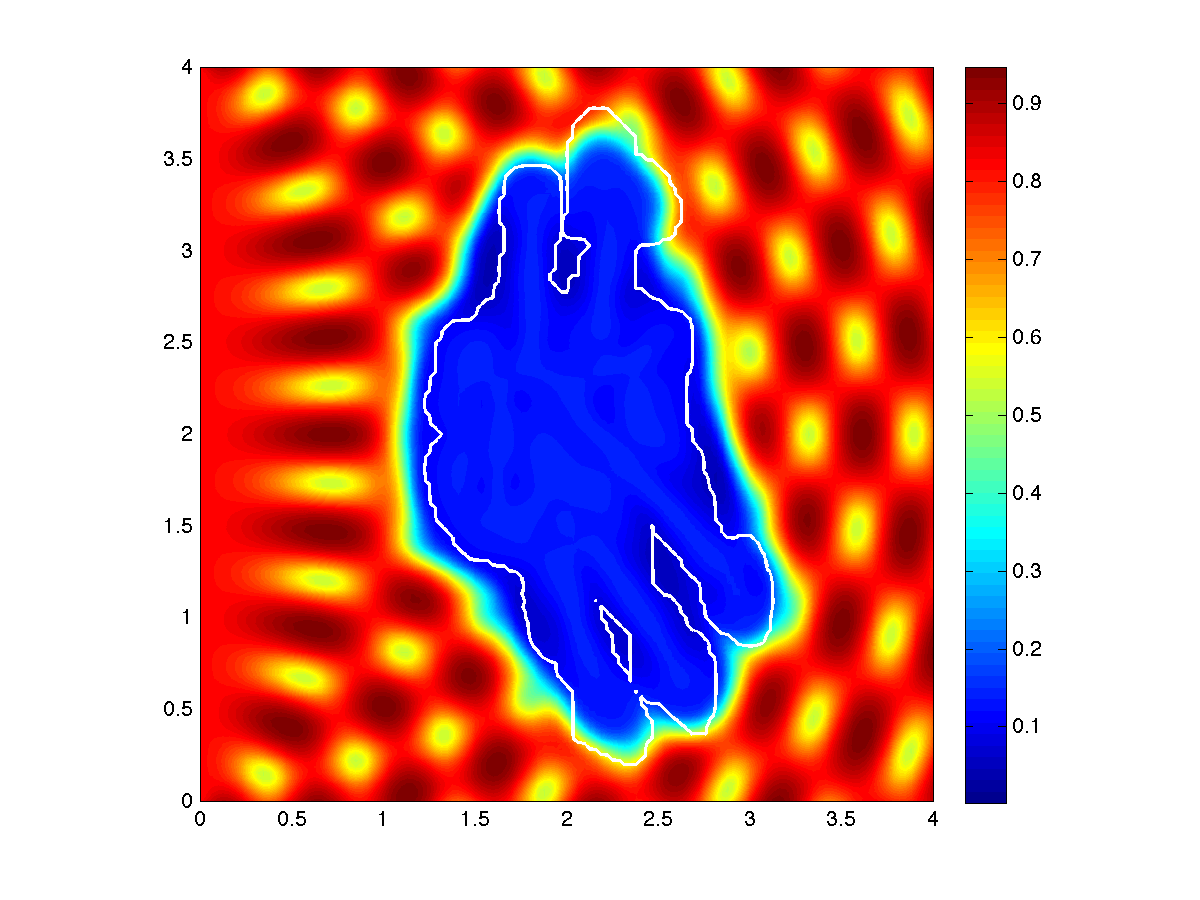}}
\end{tabular}
\vspace{-15pt}
\caption{Simulation results showing distributions of cancer cells (left column) and ECM (right column) and the invasive boundary of the tumour (white line) at various macro-micro stages: Stage 20, 40, 60. \revdt{Starting from the heterogeneous initial conditions shown in Figure \ref{fig:inicond2016}, these results were obtained for $D_c = 4.3 \times 10^{-3}$, $\beta = 0.775$, $\mu_2 = 0.01$ and $\delta = 1.5$.}}
\label{fig:hetero}
\end{figure}

\paragraph{\it{ECM initial condition.}}

To investigate what effect do different ECM initial conditions have on the whole dynamics of the model, we apply the same threshold function for both homogeneous and heterogeneous ECM scenario. 

\revdt{The ECM heterogeneity gives rise to a tumour-tissue interaction that is naturally more complex than in the case of homogeneous ECM. A direct consequence of this is that the ECM heterogeneity triggers a corresponding intrinsic variability in the tissue thresholds, which explore the peritumoural tissue conditions that the cancer interacts with during invasion. This results in a higher degree of complexity in the spatial structure of the regions with most favourable conditions for the tumour to progress further in the surrounding region, within the directions and displacement magnitudes specified by the microscale dynamics. Ultimately, this gives rise to a higher level of fingering and infiltrative patterns in the heterogeneous case as opposed to the homogeneous case. This is confirmed by our results presented in Figure ~\ref{fig:hetero} that are obtained for the heterogeneous ECM, which exhibit more ``fingered" and infiltrative spreading of the tumour compared with those obtained in the homogeneous ECM case shown in Figure~\ref{fig:homo}. This type of fingering patterns are often observed in medical imaging data, such as the one reported in the case of oesophageal and lung cancer by the authors in \cite{Ito_et_al_2012} and \cite{Masuda2012}, respectively.}

\paragraph{\it{Cancer cell diffusion coefficient $D_c$.}}

\begin{figure}[htp]
\vspace{-15pt}
\hspace{-5em}
\begin{tabular}{cc}
\subfloat{\label{circle_n} \hspace{5pt}\includegraphics[scale = 0.375]{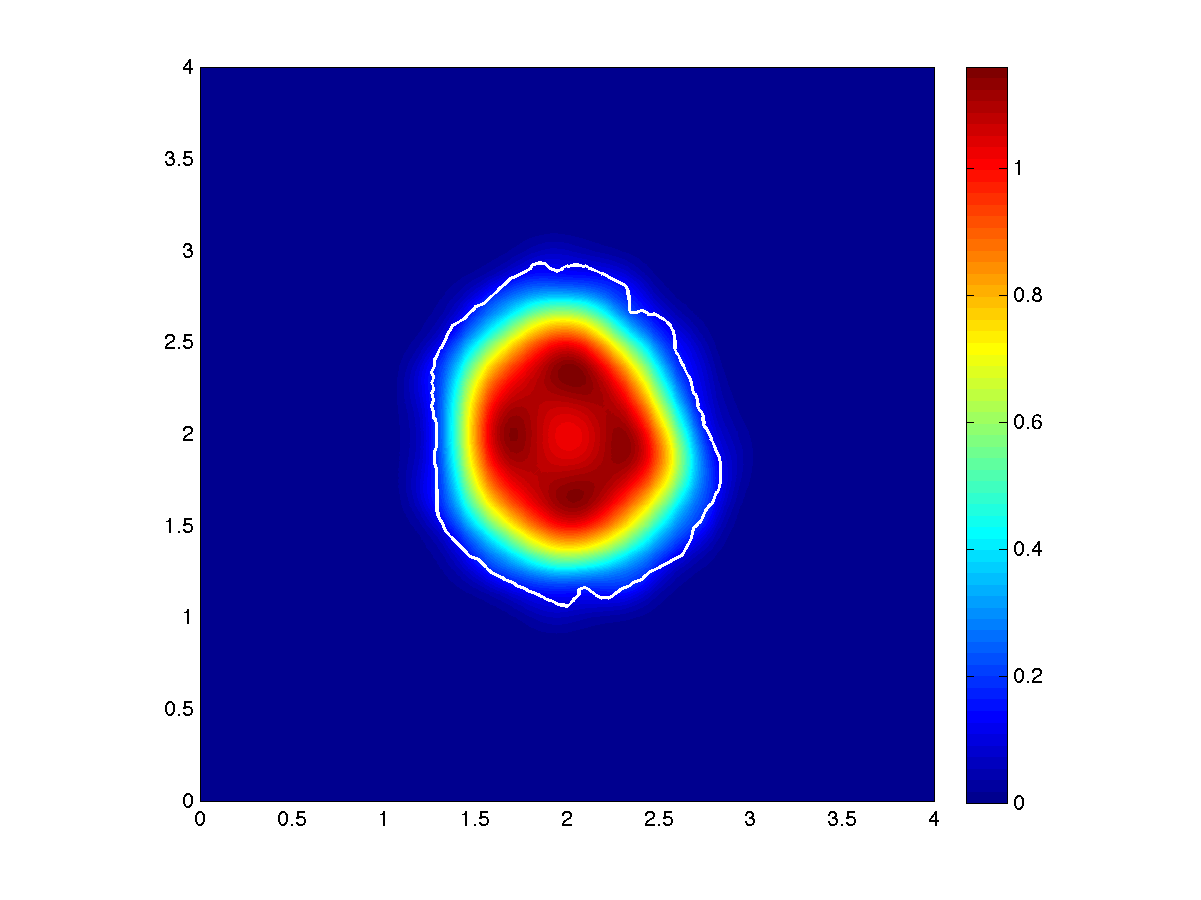}} \hspace{-35pt}
\subfloat{\label{circle_v} \includegraphics[scale = 0.375]{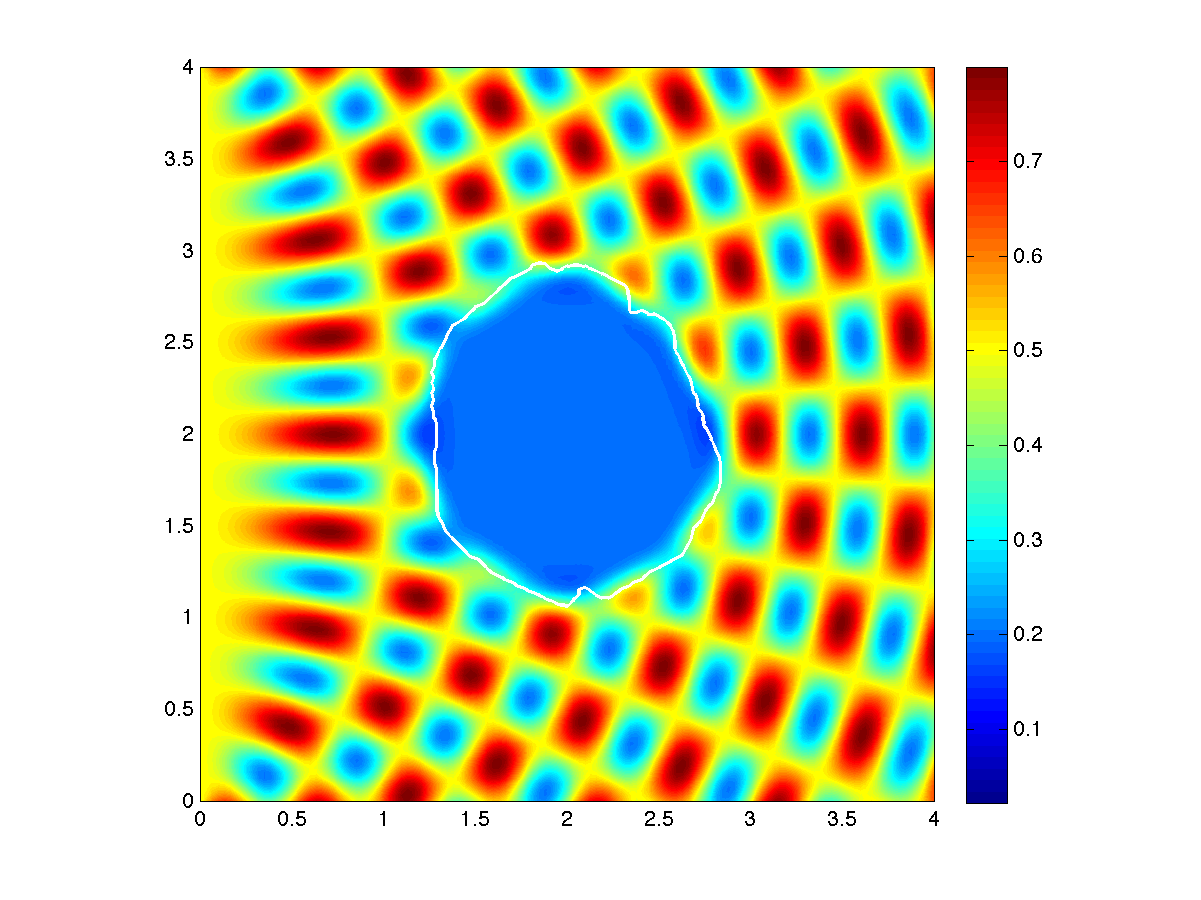}}\\[-20pt]
\subfloat{\label{circle_n} \hspace{5pt}\includegraphics[scale = 0.375]{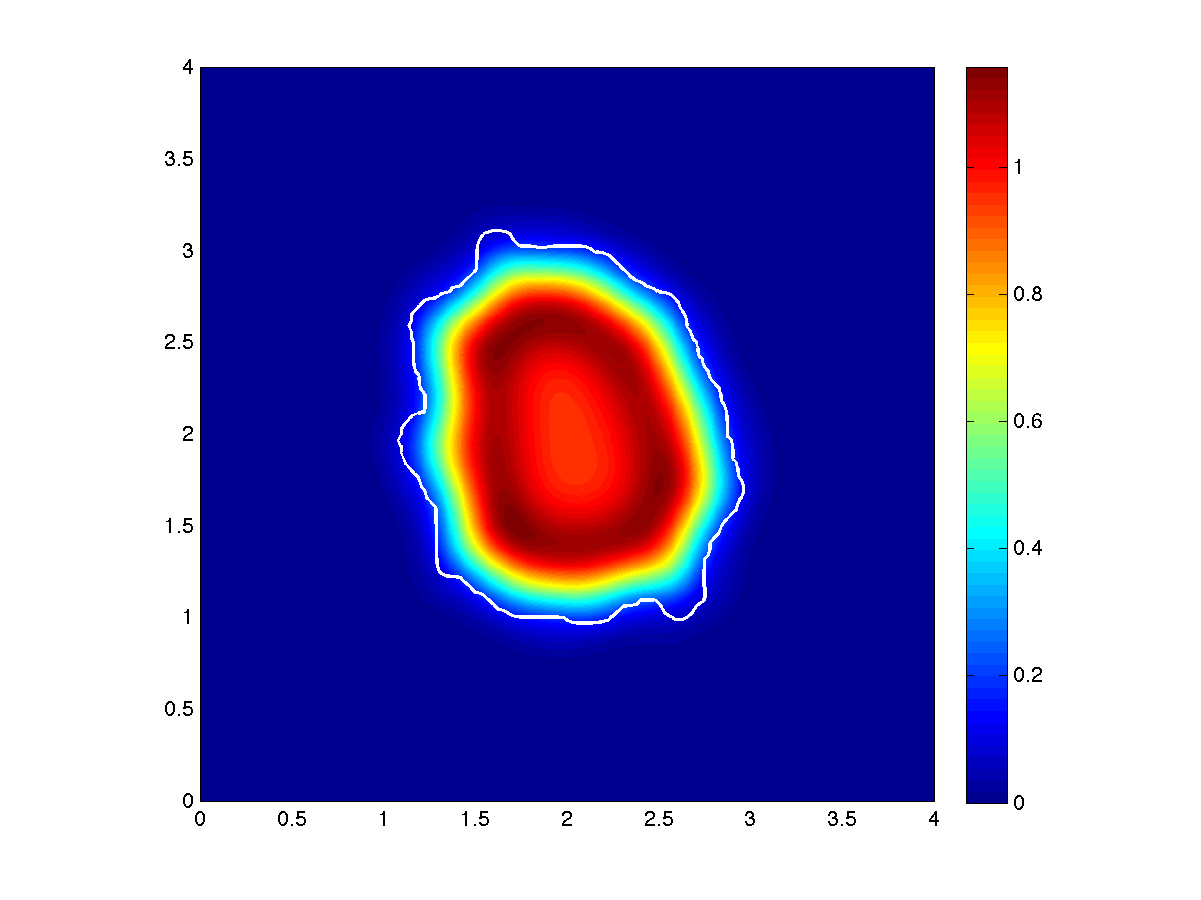}} \hspace{-35pt}
\subfloat{\label{circle_v} \includegraphics[scale = 0.375]{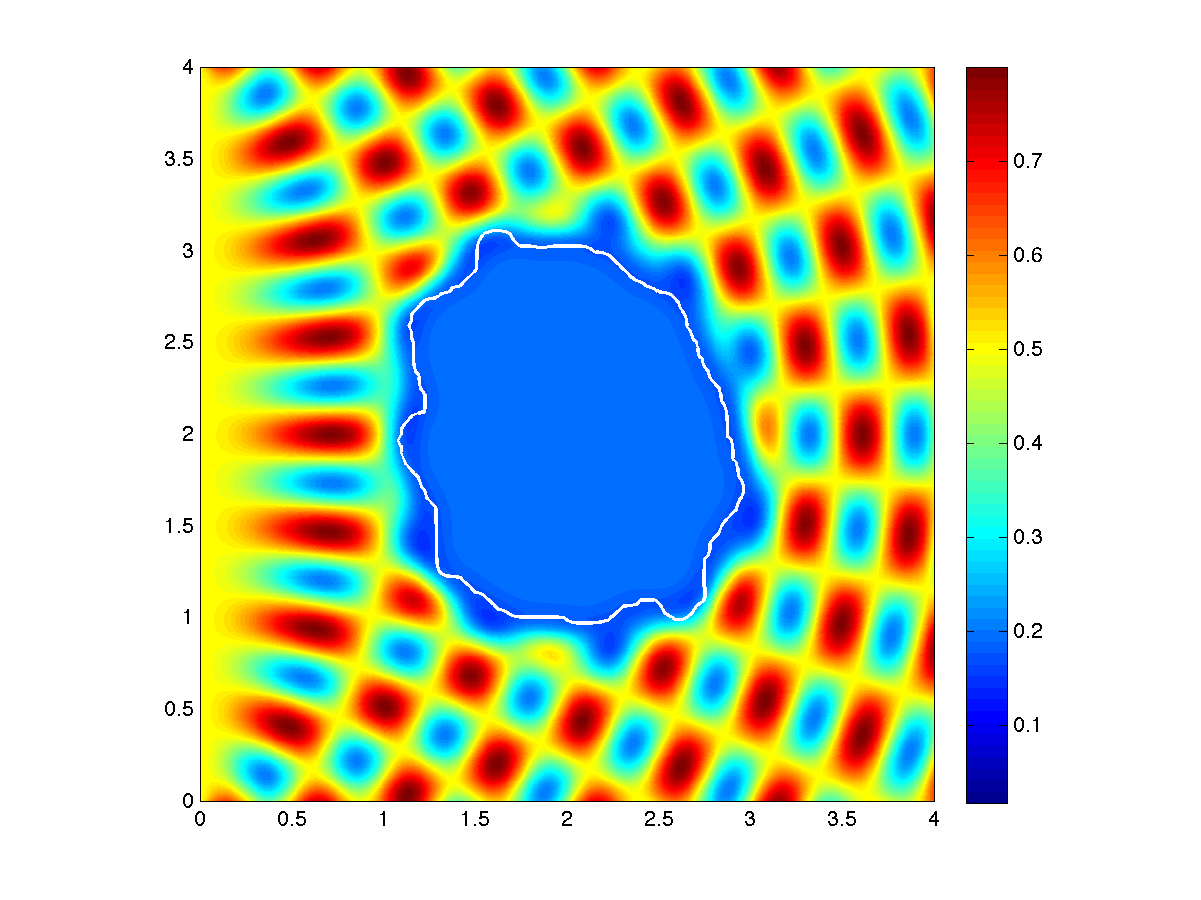}}\\[-20pt]
\subfloat{\label{circle_n} \hspace{5pt}\includegraphics[scale = 0.375]{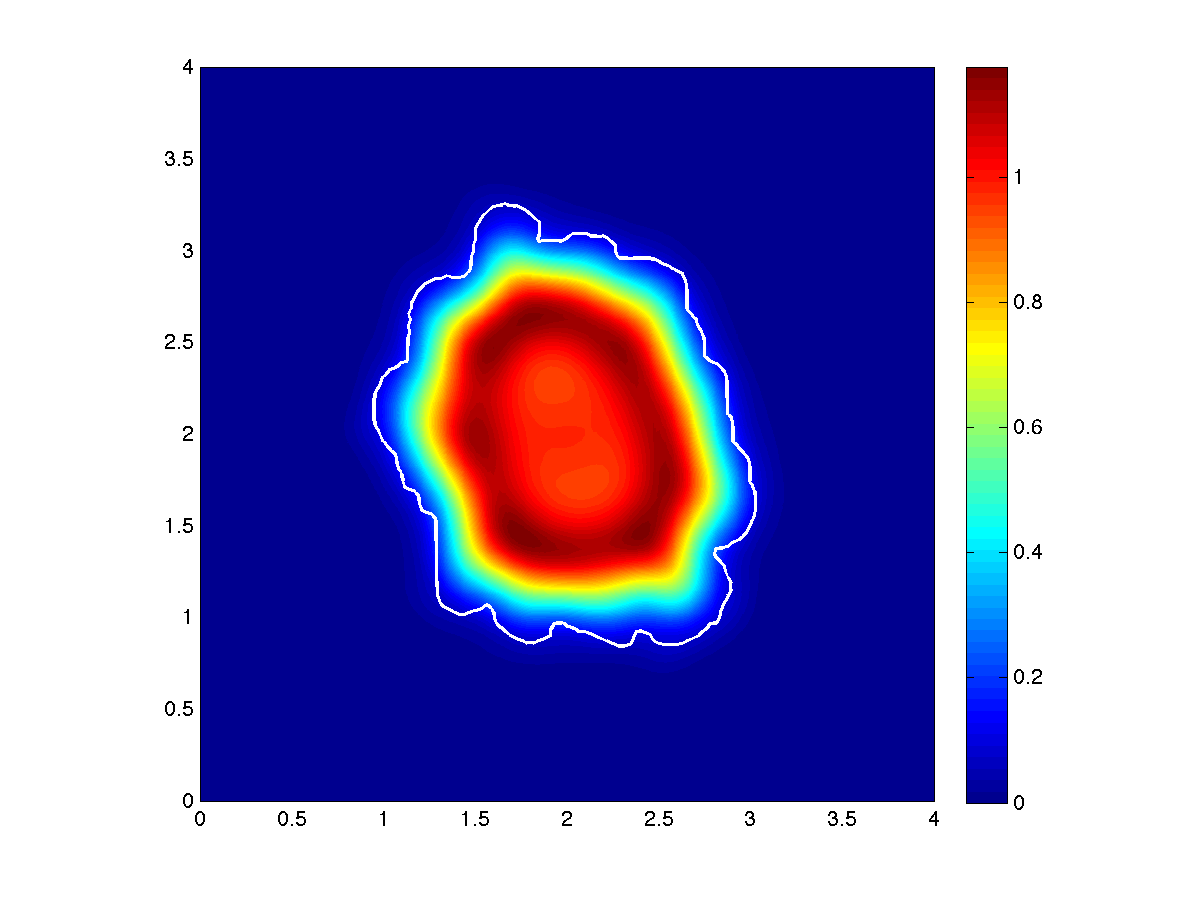}} \hspace{-35pt}
\subfloat{\label{circle_v} \includegraphics[scale = 0.375]{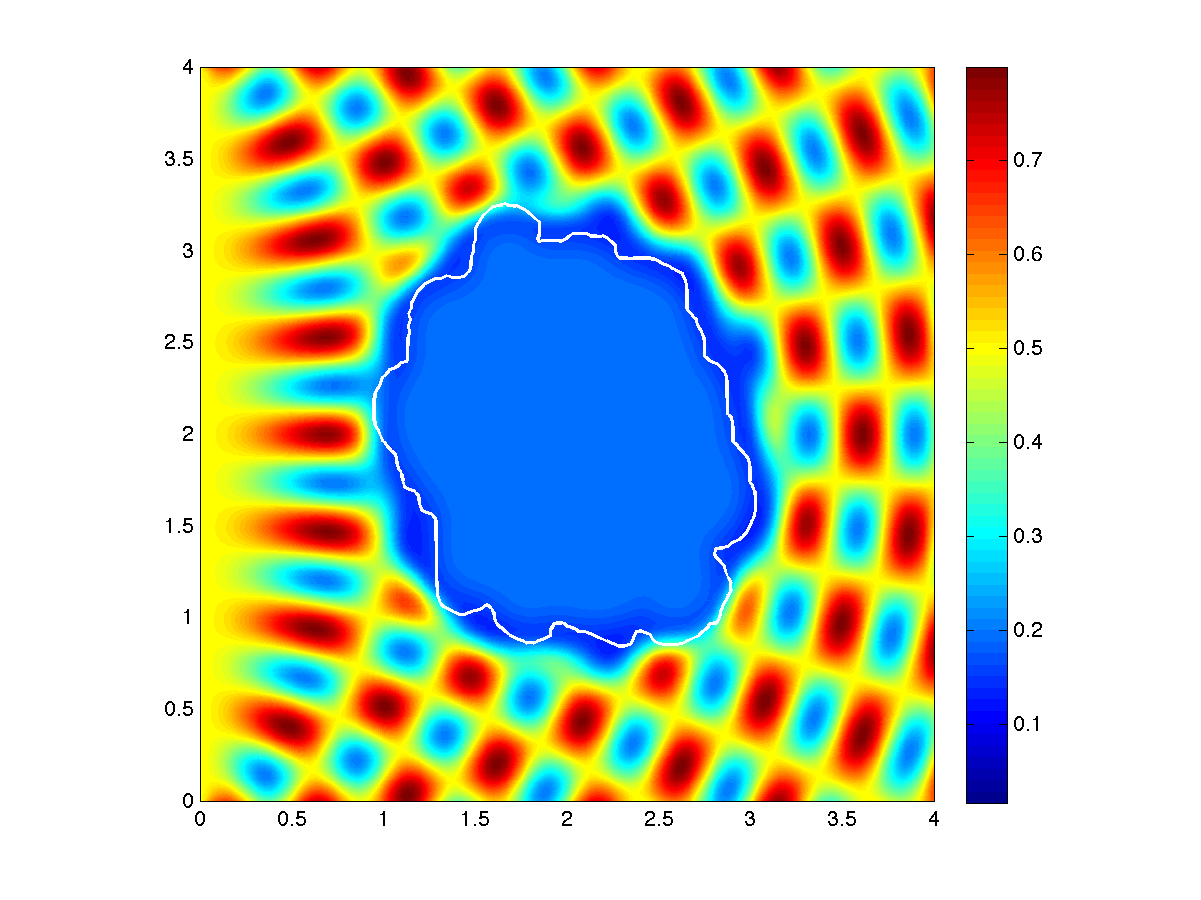}}
\end{tabular}
\vspace{-15pt}
\caption{Simulation results showing distributions of cancer cells (left column) and ECM (right column) and the invasive boundary of the tumour (white line) at various macro-micro stages: Stage 20, 40, 60. Starting from the heterogeneous initial conditions shown in Figure \ref{fig:inicond2016}, these results were obtained for $D_c = 1.4\times 10^{-2}$, $\beta = 0.775$, $\mu_2 = 0$ and $\delta = 0.75$.}
\label{fig:Dn1}
\end{figure}

\begin{figure}[htp]
\vspace{-15pt}
\hspace{-5em}
\begin{tabular}{cc}
\subfloat{\label{circle_n} \hspace{5pt}\includegraphics[scale = 0.375]{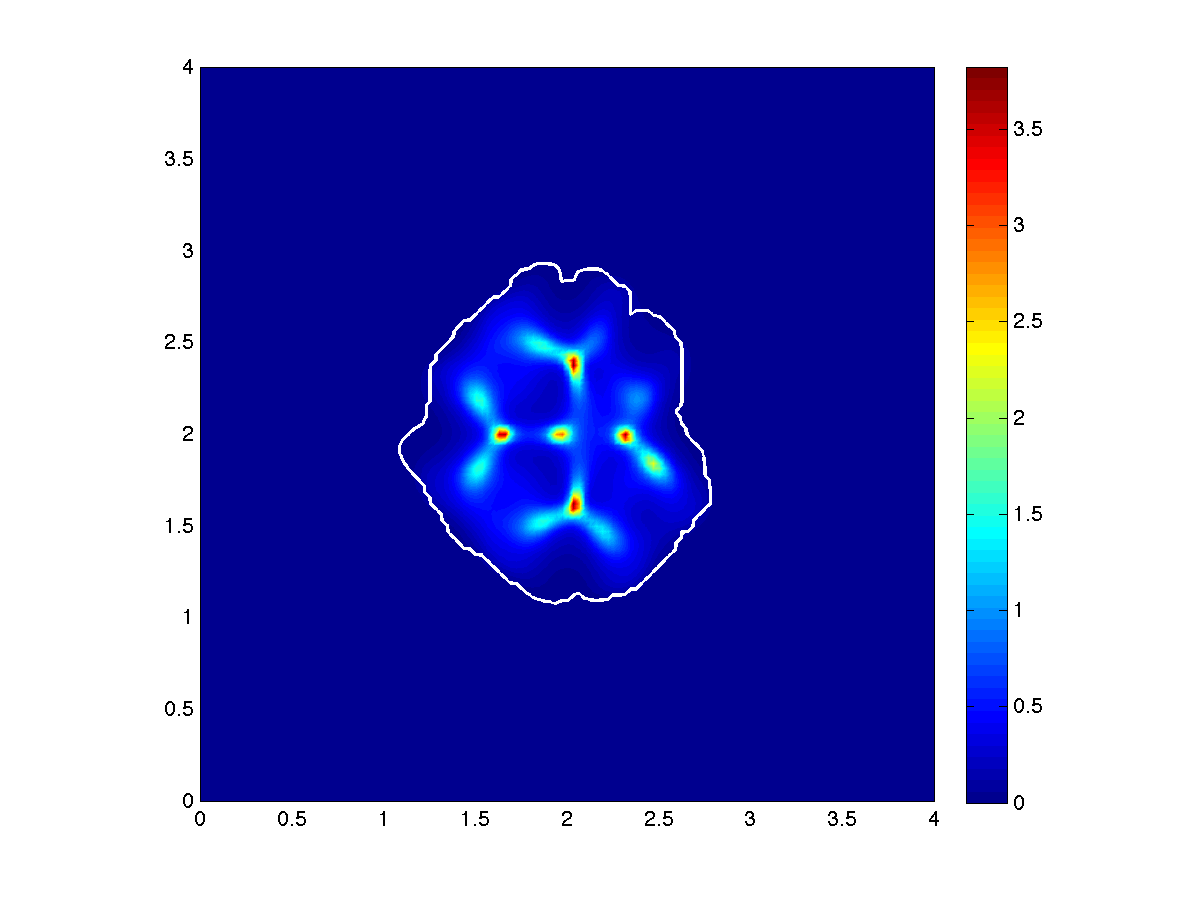}} \hspace{-35pt}
\subfloat{\label{circle_v} \includegraphics[scale = 0.375]{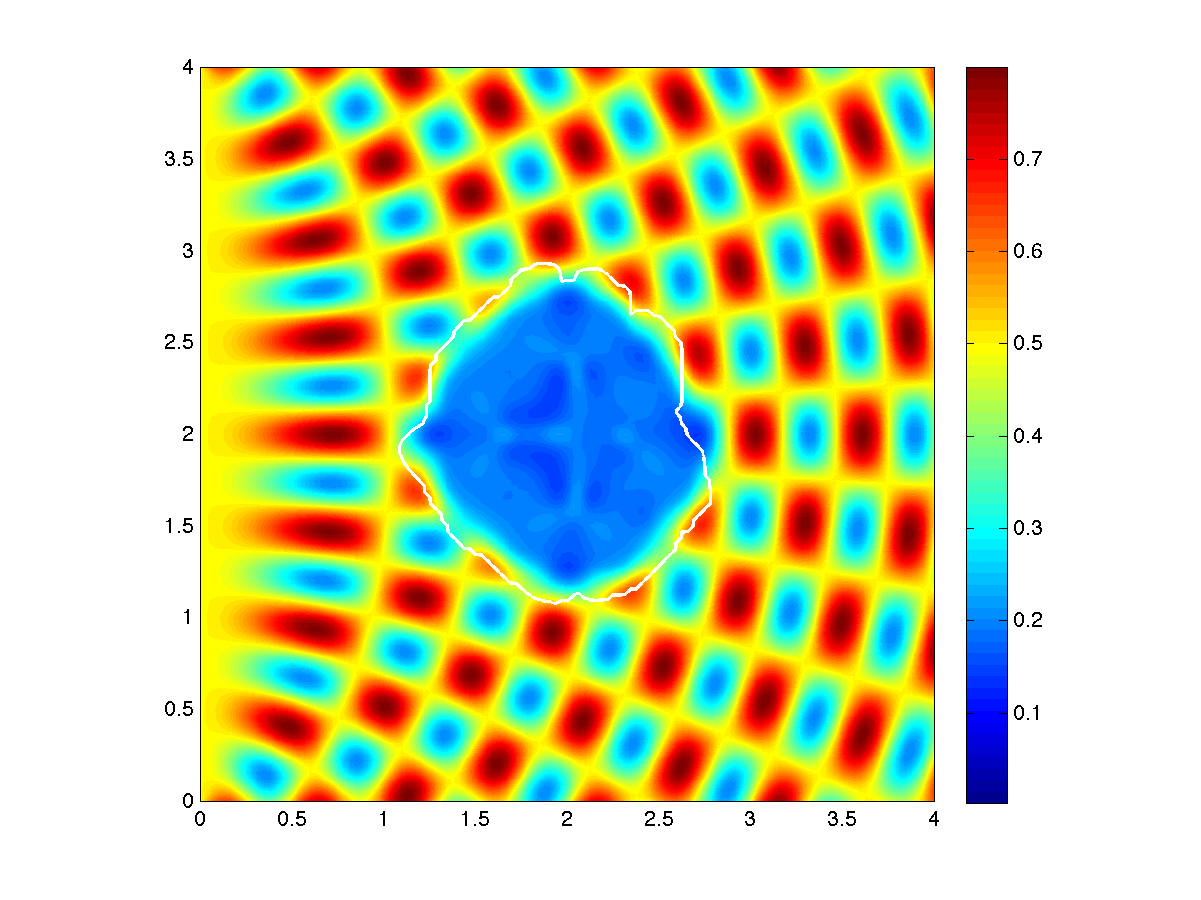}}\\[-20pt]
\subfloat{\label{circle_n} \hspace{5pt}\includegraphics[scale = 0.375]{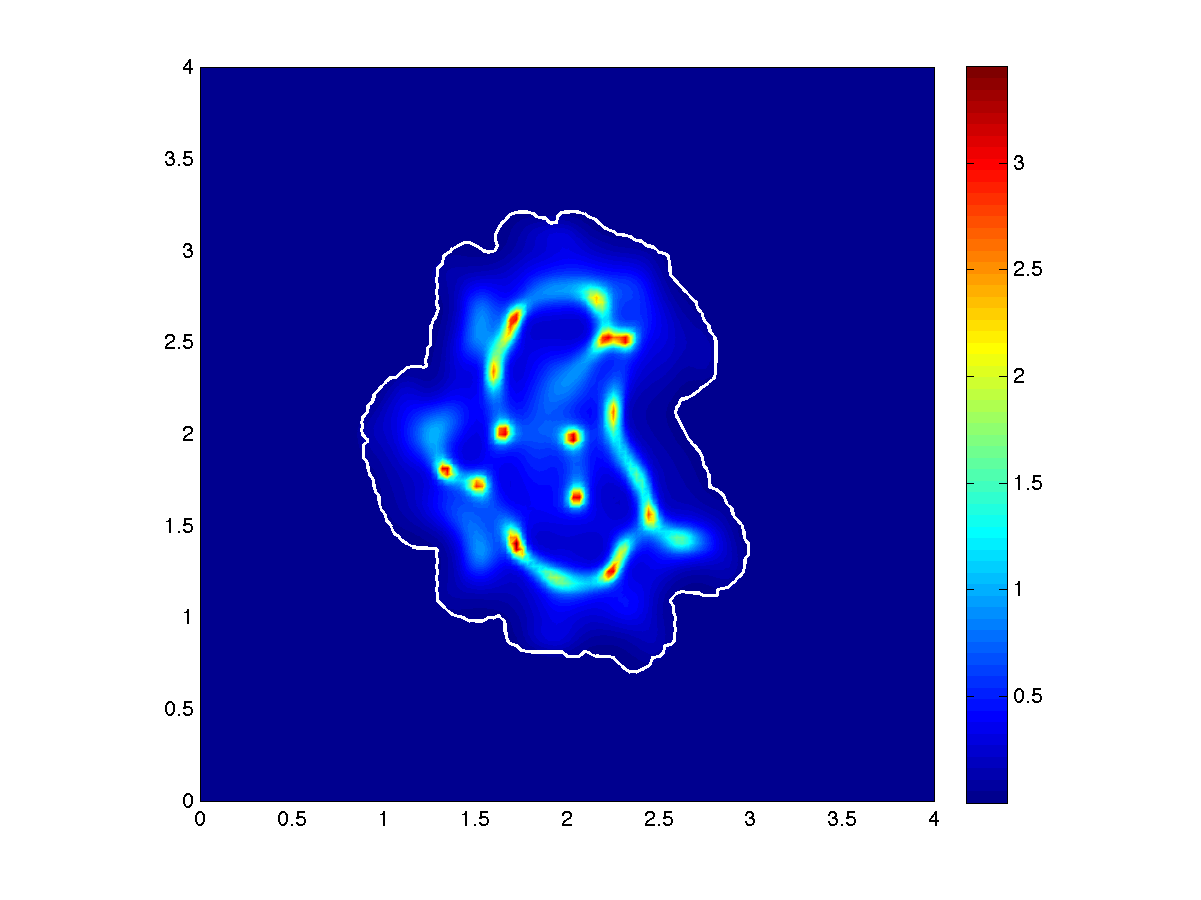}} \hspace{-35pt}
\subfloat{\label{circle_v} \includegraphics[scale = 0.375]{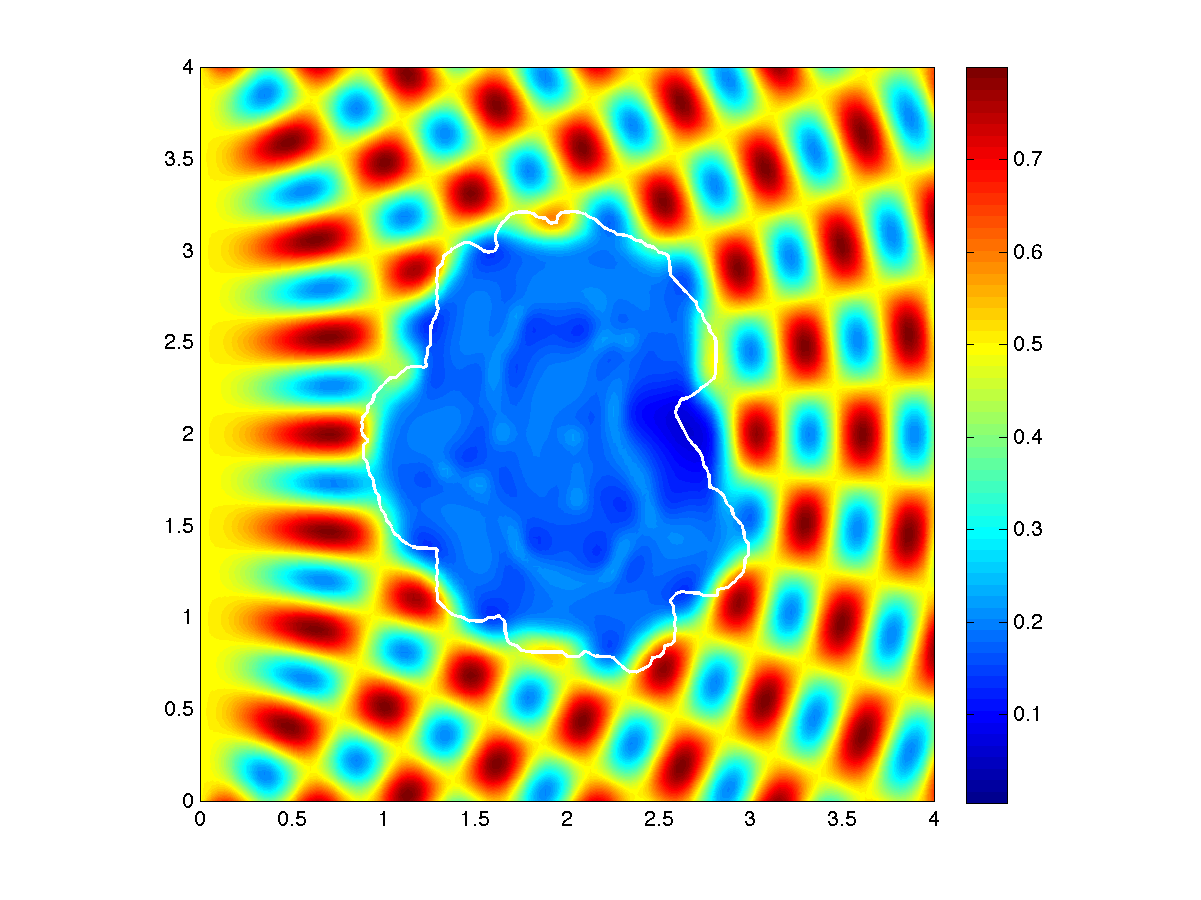}}\\[-20pt]
\subfloat{\label{circle_n} \hspace{5pt}\includegraphics[scale = 0.375]{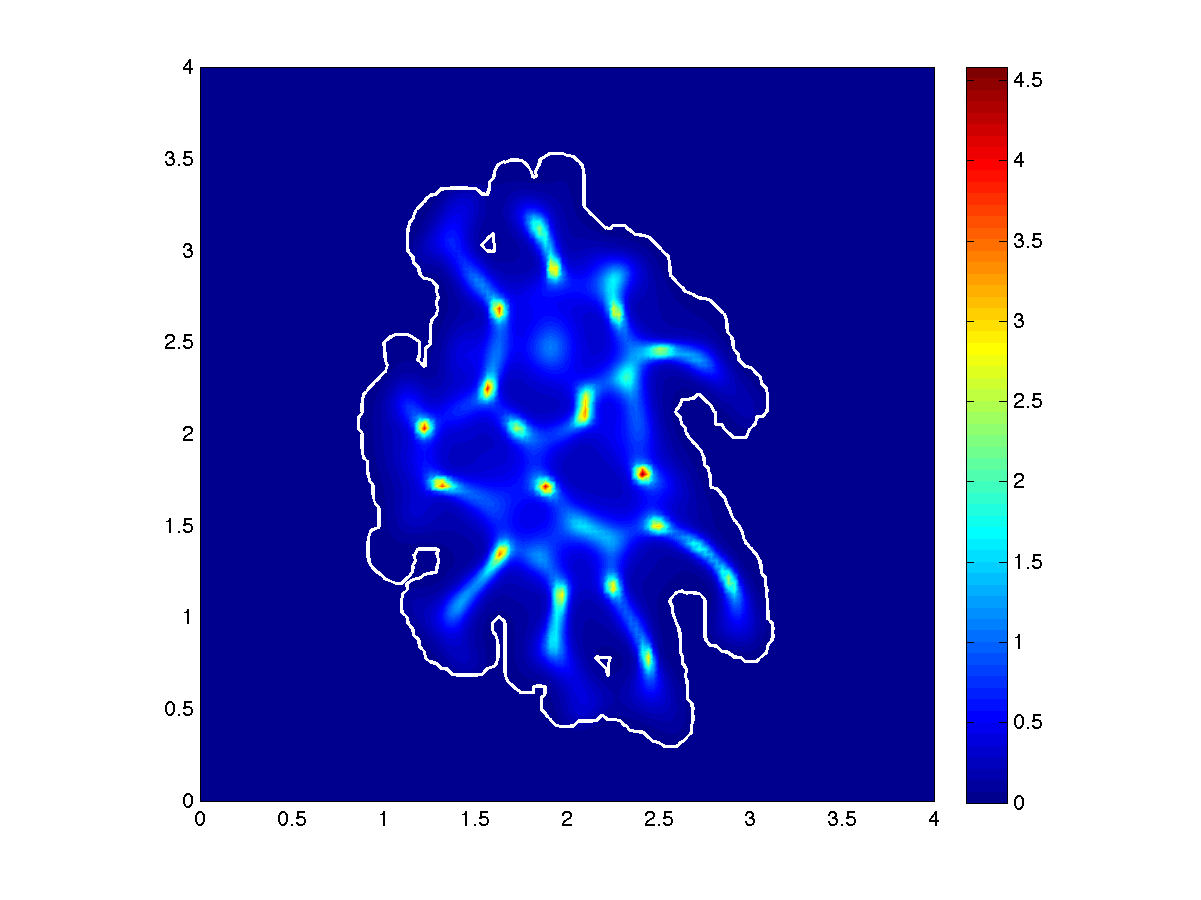}} \hspace{-35pt}
\subfloat{\label{circle_v} \includegraphics[scale = 0.375]{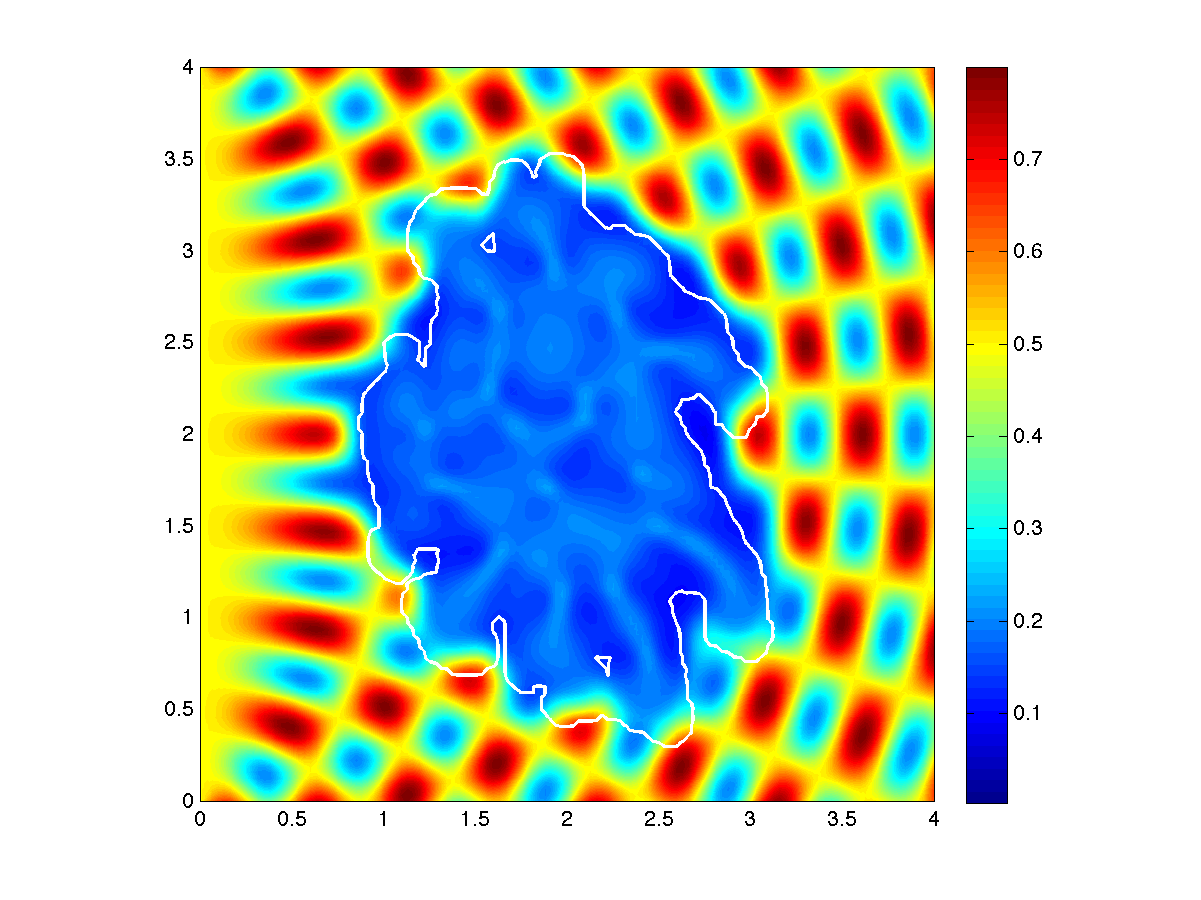}}
\end{tabular}
\vspace{-15pt}
\caption{Simulation results showing distributions of cancer cells (left column) and ECM (right column) and the invasive boundary of the tumour (white line) at various macro-micro stages: Stage 20, 40, 60. \revdt{Starting from the heterogeneous initial conditions shown in Figure \ref{fig:inicond2016}, these results were obtained for $D_c = 4.3 \times 10^{-3}$, $\beta = 0.775$, $\mu_2 = 0$ and $\delta = 0.75$.}}
\label{fig:Dn2}
\end{figure}

As was demonstrated in \cite{Hillen_Painter_2013} and \cite{Painter_Hillen_2011}, the chemotaxis terms in the cancer cell equation~\eqref{equ:Cancer} are the main causes of the occurrence of heterogeneous patterns inside the tumour domain. Therefore, if the cancer cell diffusion coefficient ($D_c$) is increased to be one order magnitude larger than the chemotaxis coefficients ($\chi_u$ and $\chi_p$), which becomes the dominant mechanism of cell movement, then no heterogeneous dynamics will occur inside the tumour as shown in Figure~\ref{fig:Dn1}. In Figure~\ref{fig:Dn2}, the chemotaxis coefficients are one order magnitude larger than the diffusion coefficient, and as a consequence we obtain heterogeneous pattern formation of cancer cells, which leads to a more dynamic tumour boundary deformation. 

\begin{figure}[htp]
\vspace{-15pt}
\hspace{-5em}
\begin{tabular}{cc}
\subfloat{\label{circle_n} \hspace{5pt}\includegraphics[scale = 0.375]{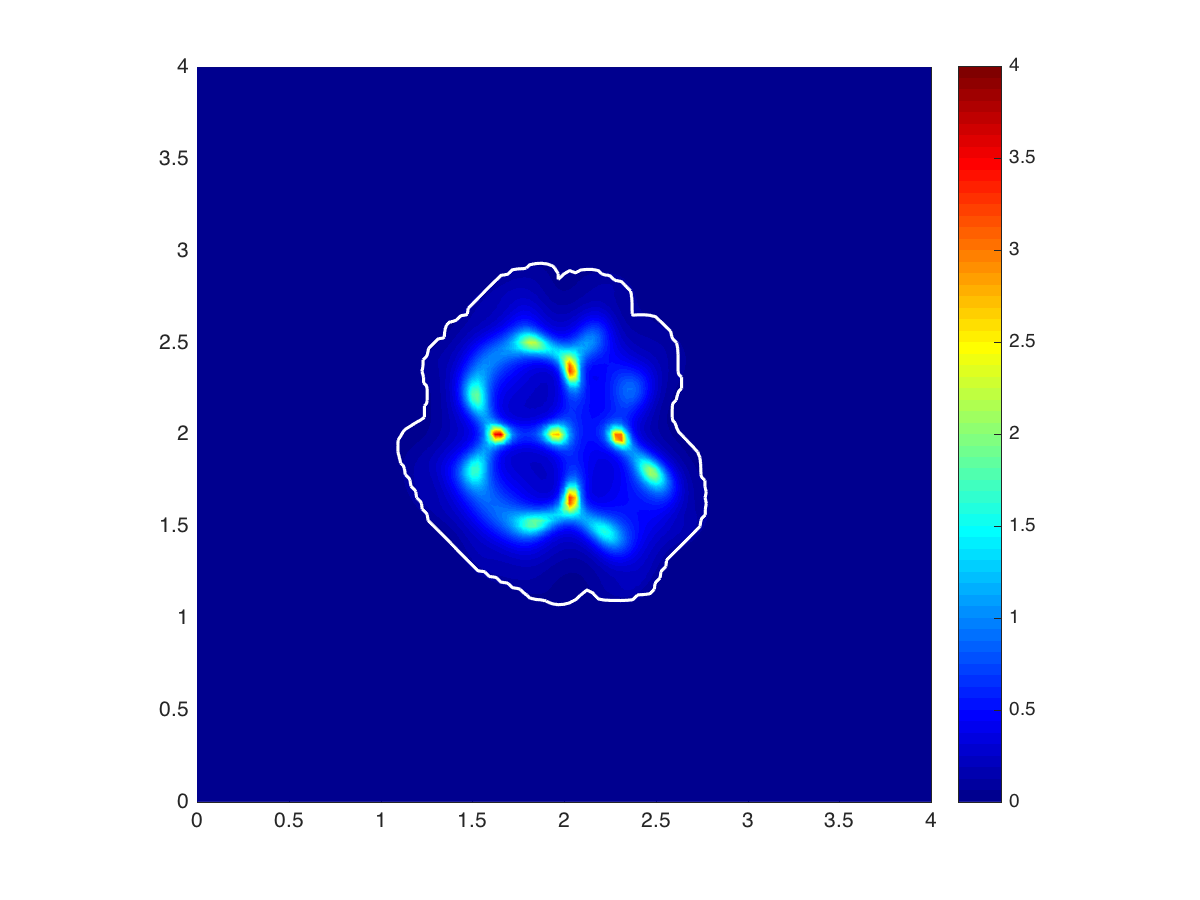}}  \hspace{-35pt}
\subfloat{\label{circle_v} \includegraphics[scale = 0.375]{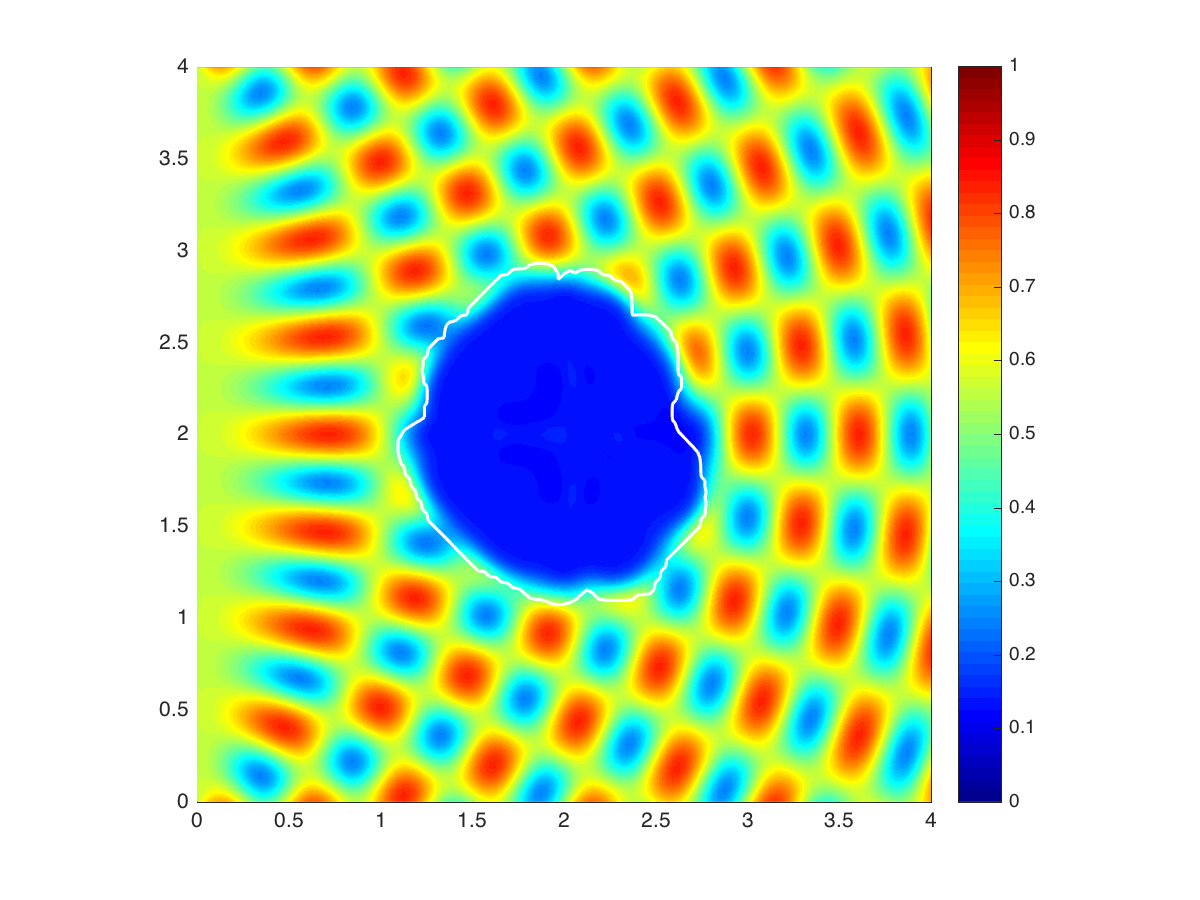}}\\[-20pt]
\subfloat{\label{circle_n} \hspace{5pt}\includegraphics[scale = 0.375]{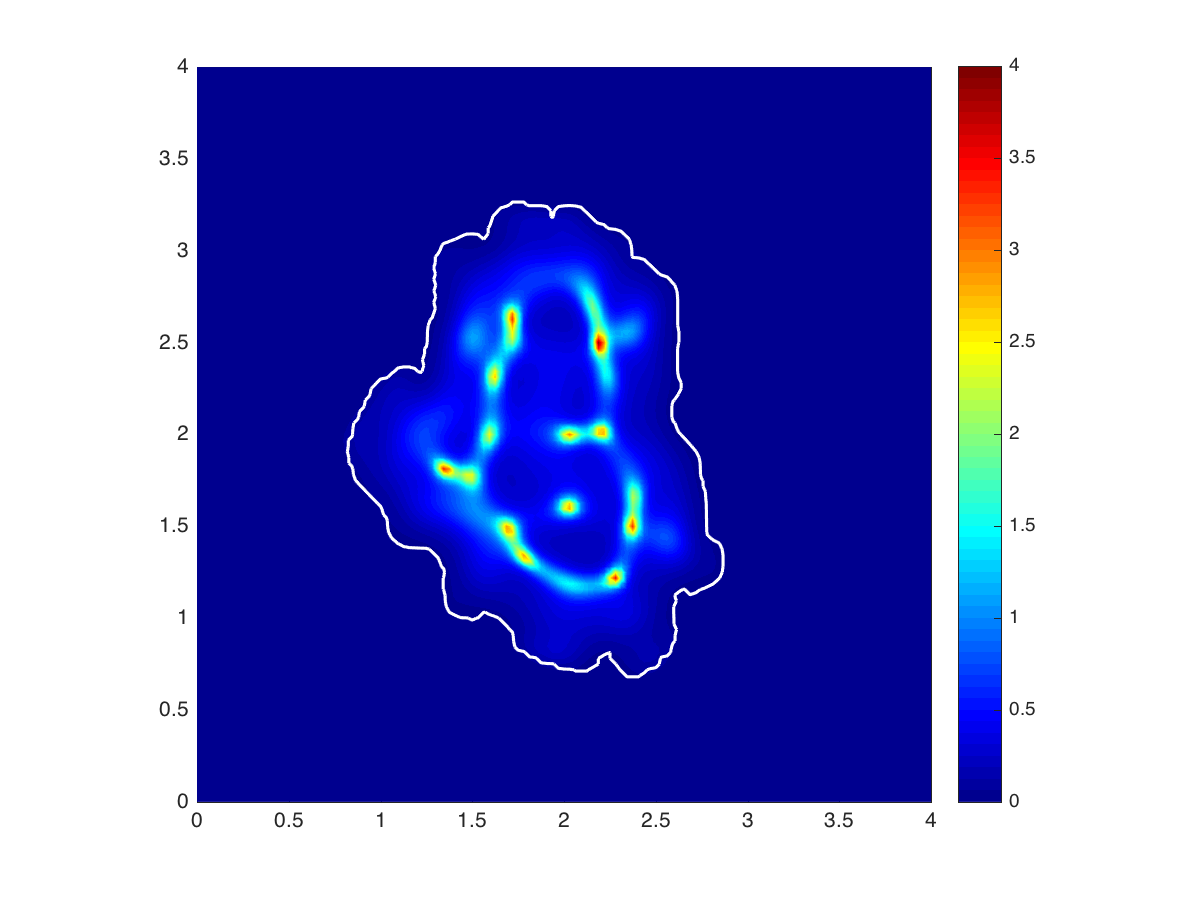}}  \hspace{-35pt}
\subfloat{\label{circle_v} \includegraphics[scale = 0.375]{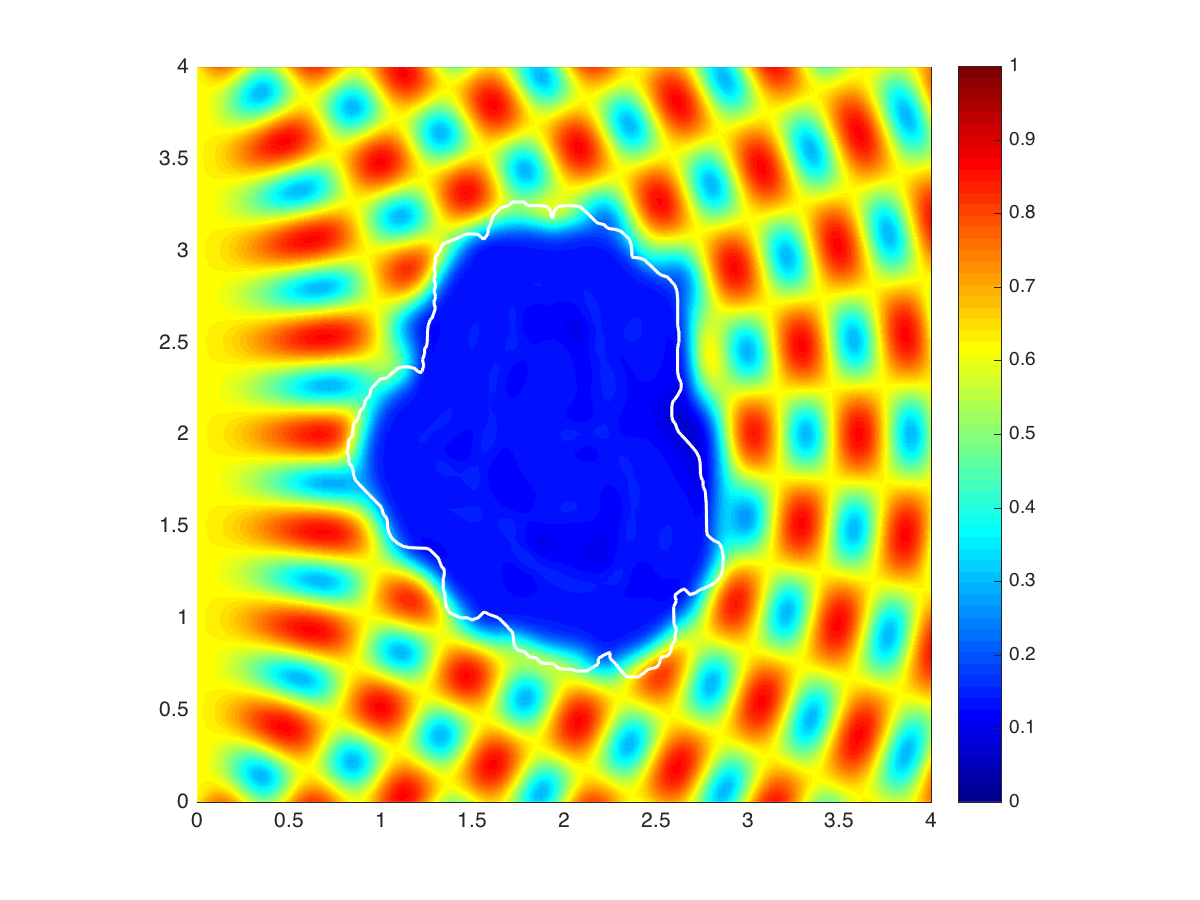}}\\[-20pt]
\subfloat{\label{circle_n} \hspace{5pt}\includegraphics[scale = 0.375]{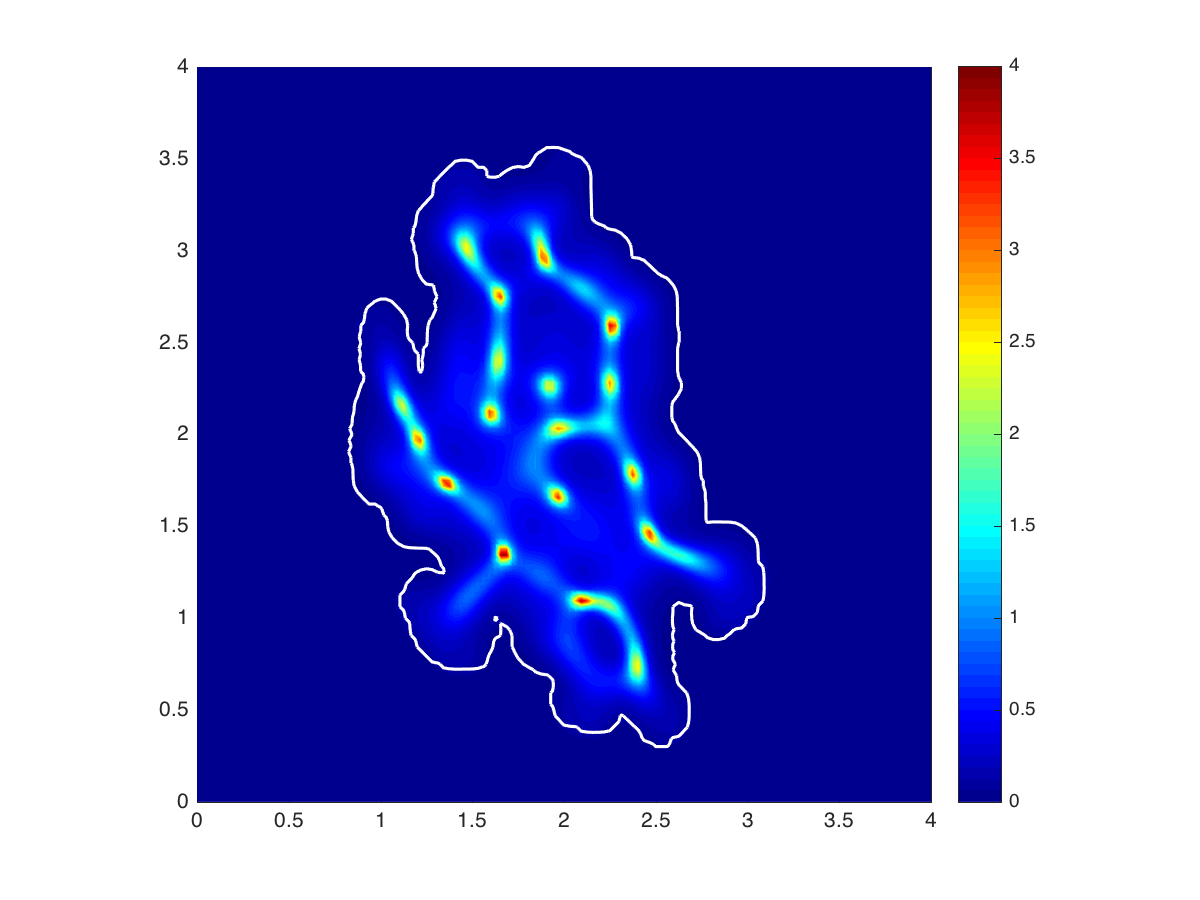}}   \hspace{-35pt}
\subfloat{\label{circle_v} \includegraphics[scale = 0.375]{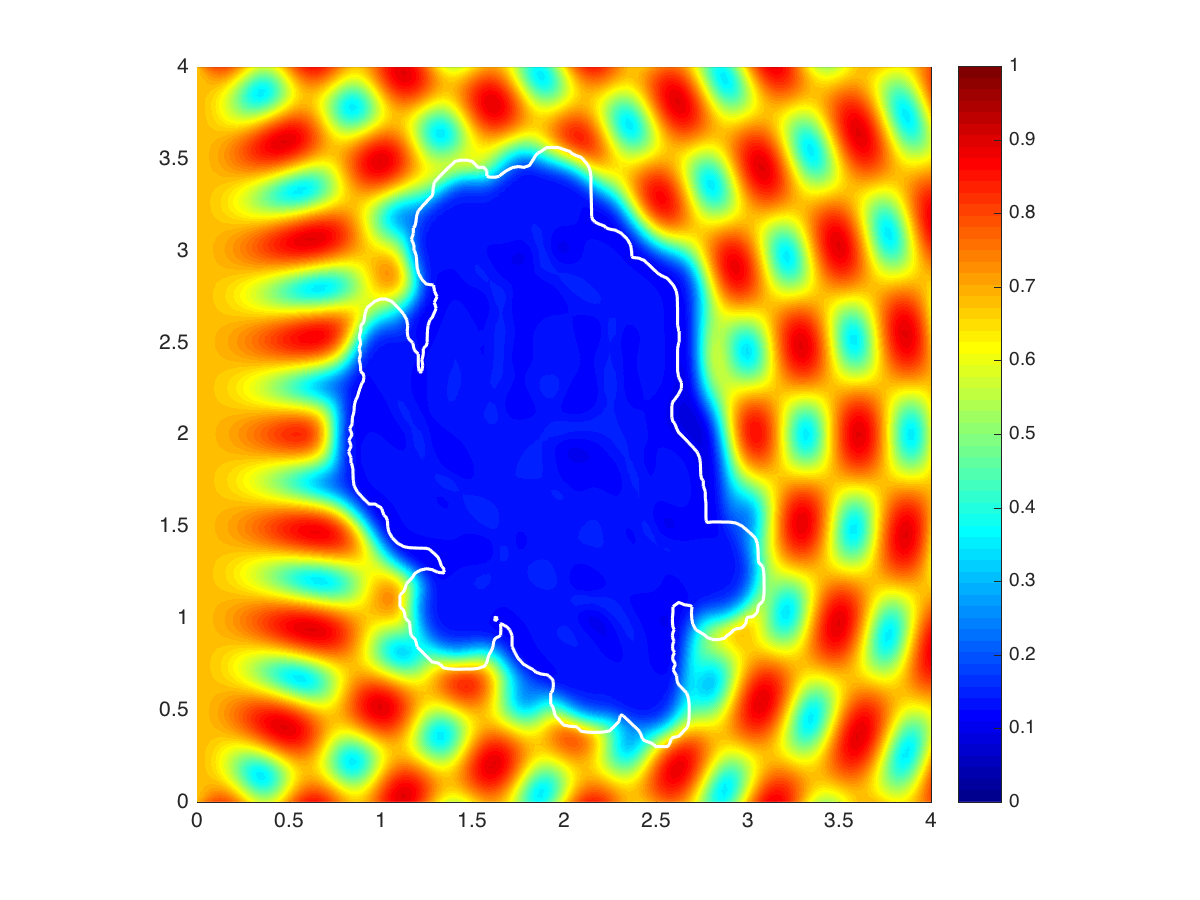}}
\end{tabular}
\vspace{-15pt}
\caption{Simulation results showing distributions of cancer cells (left column) and ECM (right column) and the invasive boundary of the tumour (white line) at various macro-micro stages: Stage 20, 40, 60. \revdt{Starting from the heterogeneous initial conditions shown in Figure \ref{fig:inicond2016}, these results were obtained for $D_c = 4.3 \times 10^{-3}$, $\beta = 0.7625$, $\mu_2 = 0.005$ and $\delta = 1.5$.}}
\label{fig:big}
\end{figure}

\begin{figure}[htp]
\vspace{-15pt}
\hspace{-5em}
\begin{tabular}{cc}
\subfloat{\label{circle_n} \hspace{5pt}\includegraphics[scale = 0.375]{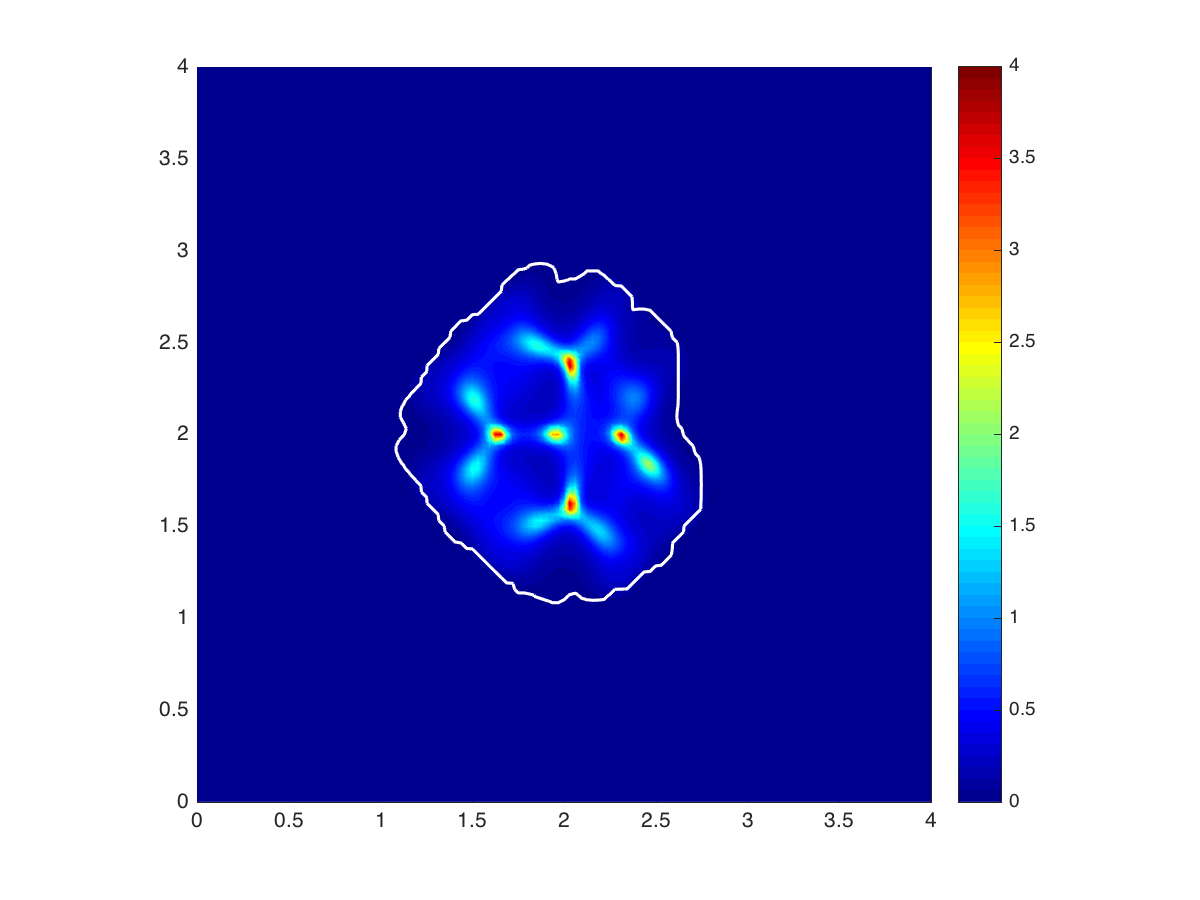}}  \hspace{-35pt}
\subfloat{\label{circle_v} \includegraphics[scale = 0.375]{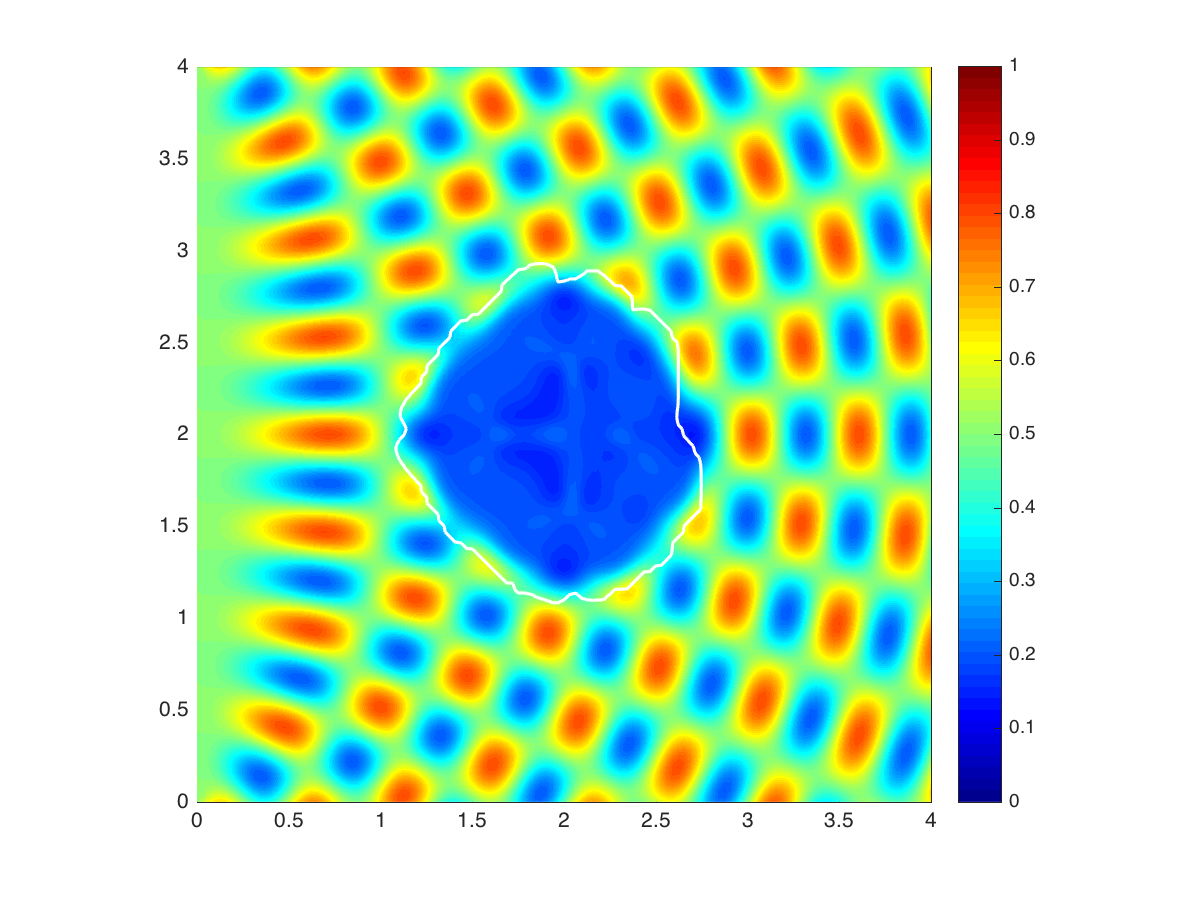}}\\
\subfloat{\label{circle_n} \hspace{5pt}\includegraphics[scale = 0.375]{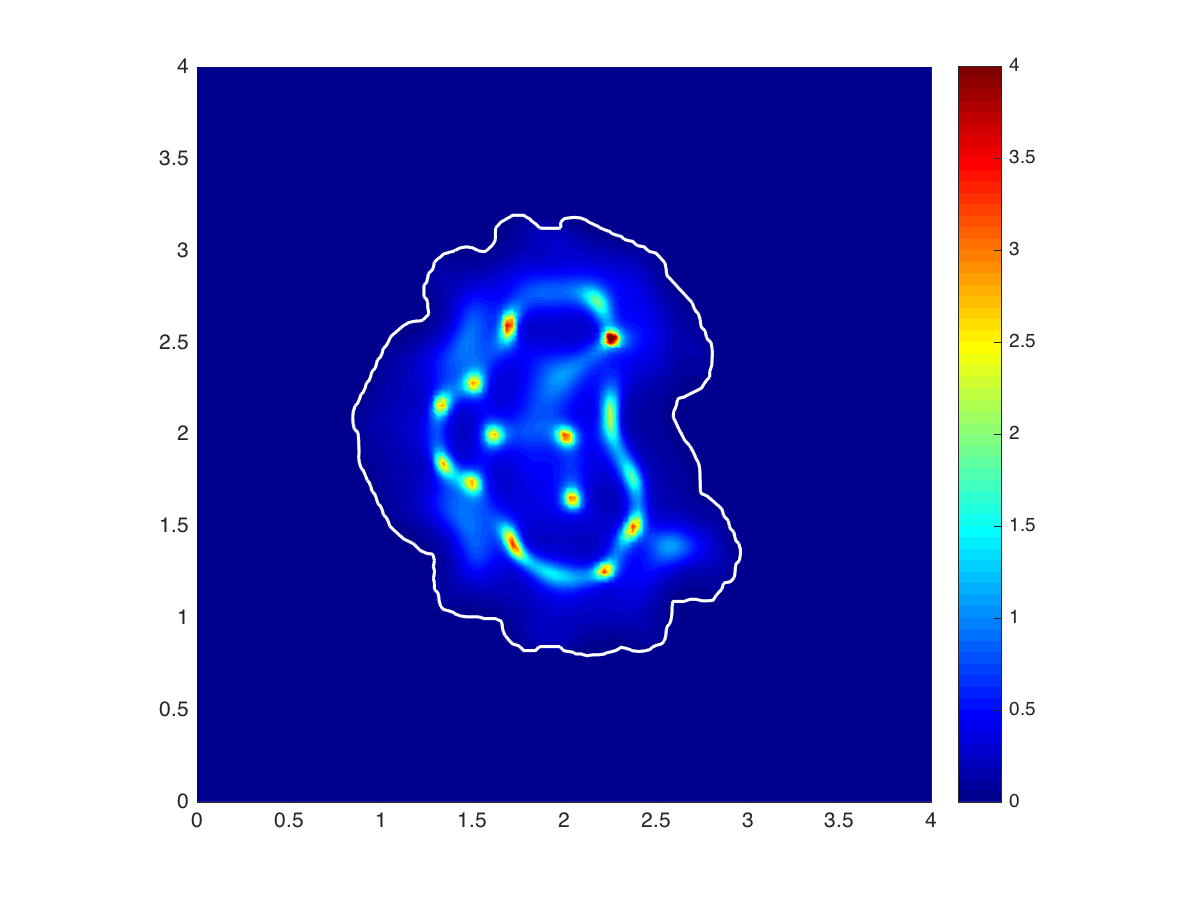}}  \hspace{-35pt}
\subfloat{\label{circle_v} \includegraphics[scale = 0.375]{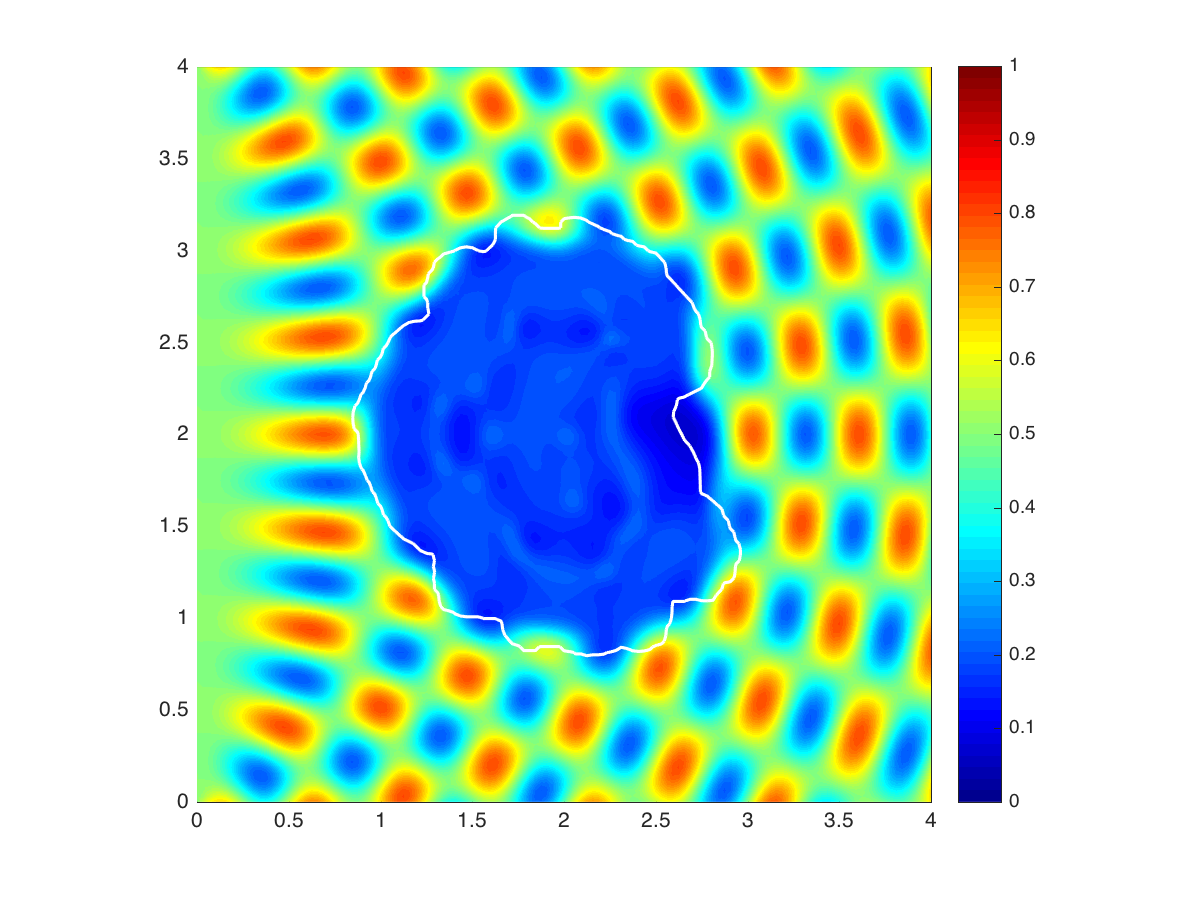}}\\
\subfloat{\label{circle_n} \hspace{5pt}\includegraphics[scale = 0.375]{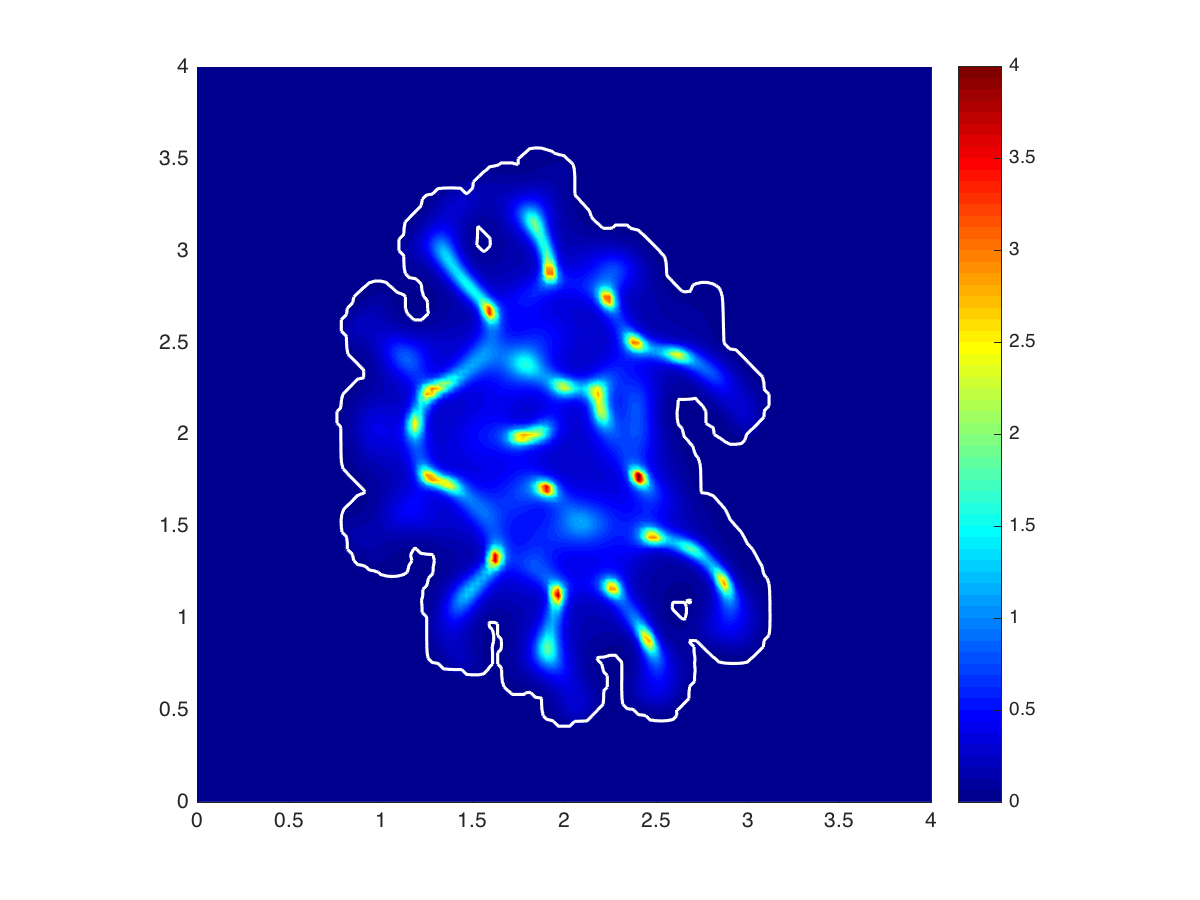}}   \hspace{-35pt}
\subfloat{\label{circle_v} \includegraphics[scale = 0.375]{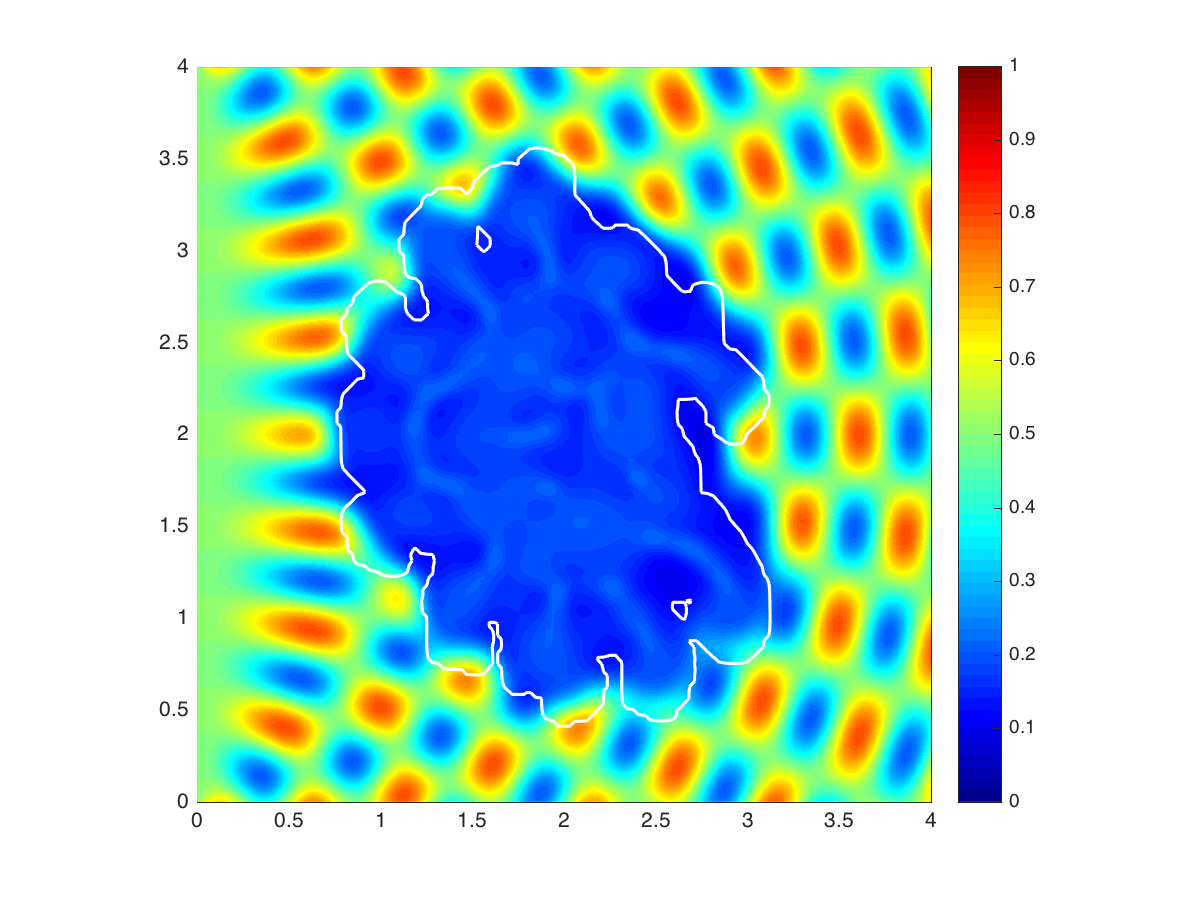}}
\end{tabular}
\vspace{-15pt}
\caption{Simulation results showing distributions of cancer cells (left column) and ECM (right column) and the invasive boundary of the tumour (white line) at various macro-micro stages: Stage 20, 40, 60. \revdt{Starting from the heterogeneous initial conditions shown in Figure \ref{fig:inicond2016}, these results were obtained for $D_c = 4.3 \times 10^{-3}$, $\beta = 0.7625$, $\mu_2 = 0$ and $\delta = 0.75$.}}
\label{fig:small}
\end{figure}

\paragraph{\it{ECM proliferation rate $\mu_2$ \& degradation rate $\delta$.}}
From all the simulation results presented so far, we conclude that the degradation of ECM facilitates cancer invasion. However, \revdt{as shown also in several interdisciplinary biological investigations \cite{Gatenby_Gawlinski_1996,Gatenby_et_al_2006,StetlerStevenson1993}}, the invasion process will stop where ECM is degraded a lot. This captures the biological scenario in 2D that when cell-matrix adhesion is too low, no focal adhesions or stress fibres are formed, and the cells do not move. In order to investigate the effect of ECM proliferation and degradation on the invasion process, we compare two groups of parameters: 1) $\mu_2 = 0.01$, $\delta = 1.5$ (non-zero proliferation rate with relatively large degradation rate, Figure~\ref{fig:big}); 2) $\mu_2 = 0$, $\delta = 0.75$ (no proliferation with relatively small degradation rate, Figure~\ref{fig:small}). From these two groups of images, we observe that when proliferation is present and the degradation rate is relatively large, deformations of the boundary is not as dynamic as that when the proliferation term is absent with a relatively small degradation rate. The reason could be that the proliferation term will reduce the degree of heterogeneity of the distribution of ECM, which leads to a less fingered spreading of the cancer cell population.

\begin{figure}[htp]
\vspace{-15pt}
\hspace{-5em}
\begin{tabular}{cc}
\subfloat{\label{circle_n} \hspace{5pt}\includegraphics[scale = 0.375]{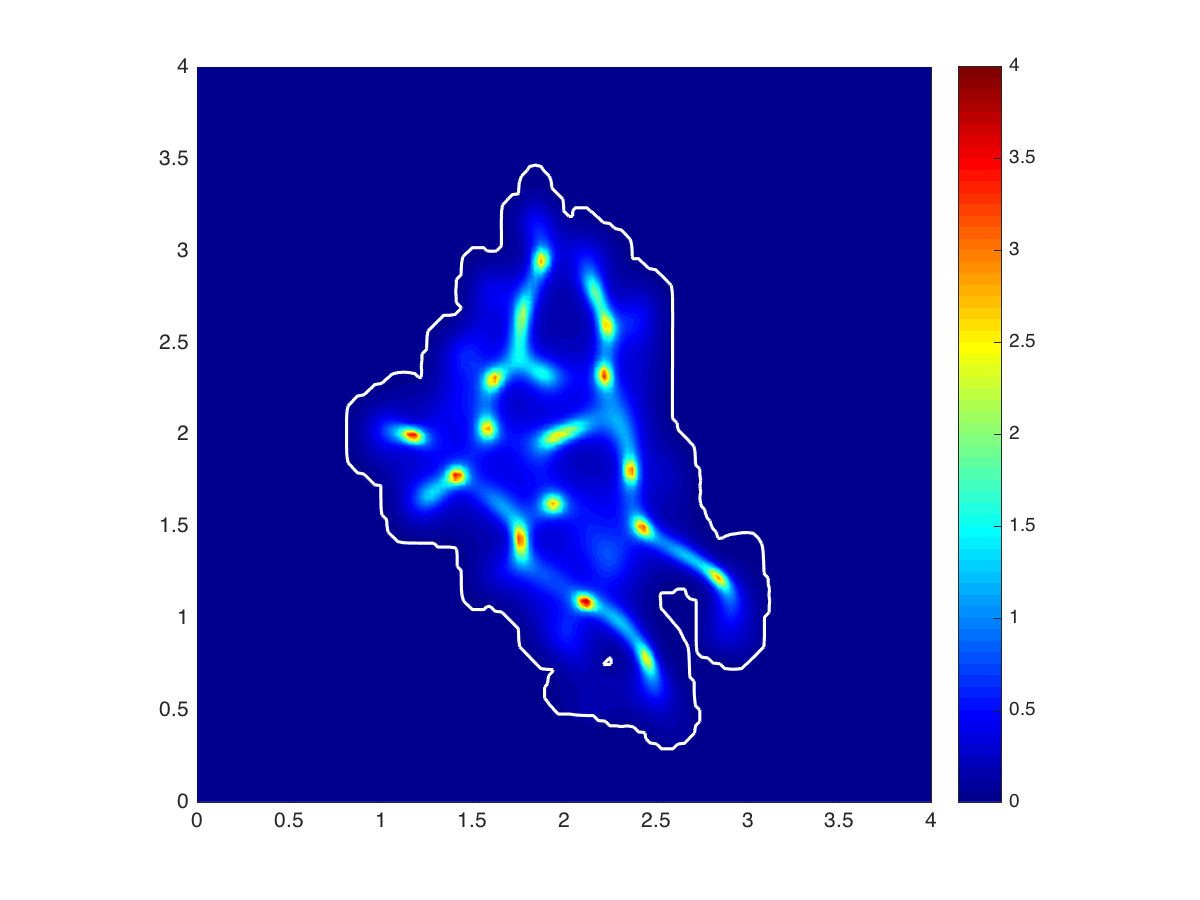}} \hspace{-35pt}
\subfloat{\label{circle_v} \includegraphics[scale = 0.375]{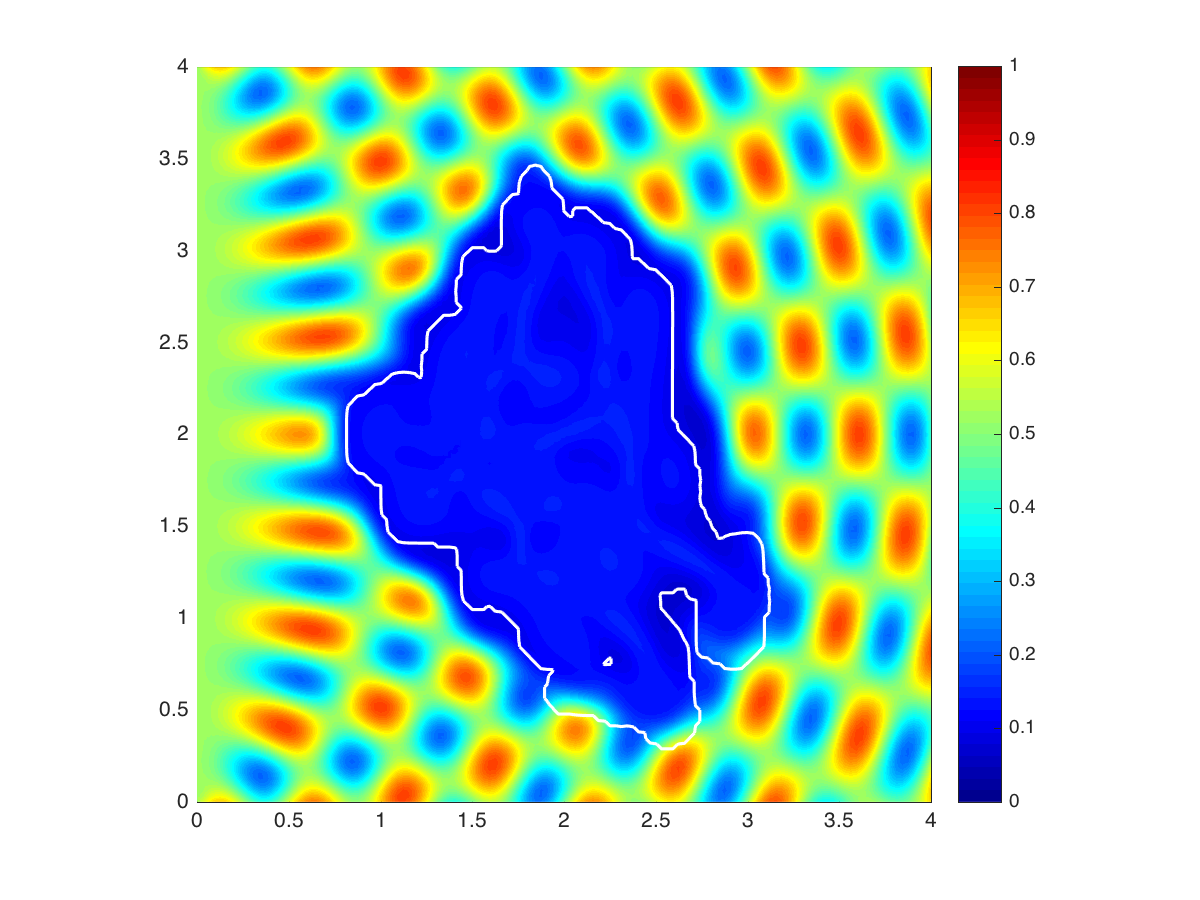}}\\
\subfloat{\label{circle_n} \hspace{5pt}\includegraphics[scale = 0.375]{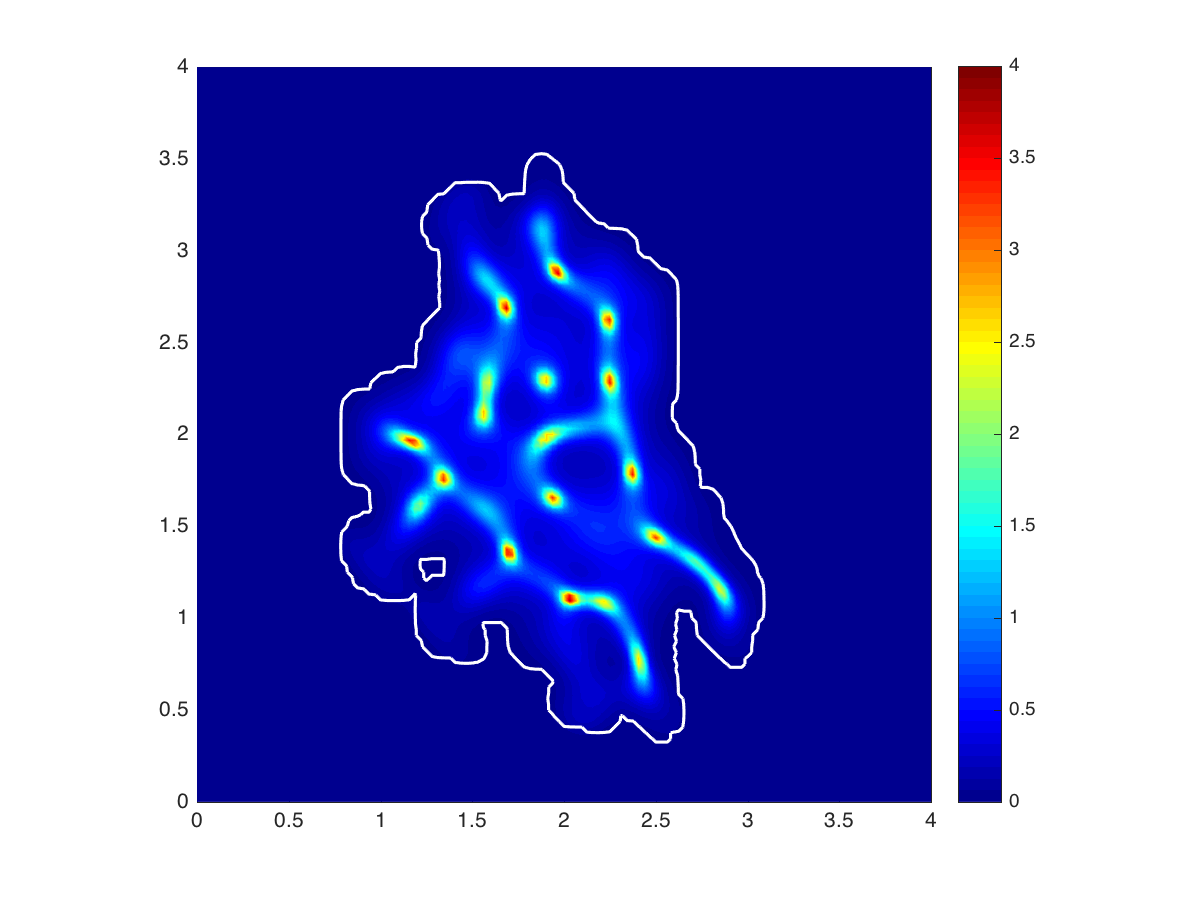}}  \hspace{-35pt}
\subfloat{\label{circle_v} \includegraphics[scale = 0.375]{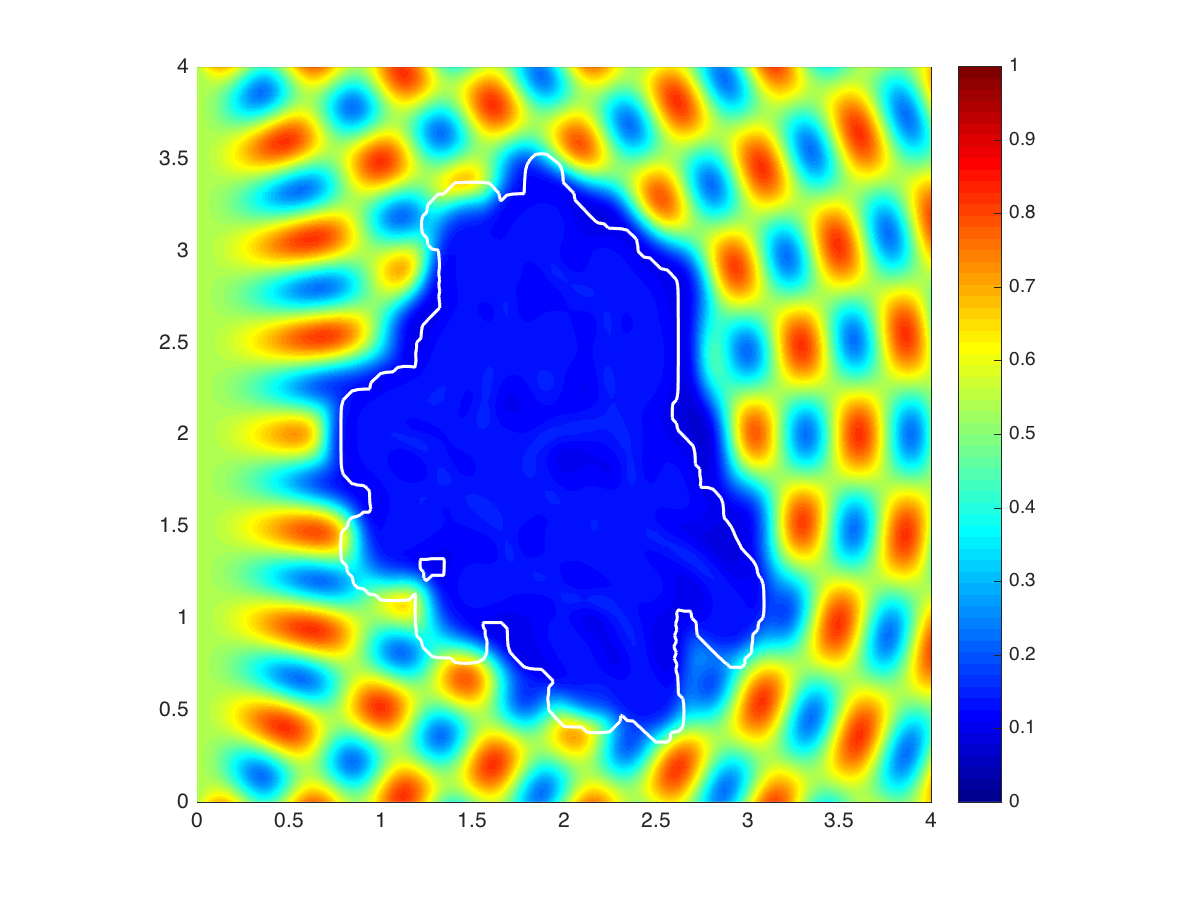}}\\
\subfloat{\label{circle_n} \hspace{5pt}\includegraphics[scale = 0.375]{Cancer60_mu2_0_005}}   \hspace{-35pt}
\subfloat{\label{circle_v} \includegraphics[scale = 0.375]{ECM60_mu2_0_005}}
\end{tabular}
\vspace{-15pt}
\caption{Simulation results showing distributions of cancer cells (left column) and ECM (right column) and the invasive boundary of the tumour (white line) at macro-micro stage 60. \revdt{Starting from the heterogeneous initial conditions shown in Figure \ref{fig:inicond2016}, these results were obtained for $D_c = 4.3 \times 10^{-3}$, $\beta = 0.7625$, $\delta = 1.5$, and for rows 1 to 3 of images we consider $\mu_2=0.0005$, $\mu_2=0.001$, and $\mu_2=0.005$, respectively.}}
\label{fig:mu2_15mar2016}
\end{figure}

\begin{figure}[htp]
\vspace{-15pt}
\hspace{-5em}
\begin{tabular}{cc}
\subfloat{\label{circle_n} \hspace{5pt}\includegraphics[scale = 0.375]{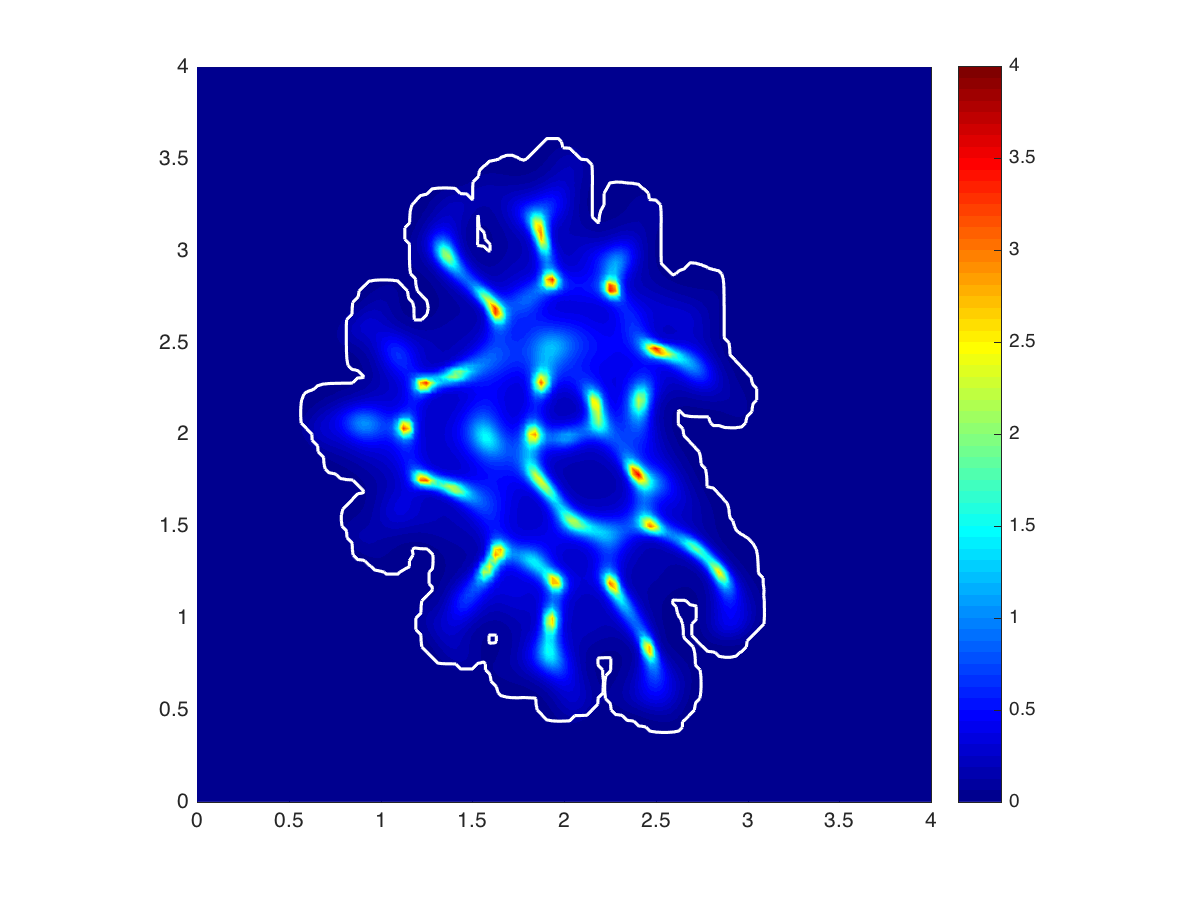}}   \hspace{-35pt}
\subfloat{\label{circle_v} \includegraphics[scale = 0.375]{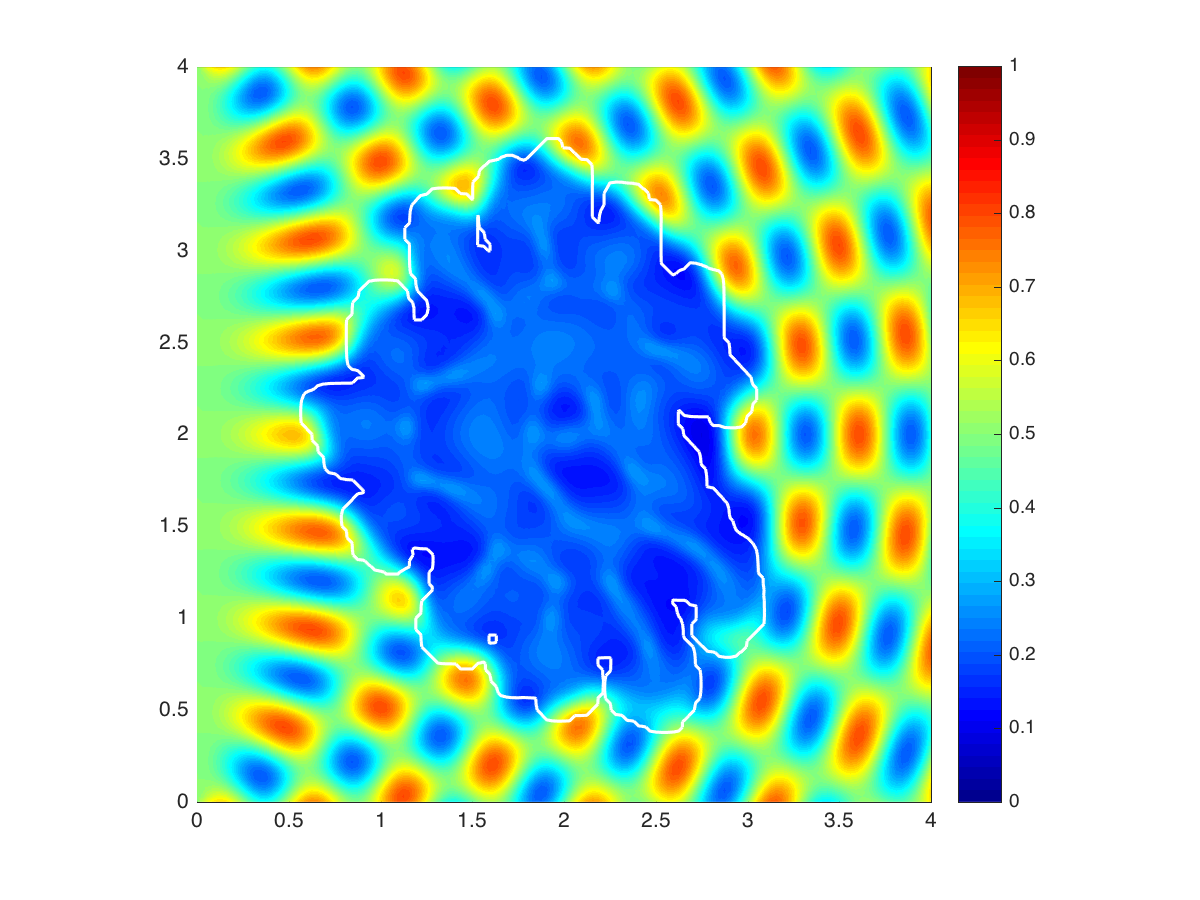}}\\
\subfloat{\label{circle_n} \hspace{5pt}\includegraphics[scale = 0.375]{Cancer60_delta_0_75}}  \hspace{-35pt}
\subfloat{\label{circle_v} \includegraphics[scale = 0.375]{ECM60_delta_0_75}}\\
\subfloat{\label{circle_n} \hspace{5pt}\includegraphics[scale = 0.375]{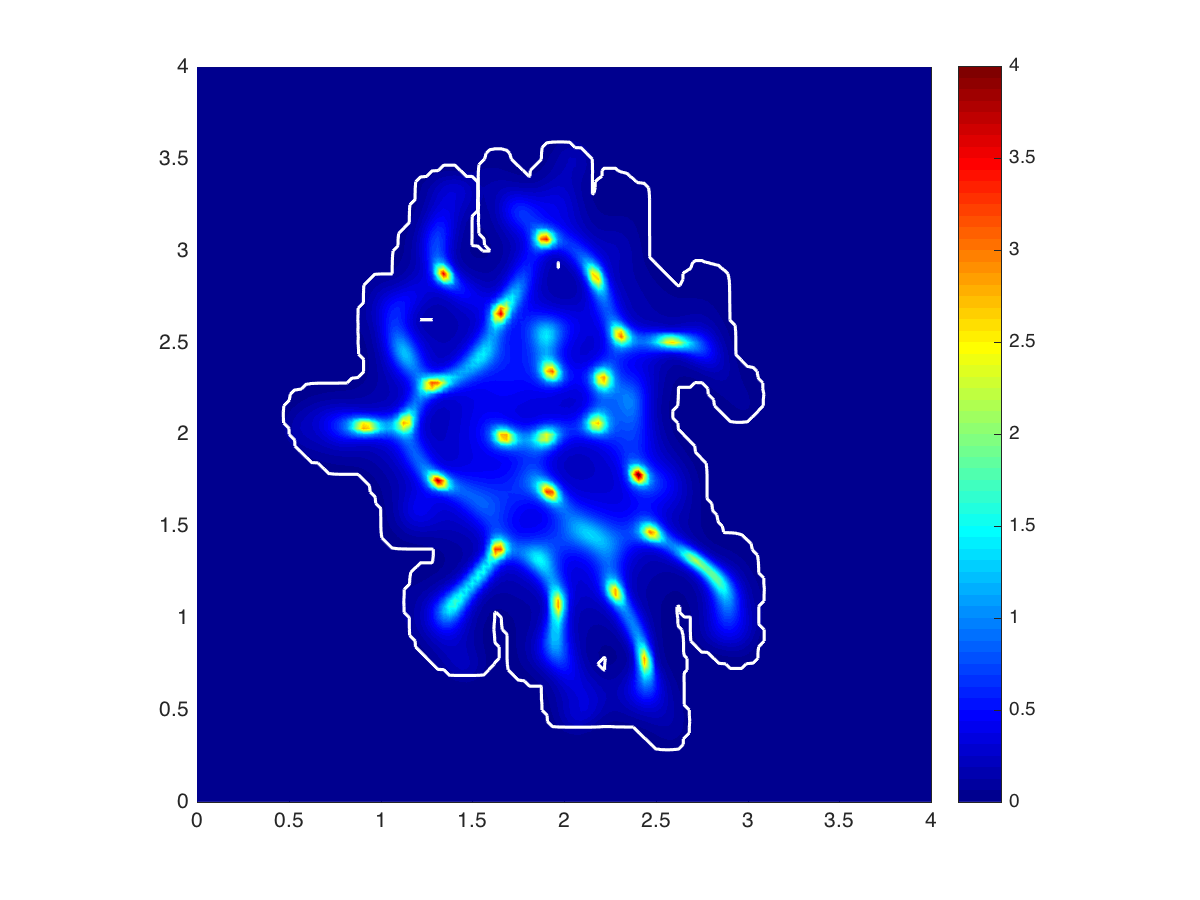}} \hspace{-35pt}
\subfloat{\label{circle_v} \includegraphics[scale = 0.375]{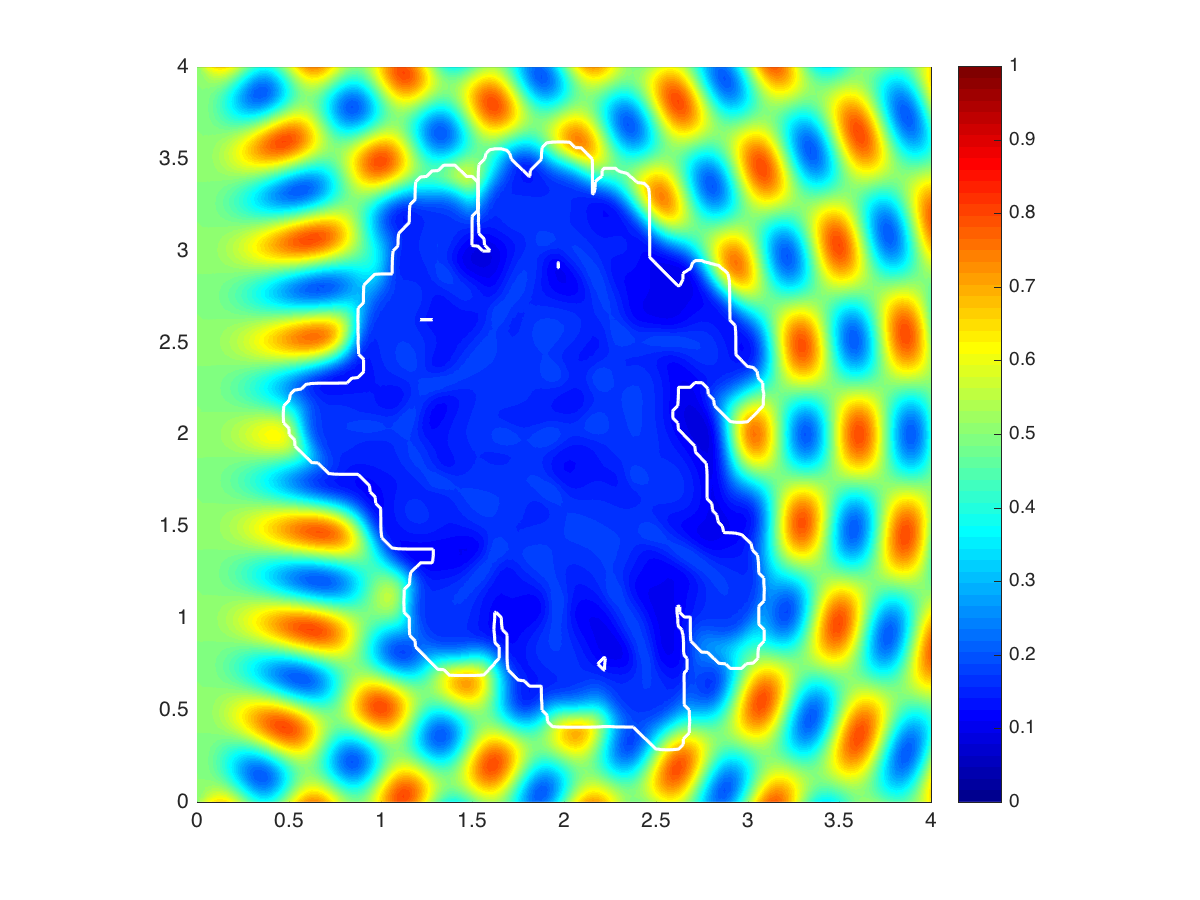}}
\end{tabular}
\vspace{-15pt}
\caption{Simulation results showing distributions of cancer cells (left column) and ECM (right column) and the invasive boundary of the tumour (white line) at macro-micro stage 60. \revdt{Starting from the heterogeneous initial conditions shown in Figure \ref{fig:inicond2016}, these results were obtained for $D_c = 4.3 \times 10^{-3}$, $\mu_2=0$, $\beta = 0.7625$, and for rows 1 to 3 of images we consider $\delta = 0.5$, $\delta = 0.75$, and $\delta = 1$, respectively.}}
\label{fig:delta_15mar2016}
\end{figure}

\revdt{Figure \ref{fig:mu2_15mar2016} shows the results of simulations at the macro-micro stage $60$ where the ECM proliferation rate parameter $\mu_2$ was increased over a range of values in the interval $[0.0005, 0.005]$. The figures show that, over the range chosen, there is little difference in either the overall extent of invasion or the morphology of the invading cancer. Figure \ref{fig:delta_15mar2016} shows the results of simulations at the same macro-micro stage $60$ where the degradation parameter $\delta$ was increased over the interval $[0.5,1]$. The figures show that, over the range chosen, the extent of the invasion is similar, but the morphology of the invading cancer changes slightly from a more fingered boundary to a less fingered boundary. Finally, as the two parameters $\mu_2$ and $\delta$ are varied, both Figure \ref{fig:mu2_15mar2016} and Figure \ref{fig:delta_15mar2016} exhibit spatial consistency aspects in the invasion pattern. 
}

\paragraph{\it{Threshold coefficient $\beta$.}}

\begin{figure}[htp]
\vspace{-15pt}
\hspace{-5em}
\begin{tabular}{cc}
\subfloat{\label{circle_n} \hspace{5pt}\includegraphics[scale = 0.375]{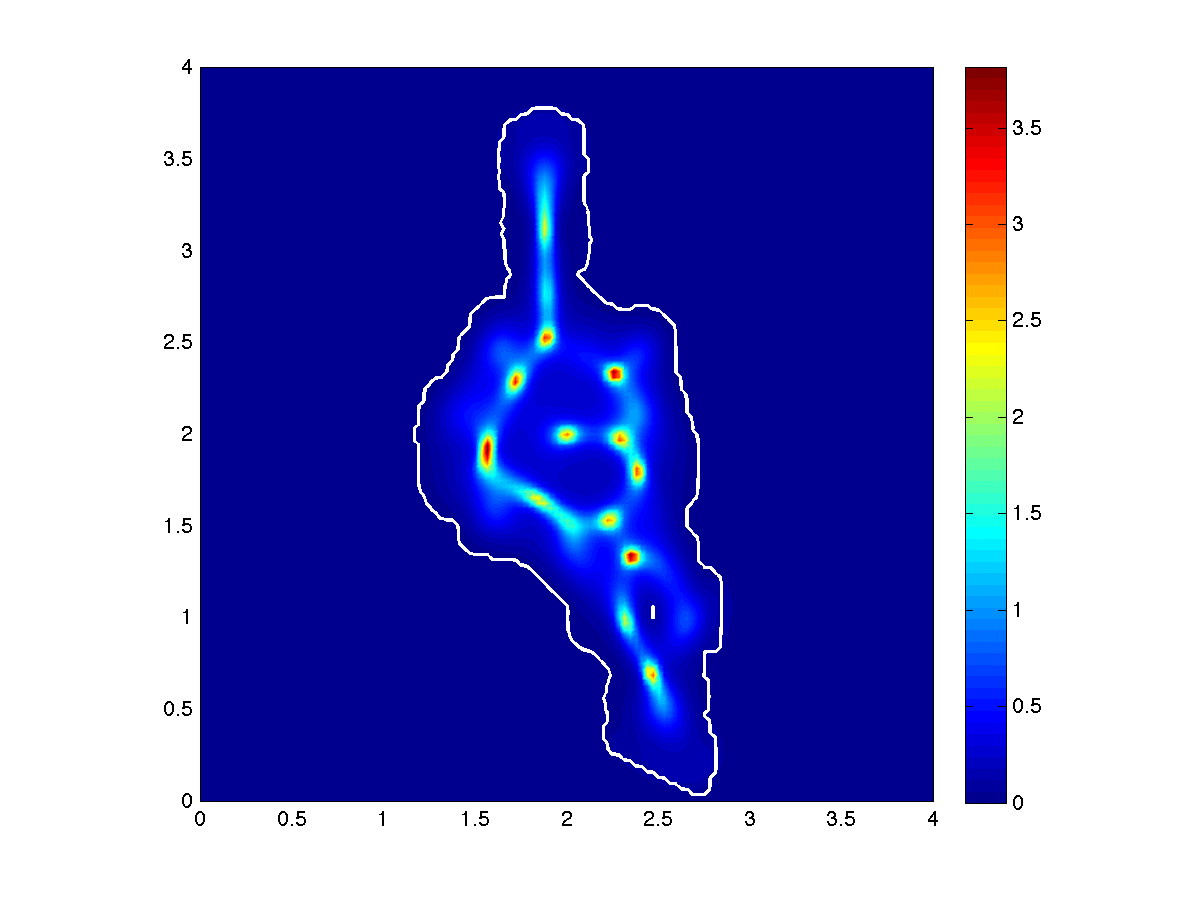}}   \hspace{-35pt}
\subfloat{\label{circle_v} \includegraphics[scale = 0.375]{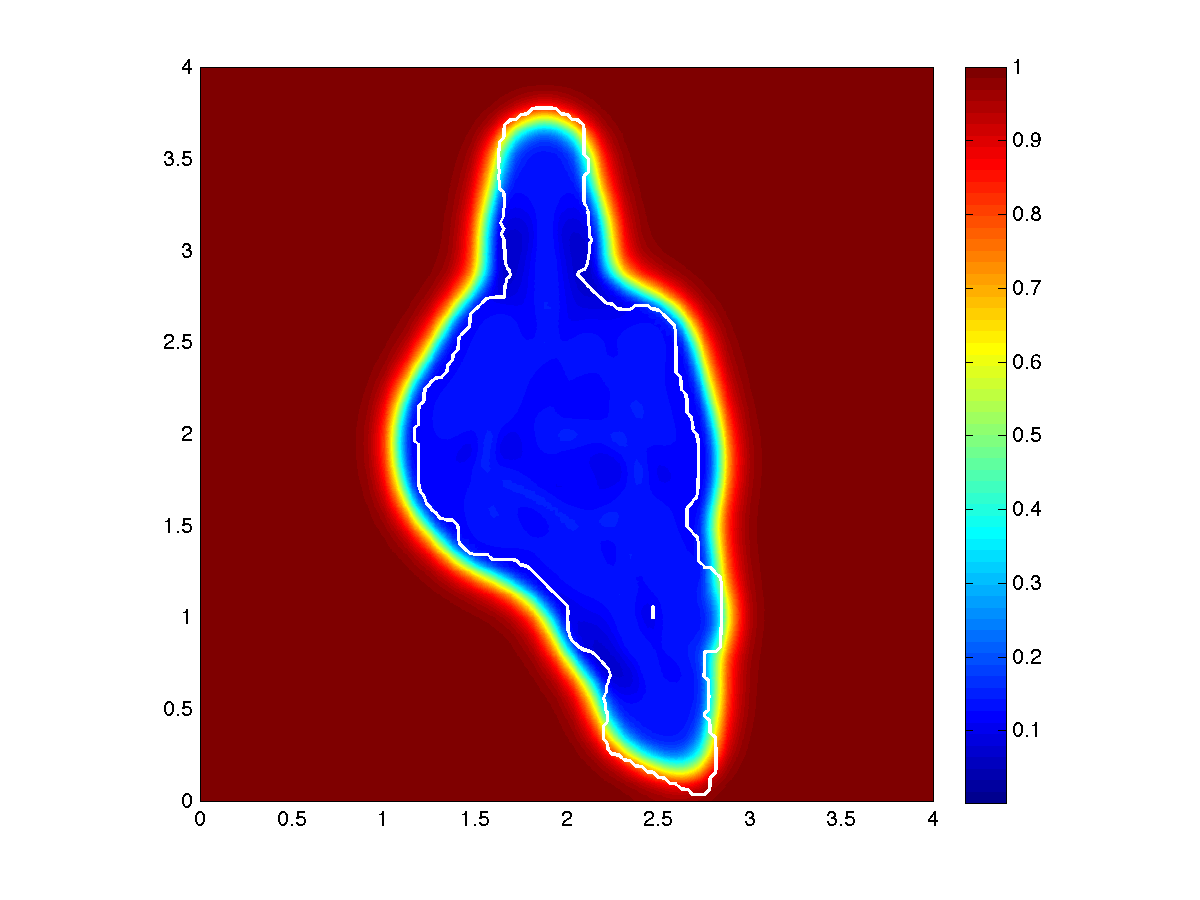}}\\
\subfloat{\label{circle_n} \hspace{5pt}\includegraphics[scale = 0.375]{Cancer60_beta0775_Dn2_homo_big}}  \hspace{-35pt}
\subfloat{\label{circle_v} \includegraphics[scale = 0.375]{ECM60_beta0775_Dn2_homo_big}}\\
\subfloat{\label{circle_n} \hspace{5pt}\includegraphics[scale = 0.375]{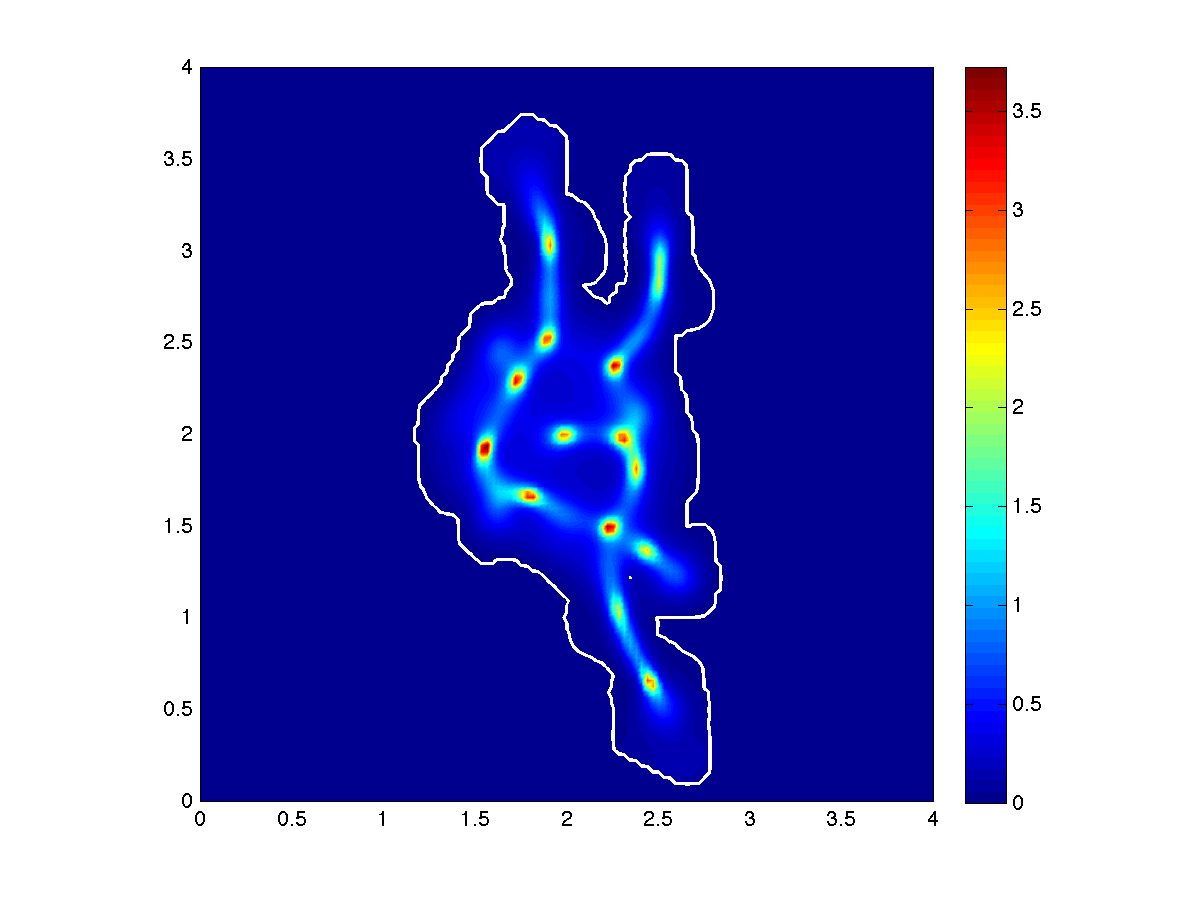}} \hspace{-35pt}
\subfloat{\label{circle_v} \includegraphics[scale = 0.375]{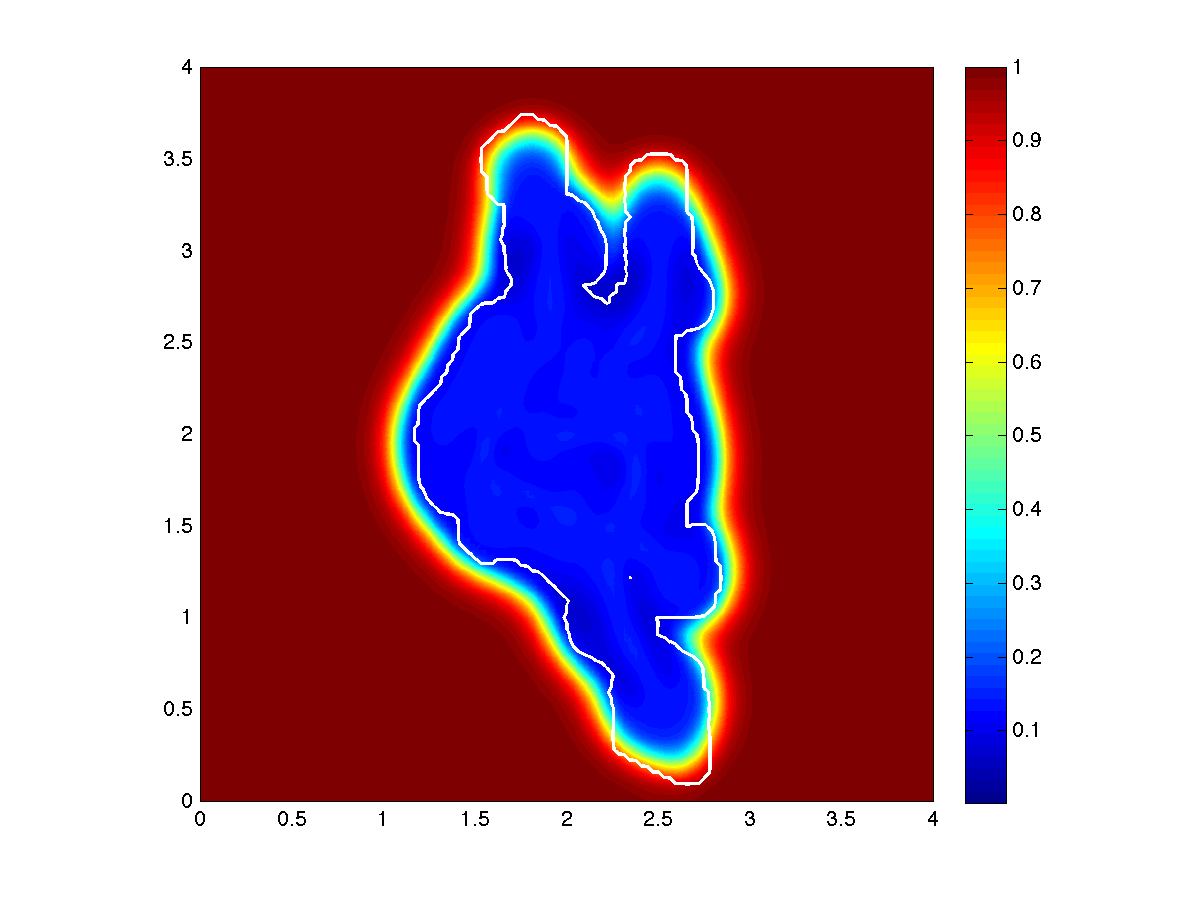}}
\end{tabular}
\vspace{-15pt}
\caption{Simulation results showing distributions of cancer cells (left column) and ECM (right column) and the invasive boundary of the tumour (white line) at macro-micro stage 60. \revdt{Starting from the homogeneous initial conditions shown in Figure \ref{fig:inicond2016}, these results were obtained for $D_c = 4.3 \times 10^{-3}$, $\delta = 1.5$, $\mu_2 = 0.01$, and for rows 1 to 3 of images we consider $\beta=0.7625$, $\beta=0.775$, and $\beta=0.7875$, respectively.}}
\label{fig:beta07625}
\end{figure}

Finally, Figure~\ref{fig:beta07625} shows us comparative results at macro-micro stage $60$ for several values of the threshold coefficients $\beta$ in the interval $[0.7625, 0.7875]$. Since $\beta$ controls the ``optimal level" of ECM density for cancer cells to migrate, its variation gives us different invasion morphologies, as expected. \revdt{Again, as we noticed also in the cases of Figure~\ref{fig:mu2_15mar2016} and Figure~\ref{fig:delta_15mar2016},} a quick comparison between consecutive rows of images (from top to bottom) in Figure~\ref{fig:beta07625} seems to indicate a certain degree of consistency in the changes occurring in the tumour morphology and ``fingered" boundary deformations with respect to increasing $\beta$ parameter. This consistency aspect is currently under investigation and will form the topic of a separate research work.  

\section{Conclusion}
\label{sec:conclusion}

In this paper, we presented and developed a mathematical model of cancer invasion based on the previous  work in  \cite{Andasari_et_al_2011,Chaplain_Lolas_2005} and \cite{Dumitru_et_al_2013}. We adapted and extended the two-scale technique in \cite{Dumitru_et_al_2013} to cope with the different settings that enable the coupling of the finite-difference uPA macro-solver and a new finite element micro-solver developed for the newly proposed leading edge micro-dynamics. This enabled us to simulate the multiscale process of cancer invasion by exploring the link between the macroscopic dynamics of the spatio-temporal distribution of cancer cells and ECM taking place on a macroscopic domain, and the matrix degrading enzymes microdynamics developed on the microscopic domains that are in close proximity to tumour boundary.   

We derived a new governing law for the micro-dynamics based exclusively on the molecular mechanics of the uPA system occurring in close proximity to the tumour boundary. This is based on the dynamics of the uPA system including uPA, the inhibitor uPA-1 and plasmin, and considers the source terms related to macroscopic components, i.e. cancer cells and ECM molecules, in a collective fashion and we solve this new microscopic uPA system by finite element method. While this is sourced from within the macrodynamics via a top-down link, in our multiscale method (described in Appendix \ref{sec:The_MovingBoundaryMethod}) the microdynamics occurring at the cell-scale neighborhood of the tumour is represented back at the macroscale through a bottom-up feedback by defining the movement direction and displacement magnitude of the tissue-scale tumour boundary. 
By coupling this new microscopic governing law for the leading edge microdynamics with the macroscopic model for the uPA system and cancer invasion proposed in \cite{Andasari_et_al_2011,Chaplain_Lolas_2005}, we are able to capture an important class of multiscale dynamic interactions between cancer cells, ECM molecules, cancer associated matrix degrading enzymes, and the peritumoural tissue conditions, leading to significant changes in tumour morphology during invasion. 

From the computational simulation results of our model, we can see that the extended two-scale technique coupled with the uPA system and more specific modelling of the pericellular proteolytic activities, gives more adverse dynamics in the invading cancer. Values of the ECM initial condition, the cancer cell diffusion coefficient, the threshold coefficient, and the ECM proliferation \& degradation rates all have an impact on the deformations of the tumour boundary. \revdt{While for simulating the proposed non-dimensional model we considered functional formulations for the initial conditions and tissue thresholds, future work will attempt to assimilate the heterogeneity of ECM and as well as the peritumoural tissue conditions from imaging data. Finally, concerning the transitional probability that arises naturally at the micro-scale and intervenes in the bottom-up link between micro and macro scales, while the simulations presented in this paper are deterministic, future work will assess the stochastic character of the overall model.} 

The conclusions \revdt{that can be drawn from the qualitative results presented in this paper} are: 1) a heterogeneous ECM initial condition leads to more fingered spreading of the tumour compared with that in the homogeneous ECM; 2) in order to obtain heterogeneous patterns of cancer cells inside the tumour region, chemotaxis must be dominant to drive the cells migration; 3) the changes of threshold coefficient will definitely affect the boundary deformations, and there is a tendency that the increase of $\beta$ reduces the number of `fingers' of the interface; 4) without a proliferation term of ECM coupled with a relatively small degradation rate, deformations of the boundary show more fingering. \revdt{However, further investigations are required to analyse the observed fingering. For instance, the dependence of the width of the fingers on the size of the microscale and potentially on other regulatory parameters remains an open question and is an important objective of a future work on the propose modelling framework.} 

\revdt{It is useful to compare these results with other models of cancer growth and invasion, particularly those that adopt an alternative modelling approach such as hybrid continuum-discrete, cellular-Potts or cellular automaton, as well as partial differential equations (PDE) models. Many PDE models of solid tumour growth and development (including invasion) have modelled cell-cell adhesion at the outer boundary (or invading edge) of the tumour as a surface-tension-like force \cite{Byrne_Chaplain_1996,Byrne_Chaplain_1997}. A reduction in surface tension i.e. interpreted as a loss of cell-cell adhesion then leads to an instability at the invading edge which manifests itself in a subsequent growth consisting of finger-like protrusions. This effect was shown computationally by \cite{Cristini_et_al_2003} and \cite{macklin07} who demonstrated a range of fingering patterns as a cell adhesion parameter was varied. Similar spectra of invasive patterns have been observed when adopting either a cellular-Potts approach \cite{Popawski2009} or a hybrid continuum-discrete approach \cite{Anderson_2005} i.e., by varying a key cell-cell adhesion parameter of the model, invasive fingering patterns can either be enhanced or suppressed. In the specific case of glioma invasion, four different approaches -- PDE, cellular-Potts, lattice-gas automaton, cellular automaton -- have each investigated the role of cell-cell adhesion and compared computational simulation results with experimental data \cite{Aubert2006,Frieboes2007,Rubenstein2008,Tektonidis2011}. Given these results, we are currently extending our current model to include cell-cell adhesion. 
}

\revdt{As the main purpose of this work was to formulate a novel non-dimensional two-scale modelling platform that simultaneously explores spatio-temporal dynamics at both macro-scale (cell population level)  and micro-scale alongside the links in between the two scales, the simulations presented here have a qualitative character. Future work will explore the possibilities of dimensionalising and calibrating the proposed multiscale model with measured data at both macroscopic (tissue-scale) and microscopic (cell-scale) levels to obtain quantitative simulation that we could compare with clinical observations}.

The two-scale modelling modelling described provides a useful mathematical platform to capture and investigate processes at different levels during cancer invasion. However, further details remain to be explained on the dynamics and interactions of the tumour cell community at the macro-level and the micro-level \revdt{as well as regarding the links in between these two scales. Various avenues could be pursued to extend our modelling }by exploring for instance other ways to determine the cancer cell population macro-dynamics (for example, change the definition of the threshold function $\omega$ characterising the interaction with the peritumoural tissue) or by accounting for more complex signalling mechanisms in the establishment of more detailed micro-dynamics as well as the appropriate top-down and bottom-up links between these different levels (scales) of the invasion process.

\section{Appendix: The Two-Scale Computational Modelling Method}\label{sec:Appendix}
In the following we will briefly present the technique introduced in \dt{\cite{Dumitru_et_al_2013}} and adjust this with all the details to our new situation. For completion, we introduce here all the necessary notations and describe the defining principles that are referred to in the paper, as well as the relevant considerations and explanations concerning our new model.
\subsection{Preliminary considerations and notations}
\label{sec:notations}
\begin{figure}[!h]
\centering
\includegraphics[scale = 0.45]{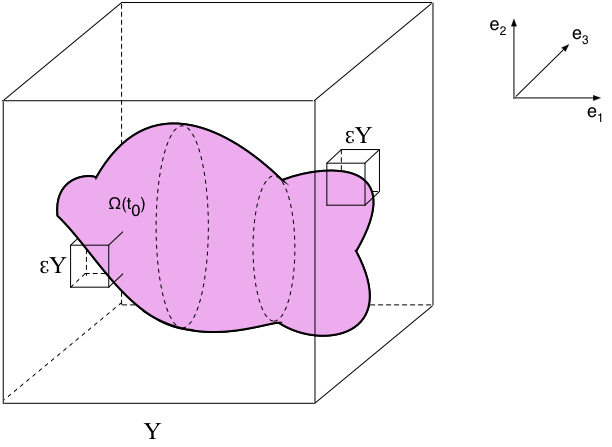}
\caption{Schematic diagram showing the cubic region Y centred at the origin $\in\mathbb{R}^3$. The dashed blue lines represent the Euclidean directions $\{e_1,e_2,e_3\}$, the pink region illustrates the cancer cluster $\Omega(t_0 )$, and the solid blue line represents the family of microscopic cubic domains $\epsilon Y$ placed at the boundary $\partial \Omega(t_0 )$.}
\label{fig:regionY}
\centering
\end{figure}

It is assumed that the domain within which the cancer and extracellular matrix exists is a maximal reference spatial cube $Y \subset \mathbb{R}^n (n=2,3) $  with its centre at  the origin. \dt{Given a fixed $ \epsilon$ representing a negative power of $2$ (\emph{i.e.}, $0<\epsilon<1$), the initial $Y$ is uniformly decomposed $\epsilon$-size cubes, $\epsilon Y$, whose union will be referred to as an $\epsilon$-resolution of $Y$. For any $\epsilon Y$ from the decomposition, the ``half-way shifted" cubes in the direction $i \bar{e}_1 + j \bar{e}_2 + k \bar{e}_3$ given by any triplet $(i,j,k) \in \{(i,j,k)|i,j,k \in\{-1,0,1\}\}$ are defined as}
\begin{equation}
\epsilon Y_{\frac{i}{2}, \frac{j}{2}, \frac{k}{2}} =  \epsilon Y + \frac{\epsilon(i \bar{e}_1 + j \bar{e}_2 + k \bar{e}_3) }{2} ,
\end{equation}
where,
\begin{align}
\bar{e}_1:&=e_1, & \bar{e}_2:&=e_2,  &\textrm{and}&, &\bar{e}_3:&=\left\{\begin{array}{ll} e_3 &  \textrm{for} \;\;\;\;N=3,   \\ 0 &  \textrm{for} \;\;\;\;N=2 ,\end{array} \right.
\end{align}
and $\{e_1,e_2,e_3\}$ is the standard Euclidean basis of $\mathbb{R}^3$. 
\dt{The family of all these $\epsilon-$cubes is denoted by $\mathcal{F}$, i.e.,}
\begin{equation}
\mathcal{F}:= \underset{i,j,k \in\{-1,0,1\}}{ \bigcup} \big\{  \epsilon Y_{\frac{i}{2}, \frac{j}{2}, \frac{k}{2}} \big |  \epsilon Y \; \textrm{is in the} \; \epsilon\textrm{-resolution of} \; Y \big \}.
\end{equation}
In Figure~\ref{fig:regionY}, the notations mentioned so far are illustrated schematically.

In order to capture mathematically the microdynamics that occur in a cell-scale neighbourhood of the tumour boundary $\partial \Omega(t_0)$, \dt{out of the initial family $\mathcal{F}$,  we will focus our attention of the subfamily denoted by $\mathcal{F}_{\Omega(t_0)}$ which consists of only the $\epsilon-$cubes that cross the interface $\partial \Omega(t_0)$ and have exactly one face included in the interior of $\Omega(t_0)$, namely}

\begin{align}
\mathcal{F}_{\Omega(t_0)} : = \{ & \epsilon Y \in \mathcal{F} | \epsilon Y \cap (Y\backslash \Omega(t_0)) \neq \emptyset , \notag \\ 
                                                          & \textrm{and} \; \epsilon Y \; \textrm{has only one face included in int}(\Omega(t_0)) \}, 
\end{align}
where int($\Omega(t_0)$) is the topological interior of $\Omega(t_0)$ with respect to the natural topology on $\mathbb{R}^n$. 
\begin{figure}[h!]
\begin{center}
\includegraphics[scale = 0.42]{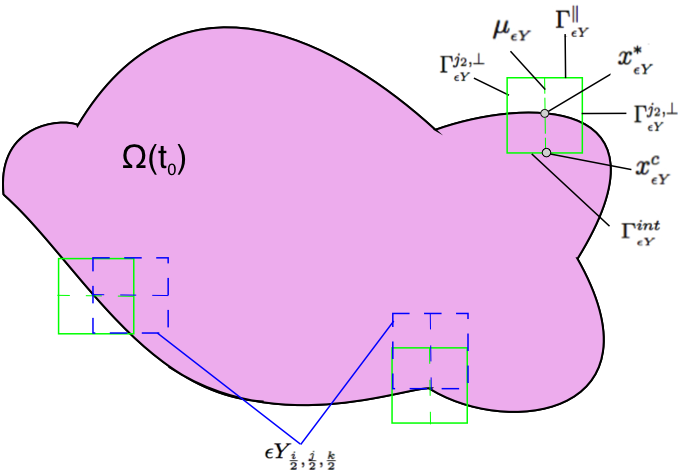}
\caption{Schematic diagram illustrating the notations introduced in \eqref{equ:faces}, \eqref{equ:points}. 
For the arbitrary microdomain $\epsilon Y \in \mathcal{P}_{\epsilon}$, we indicate with a black arrow the features: $\Gamma_{\epsilon Y}^{int}$ , $ \Gamma_{\epsilon Y}^{j_1,\bot}$, and $\Gamma_{\epsilon Y}^{j_2,\bot},\; j_1,j_2\in\{1,...,2^{N-1} \}$, $ \Gamma_{\epsilon Y}^{\|}$ , $x_{\epsilon Y}^{c}$, $\mu_{\epsilon Y}$ , and $x_{\epsilon Y}^*$ . The arbitrary cube $\epsilon Y \in \mathcal{P}_{\epsilon}^*$ is shown in green, while the corresponding half-way shifted  $\epsilon Y_{\frac{i}{2}}^{sign} \in \mathcal{P}_{\epsilon}$that are not chosen in $\mathcal{P}_{\epsilon}^*$ are shown in the blue dashed line.}
\label{fig:Omega_t0}
\end{center}
\end{figure}
\noindent \dt{In this context, for each $\epsilon Y \in \mathcal{F}_{\Omega(t_0)}$, we have the following face-notations:}
\begin{equation}
\left\{\begin{array}{l}
\Gamma_{\epsilon Y}^{int} \; \textrm{denotes the face of}\;\epsilon Y \;\textrm{that is included in int}(\Omega(t_0)),\\ \\
\Gamma_{\epsilon Y}^{j,\bot}, j=1,...,2^{N-1}, \; \textrm{denote the faces of} \;\epsilon Y \;\textrm{that are perpendicular to}\;  \Gamma_{\epsilon Y}^{int} \\ \\
\Gamma_{\epsilon Y}^{\|} \; \textrm{denotes the face of}  \;\epsilon Y \;\textrm{that is parallel to}\;  \Gamma_{\epsilon Y}^{int}.
\end{array}
\right.
\label{equ:faces}
\end{equation}
These are illustrated schematically in Figure~\ref{fig:Omega_t0}.

Furthermore, for each $\epsilon Y \in \mathcal{F}_{\Omega(t_0)}$, the topological closure of the only connected component of $\Omega(t_0) \cap \epsilon Y$ that is confined between $[\partial \Omega(t_0)]_{\epsilon Y}$ and $\Gamma_{\epsilon Y}^{int}$ is denoted by $[\Omega(t_0)]_{\epsilon Y}$. Moreover, denoting by $[\partial \Omega(t_0)]_{\epsilon Y}$ the connected component part of $\partial \Omega(t_0) \cap \epsilon Y$ with the property that 
\begin{equation}\label{non_void_inter}
[\partial \Omega(t_0)]_{\epsilon Y} \cap \Gamma_{\epsilon Y}^{j,\bot} \neq \emptyset\qquad \textrm{for any $j=1,2,...,2^{n-1}$},
\end{equation}
we can observe that $[\partial \Omega(t_0)]_{\epsilon Y}$ represents the part of $\partial \Omega(t_0) \cap \epsilon Y$ that corresponds to  $[\Omega(t_0)]_{\epsilon Y}$, and is actually the only connected component of this intersection that has property (\ref{non_void_inter}). Finally, using this observation, for the currently fixed $\epsilon$, the subfamily denoted by $\mathcal{P}_{\epsilon}$ consisting of all those $\epsilon-$cubes that have $[\Omega(t_0)]_{\epsilon Y}$ not touching $\Gamma_{\epsilon Y}^{\parallel}$ is selected as follows:
\begin{equation}
\mathcal{P}_{\epsilon} := \{ \epsilon Y \in \mathcal{F}_{\Omega(t_0)} |\; [\Omega(t_0)]_{\epsilon Y} \subset \epsilon Y \;\textrm{and}\; [\partial \Omega(t_0)]_{\epsilon Y} \cap \Gamma_{\epsilon Y}^{\parallel} = \emptyset  \}.
\end{equation}
Leaving now $\epsilon$ to take all the negative powers of $2$, the union
 \[\underset{\epsilon\in\{2^{-k}\,|\,k\in\N\}}{\bigcup}  \mathcal{P}_{\epsilon}
 \]
provides an infinite covering of  $\partial \Omega(t_0)$. Since $\partial \Omega(t_0)$ is compact, using standard compactness arguments, a finite complete sub-covering of $\partial \Omega(t_0)$ that consist only of small cubes an equal size $\epsilon^*$ is denoted by $\mathcal{P}_{\epsilon}^*$, i.e., 
\begin{equation}
\partial \Omega(t_0) \subset \underset{\epsilon Y \in \mathcal{P}_{\epsilon}^*}{\bigcup} \epsilon Y.
\end{equation}
Together with this finite complete covering $\mathcal{P}_{\epsilon}^*$ of the tumour interface $\partial \Omega(t_0)$, at each time of the tumour evolution we obtain also the size of the micro-scale $\epsilon^*$ \cite{Dumitru_et_al_2013}. For simplicity, in this paper, the size of the cell-scale $\epsilon^*$ will still be denoted by $\epsilon$. 
Finally, for each $\epsilon Y \in \mathcal{P}_{\epsilon}^*$, we distinguish the following topological details:
\begin{equation}
\left\{\begin{array}{l}
x_{\epsilon Y}^{c} \; \textrm{denotes the centre of the face}\;\Gamma_{\epsilon Y}^{int} ,\\ \\
\mu_{\epsilon Y}, \; \textrm{is the line that passes through} \; x_{\epsilon Y}^{c} \;\textrm{and is perpendicular on}\;  \Gamma_{\epsilon Y}^{int} \\ \\
x_{\epsilon Y}^{*} \in  [\partial \Omega(t_0)]_{\epsilon Y} \; \textrm{which will be referred to as the ``midpoint" of}  \; [\partial \Omega(t_0)]_{\epsilon Y},  \\
\textrm{represents the point from the intersection}\; \mu_{\epsilon Y} \cap [\partial \Omega(t_0)]_{\epsilon Y} \; \textrm{that is located} \\
\textrm{at the smallest distance with respect to} \; x_{\epsilon Y}^{c}.
\end{array}
\right.
\label{equ:points}
\end{equation}
\dt{The well-posedness of these topological features is discussed in \cite{Dumitru_et_al_2013}, and these are illustrated in Figure~\ref{fig:Omega_t0}.}


\subsection{The multiscale moving boundary approach for the proposed cancer invasion model}
\label{sec:The_MovingBoundaryMethod}

In the following, we will \dt{explain} how the set of \dt{midpoints $\{x_{\epsilon Y}^* \}_{\epsilon Y \in \mathcal{P}_{\epsilon}^*}$ defined} on the boundary of tumour at the current time moves to a set of new spatial positions \dt{$\{\widetilde{x_{\epsilon Y}^* }\}_{\epsilon Y \in \mathcal{P}_{\epsilon}^*}$} to form the new boundary at the very next time, by describing the movement of one such midpoint $x_{\epsilon Y}^* \in [\partial\Omega(t_0)]_{\epsilon Y}$ for any $\epsilon Y \in \mathcal{P}_{\epsilon}^*$.

Based on biological observations that, on any micro domain $\epsilon Y$, provided that a sufficient amount of plasmin has been produced across the invading edge and it is the pattern of the front of the advancing spatial distribution of plasmin that characterised ECM degradation, therefore it is assumed that each boundary midpoint $x_{\epsilon Y}^* \in [\partial\Omega(t_0)]_{\epsilon Y}$ will be potentially relocated in a \dt{movement} direction \dt{and by a certain displacement magnitude} dictated by the spatial distribution of plasmin obtained via the micro process on $\epsilon Y$ at the final micro-time  $\tau_f:=\Delta t$, namely, $m(\cdot, \tau_f)$. \dt{In the following, we explain how the movement direction and displacement magnitude are defined for each $x_{\epsilon Y}^* \in [\partial\Omega(t_0)]_{\epsilon Y}$}. 

\dt{For any given threshold $\delta > 0$ and any fixed $\epsilon Y \in \mathcal{P}_{\epsilon}^*$, the regularity property of Lebesgue measure \cite{Halmos1974} is used to select the first dyadic decomposition $\{ D_j\}_{j \in \mathcal{J}_{\delta}}$ of $\epsilon Y$ 
such that
\begin{equation}
\lambda \bigg( [\epsilon Y \backslash \Omega(t_0) ] \;\backslash\; \underset{ \{j \in \mathcal{J}_{\delta}\,| \,D_{j} \subset \epsilon Y \backslash \Omega(t_0) \} }{\bigcup} \mathcal{D}_j\bigg) \leq \delta.
\end{equation}
which simply means that $\epsilon Y \backslash \Omega(t_0)$ is approximated with accuracy $\delta$ by the union of all the dyadic cubes that this includes. Once this dyadic decomposition is selected, we denote by $y_{j}$ the barycenters of $D_{j}$, for all $j\in \mathcal{J}_{\delta}$.
As discussed in \cite{Dumitru_et_al_2013} for all $\epsilon Y \in \mathcal{P}_{\epsilon}^*$, this provides a resolution at which we read the further away part of the level set $\frac{1}{\lambda(\epsilon Y \backslash \Omega(t_0))} \int_{\epsilon Y \backslash \Omega(t_0)} m(y, \cdot) dy$ in the distribution of the advancing degrading enzymes $m(\cdot,\cdot)$ outside $\Omega(t_{0})$ in radial direction with respect to the midpoint $x_{\epsilon Y}^*$. Therefore, this enable us to locate dyadic pixels $D_{l}$ that support the peaks at the tip of the plasmin front with significant contribution in degrading the ECM. Hence, at the final microscopic time $\tau_{f}$, the pixels supporting these peaks are therefore selected as}
\dt{
\begin{equation}
\mathcal{\mathcal{I}_{\delta}}\!:=\! \left\{l\in \mathcal{J}_{\delta}
\left\rvert \begin{array}{l}
\exists r\in S^{1} \textrm{ such that, if the index $i\in \mathcal{J}_{\delta}$ has the properties:}\\
1)\mathcal{D}_i\cap \{x\in\R^{n}\,| x=x_{\epsilon Y}^*+\alpha r, \alpha\in\R\} \neq \emptyset, \\
2)\mathcal{D}_i\subset \epsilon Y \backslash \Omega(t_0),\\
3) \frac{1}{\lambda(\mathcal{D}_i)} \int_{\mathcal{D}_i} m(y,\tau_f) dy\geq \frac{1}{\lambda(\epsilon Y \backslash \Omega(t_0))} \int_{\epsilon Y \backslash \Omega(t_0)} m(y, \tau_f) dy , \\ 
\textrm{then} \\
l=argmax\{d(x_{\epsilon Y}^*, y_{i})\,|\,\textrm{$i\in\mathcal{J}_{\delta}$ satisfies:} 1), 2), and 3)\}                                                                                                   
\end{array} \right.                                                                                                     
\!\!\!\right\}\!,
\label{eq:selected_peek_element}                                                                                                                                                                                                   
\end{equation}
where $S^{1}\subset \R^{n}$ is represents the unit sphere, and $d(\cdot,\cdot)$ is the Euclidean distance on $\R^{n}$.} 
Thus, \dt{cumulating the driving ECM degradation forces spanned by each front peak of plasmin given by the dyadic pixels} $\mathcal{D}_{l}$ with $l \in \dt{\mathcal{I}_{\delta}}$ in the direction \dt{of the position vectors $\overrightarrow{x_{\epsilon Y}^*, y_l}$  and appropriately representing the amount of plasmin that each $D_l$ supports, the revolving direction of movement $\eta_{\epsilon Y}$ for the potential displacement of }$x_{\epsilon Y}^* $ is given by:
\begin{equation}
\eta_{\epsilon Y} = x_{\epsilon Y}^* + \nu \underset{l \in \mathcal{I}_{\delta}}{\sum} \bigg( \int_{\mathcal{D}_l} m(y,\tau_f) dy \bigg) (y - x_{\epsilon Y}^*), \nu \in [0,\infty] .
\label{eq:direction}
\end{equation}
Further, the displacement magnitude of the point $x_{\epsilon Y}^*$ is defined as:
\begin{equation}
\xi_{\epsilon Y} := \underset{l \in \mathcal{I}_{\delta}}{\sum}  \frac{\int_{\mathcal{D}_l} m(y,\tau_f) dy  }{\underset{l \in \mathcal{I}_{\delta}}{\sum} \int_{\mathcal{D}_l} m(y,\tau_f) dy  } \big | \overrightarrow{x_{\epsilon Y}^*y_l} \big | .
\label{eq:magnitude}
\end{equation}
\dt{Finally, as debated in \cite{Dumitru_et_al_2013}, although a displacement magnitude and a moving direction is derived for each $x_{\epsilon Y}^*$, this will only exercise the movement if and only if the ECM degradation were of a certain local strength. The strength of ECM degradation within $ \epsilon Y$ is explored by the transitional probability
\[
q^* : \sum\left(\underset{\epsilon Y \in \mathcal{P}_{\epsilon}^*}{\bigcup} \epsilon Y\right) \rightarrow \R_{+}
\] defined as
\begin{equation}
q^*(G):=  \frac{1}{\int_{G} m(y,\tau_f)dy} \int_{G\backslash \Omega(t_0)} m(y,\tau_f) dy, \qquad \textrm{ for all } G\in \sum\left(\underset{\epsilon Y \in \mathcal{P}_{\epsilon}^*}{\bigcup} \epsilon Y\right)
\label{eq:probability}
\end{equation}
where $\sum\left(\underset{\epsilon Y \in \mathcal{P}_{\epsilon}^*}{\bigcup} \epsilon Y\right)$ represents the Borel $\sigma-$algebra of $\underset{\epsilon Y \in \mathcal{P}_{\epsilon}^*}{\bigcup} \epsilon Y$.}
\dt{Locally, in each $\epsilon Y$,} equation \eqref{eq:probability} is in fact a quantification of the amount of plasmin in $\epsilon Y \backslash \Omega(t_0)$ relative to the total amount of plasmin concentration in $\epsilon Y$. \dt{In conjunction with the local tissue conditions, this} characterises whether the point $x_{\epsilon Y}^*$ is likely to relocate to the new spatial position $\widetilde{x_{\epsilon Y}^*}$ or not.

Now, by assuming that the point $x^*_{\epsilon Y}$ is moved to the position $\widetilde{x^*_{\epsilon Y}}$ if and only if $q^{*}(x^*_{\epsilon Y}):=q^{*}(\epsilon Y)$ exceeds a certain threshold $\omega_{ \epsilon Y} \in (0,1)$, we find that the new boundary $\partial \Omega(t_0+\Delta t)$ will be the interpolation of the following set of points:
\begin{equation}
\{ x^*_{\epsilon Y} | \epsilon Y \in \mathcal{P}_{\epsilon Y}^* \;\textrm{and}\; q(x^*_{\epsilon Y}) < \omega_{\epsilon Y}     \} \;\cup\; \{\widetilde{x^*_{\epsilon Y}} | \epsilon Y \in \mathcal{P}_{\epsilon Y}^* \;\textrm{and}\; q(x^*_{\epsilon Y}) \geq \omega_{\epsilon Y}   \}
\label{eq:newboundary}
\end{equation}

Finally before moving to the next time-step of the whole macro-micro two scale system, we replace the initial conditions of the macroscopic dynamics with the solution at the final time of the previous invasion step as follows:

\begin{align}
c_{\Omega(t_0+\Delta t)}(x,t_0) &:= c(x,t_0+\Delta t)( \chi_{_{\Omega(t_0)\backslash \underset{\epsilon Y \in \mathcal{P}_{\epsilon}^*}{\bigcup} \epsilon Y} }\ast \psi_{\gamma}), \notag \\
v_{\Omega(t_0+\Delta t)}(x,t_0) &:= v(x,t_0+\Delta t)(\chi_{_{Y \backslash \underset{\epsilon Y \in \mathcal{P}_{\epsilon}^*}{\bigcup} \epsilon Y}} \ast \psi_{\gamma}), \notag \\
u_{\Omega(t_0+\Delta t)}(x,t_0) &:= u(x,t_0+\Delta t)( \chi_{_{\Omega(t_0)\backslash \underset{\epsilon Y \in \mathcal{P}_{\epsilon}^*}{\bigcup} \epsilon Y} }\ast \psi_{\gamma}), \notag \\
p_{\Omega(t_0+\Delta t)}(x,t_0) &:= p(x,t_0+\Delta t)( \chi_{_{\Omega(t_0)\backslash \underset{\epsilon Y \in \mathcal{P}_{\epsilon}^*}{\bigcup} \epsilon Y} }\ast \psi_{\gamma}), \notag \\
m_{\Omega(t_0+\Delta t)}(x,t_0) &:= m(x,t_0+\Delta t)( \chi_{_{\Omega(t_0)\backslash \underset{\epsilon Y \in \mathcal{P}_{\epsilon}^*}{\bigcup} \epsilon Y} }\ast \psi_{\gamma}).
\label{eq:newinitial}
\end{align}

Here $\chi_{_{\Omega(t_0)\backslash \underset{\epsilon Y \in \mathcal{P}_{\epsilon}^*}{\bigcup} \epsilon Y} }\ast \psi_{\gamma}$ and $\chi_{_{Y\backslash \underset{\epsilon Y \in \mathcal{P}_{\epsilon}^*}{\bigcup} \epsilon Y} }\ast \psi_{\gamma}$ are the characteristic functions corresponding to the sets $\Omega(t_0)\backslash \underset{\epsilon Y \in \mathcal{P}_{\epsilon}^*}{\bigcup} \epsilon Y$ and $Y\backslash \underset{\epsilon Y \in \mathcal{P}_{\epsilon}^*}{\bigcup} \epsilon Y$, and choosing $\gamma \ll \frac{\epsilon}{3}$, $\psi_{\gamma}: \mathbb{R}^n \rightarrow \mathbb{R}_+$ is constructed as a smooth compact support function with $\textrm{sup}(\psi_{\gamma}) = \{z \in \mathbb{R}^n | || z ||_2 \leq \gamma \}$. This is defined by the standard mollifier $\psi: \mathbb{R}^n \rightarrow \mathbb{R}_+$, namely, 

\begin{equation}\label{mollifier_def1}
\psi_{\gamma}(x) := \frac{1}{\gamma^n} \psi(\frac{x}{\gamma}),
\end{equation}
and,
\begin{equation}\label{mollifier_def2}
\psi(x):= \left\{ \begin{array}{ll}
\frac{ \textrm{exp}(\frac{1}{x^2-1})    }{ \underset{\{z \in \mathbb{R}^n | || z ||_2 \leq \gamma \}}{\int}  \textrm{exp}(\frac{1}{z^2-1})   dz   } &\;\; \textrm{if}\;\; ||x||_2<1, \\
0&\;\; \textrm{if}\;\; ||x||_2 \geq 1, 
\end{array}
\right.
\end{equation}
Then, the invasion process will continue on the new expanded domain $\Omega(t_0)$ with the macroscopic system and the new initial conditions in \eqref{eq:newinitial} at macro-level followed by proteolytic microprocesses around its boundary, which again governs the movement of the boundary of the next time multiscale stage.   

\section{Appendix: Description of the multiscale numerical approach}
\label{sec:numerical method}
We compute and solve the multiscale model in a two-dimensional setting by using computational approach based on a finite difference scheme for macrodynamics and finite element approximation for the microdynamics occurring on each of the boundary $\epsilon Y$ microdomains. In the following subsections we detail the computational approach and present the steps of the overall multiscale algorithm. 

\subsection{The macroscopic stage of the numerical scheme}
\label{subsec:macro-numerical}
Since the macroscopic dynamics are taking place in the cube $Y$, we discretise the entire Y by considering a uniform spatial mesh of size $h:=\frac{\epsilon}{2}$, i.e., $\Delta x =\Delta y = h $. And, the time interval $[t_0,t_0+\Delta t]$ is discretised in $k$ uniformly distributed time steps, i.e., using the uniform time step $\delta \tau:=\frac{\Delta t}{k}$.
The temporal discretisation of the reaction-diffusion system \eqref{equ:Cancer} - \eqref{equ:plasmin_macro} that we used here is a second-order trapezoidal scheme, \eqref{equ:ECM}; 
while the diffusion term and haptotactic terms are approximated with a second-order midpoint rule. For instance, for the diffusion and haptotactic terms involved in  \eqref{equ:Cancer}, we approximate $\nabla \cdot (\nabla c)_{i,j}^{n}$ and $\nabla \cdot (c \nabla v)_{i,j}^{n} $ as follows:
\begin{align*}
\nabla\!\! \cdot\! (\nabla c)_{i,j}^{n} &= \textrm{div}(\nabla c)_{i,j}^{n} \\
                                                      & \simeq \frac{(c_x)_{i+\frac{1}{2},j}^{n} - (c_x)_{i-\frac{1}{2},j}^{n}}{\Delta x} +\frac{ (c_y)_{i,j+\frac{1}{2}}^{n} -   (c_y)_{i,j-\frac{1}{2}}^{n} }{\Delta y},
\end{align*}
and
\begin{align*}
\nabla\! \!\cdot \! (c \nabla v)_{i,j}^{n}\! &=\!\textrm{div} (c \nabla v)_{i,j}^{n} \\
                                                          &\simeq\!\! \frac{c_{i+\frac{1}{2},j}^{n} (v_x)_{i+\frac{1}{2},j}^{n}\! \!-\! c_{i-\frac{1}{2},j}^{n} (v_x)_{i-\frac{1}{2},j}^{n}       }{\Delta x}\! +\!\frac{ c_{i,j+\frac{1}{2}}^{n} (v_y)_{i,j+\frac{1}{2}}^{n} \!\! - \!   c_{i,j-\frac{1}{2}}^{n} (v_y)_{i,j-\frac{1}{2}}^{n}  }{\Delta y},
\end{align*}
where 
\[
\left\{ 
  \begin{array}{l l}
    c_{i,j+\frac{1}{2}}^{n} &:= \frac{c_{i,j}^{n} + c_{i,j+1}^{n}}{2}, \\ 
    c_{i,j-\frac{1}{2}}^{n} &:= \frac{c_{i,j}^{n} + c_{i,j-1}^{n}}{2}, \\ 
    c_{i+\frac{1}{2},j}^{n} &:= \frac{c_{i,j}^{n} + c_{i+1,j}^{n}}{2}, \\
    c_{i-\frac{1}{2},j}^{n} &:= \frac{c_{i,j}^{n} + c_{i-1,j}^{n}}{2}, 
  \end{array} \right.
\]
are the midpoint approximations for $c$ and\\

\begin{tabular}{lll}
$
\left\{ 
  \begin{array}{l l}
    (c_y)_{i,j+\frac{1}{2}}^{n} &:= \frac{c_{i,j+1}^{n} - c_{i,j}^{n}}{\Delta y}, \\ 
    (c_y)_{i,j-\frac{1}{2}}^{n} &:= \frac{c_{i,j}^{n} - c_{i,j-1}^{n}}{\Delta y}, \\ 
    (c_x)_{i+\frac{1}{2},j}^{n} &:= \frac{c_{i+1,j}^{n} - c_{i,j}^{n}}{\Delta x}, \\
    (c_x)_{i-\frac{1}{2},j}^{n} &:= \frac{c_{i,j}^{n} - c_{i-1,j}^{n}}{\Delta x}, 
  \end{array} \right.$,               
\;\;\;\;& \textrm{and} \;\;\;\;& 
$
\left\{ 
  \begin{array}{ll}
    (v_y)_{i,j+\frac{1}{2}}^{n} &:= \frac{v_{i,j+1}^{n} - v_{i,j}^{n}}{\Delta y}, \\ 
    (v_y)_{i,j-\frac{1}{2}}^{n} &:= \frac{v_{i,j}^{n} - v_{i,j-1}^{n}}{\Delta y}, \\ 
    (v_x)_{i+\frac{1}{2},j}^{n} &:= \frac{v_{i+1,j}^{n} - v_{i,j}^{n}}{\Delta x}, \qquad\\
    (v_x)_{i-\frac{1}{2},j}^{n} &:= \frac{v_{i,j}^{n} - v_{i-1,j}^{n}}{\Delta x}, 
  \end{array} \right.$
\end{tabular}
\\\\represent the central differences for spatial derivatives of $c$ and $v$.
Note that $n = 0,1,...,k$ are index of time step, and $(i,j)$ are spatial nodes where $i = 1,...q$ are the indices for the $x$-direction and $j = 1,...q$ are the indices for the $y$-direction. The diffusion terms in equations \eqref{equ:uPA_macro} -\eqref{equ:plasmin_macro} are approximated in the same way as it is in equation \eqref{equ:Cancer} and \eqref{equ:ECM}.

\subsection{The computational microscopic scheme and its relation to the macroscopic level}
\label{sec:The computational microscopic scheme and its relation to the macroscopic level}

In this section, we describe our computational scheme for the micro scale dynamics occurring on each microdomains $\epsilon Y \in \mathcal{P}_{\epsilon}^*$, which are cubes of size $\epsilon$ located at the boundary $\partial \Omega(t_0)$. We have each microdomain $\epsilon Y$ centred at a boundary point form the macroscopic mesh, with the neighbouring $\epsilon$-cubes staring from the centre of the current one (i.e. they are appropriately ``half-way shifted" copies of $\epsilon Y \in \mathcal{P}_{\epsilon}^*$), due to the purposely chosen macroscopic mesh size $h=\frac{\epsilon}{2}$ and the properties of the family $\mathcal{P}_{\epsilon Y}^*$. Moreover, the centre point of the microdomains are coincidentally the midpoint induced by $\epsilon Y$ on $[\partial \Omega(t_0)]_{\epsilon Y}$, i.e.  $x_{\epsilon Y}^*$.

In order to  compute the integrals in the source terms (i.e., $f_1^{\epsilon Y}$ and $f_2^{\epsilon Y}$) in the microscopic system \eqref{equ:uPA} - \eqref{equ:plasmin}, a midpoint rule is proposed and the constitutive details are given below.
Assuming that $K$ denotes a generic element domain in a finite element subdivision with either triangular or square elements of a given region $A \subset \mathbb{R}^2$, this ``midpoint rule" consists of approximating the integral of a function $f$ over $K$ as the product between the value of $f$ at the centre of mass of $K$, $K_{centre}$, and  the Lebesgue measure of $K$, namely,
\begin{equation}
\int\limits_{K} f = f(K_{centre}) \lambda(K).
\end{equation}

For an arbitrarily chosen $\epsilon Y \in \mathcal{P}_{\epsilon}^*$, we consider a finite element approach involving triangular elements on a uniform micro-mesh, which is maintained with identical structure for all the micro-domains. Further, we consider a time-constant approximation $\tilde{f}_1^{\epsilon Y}$ of $f_1^{\epsilon Y}$ on the time interval $[0, \Delta t]$. In this context, using the computed final-time values of $c(\cdot, t_0+\Delta t)$ at the macro-mesh points that are included on the current microdomain, $x_1,x_2,...,x_{P_{\epsilon Y}} \in \epsilon Y \cap \Omega(t_0)$, we take:
\begin{equation}\label{int_formula15March}
 \tilde{f}_1^{\epsilon Y} (x_s) =  \frac{1}{\lambda(B(x_s,2\epsilon)\cap \Omega(t_0))}  \int\limits_{B(x_s,2\epsilon)\cap \Omega(t_0)} c\;(x_s,t_0+\Delta t) \;dx,
\end{equation}
where $s = 1,..., P_{\epsilon Y}$, and the integrals are computed via the midpoint rule. For the rest of the points $y$ on the micro-mesh, the value of $\tilde{f}_1^{\epsilon Y}$ is obtained in terms of the set of finite element basis functions considered at the contact points, i.e. , $\{\phi_{x_s} | s = 1,..., P_{\epsilon Y} \}$. Finally, we observe that for any micro mesh point $y \in \epsilon Y$ we have two possibilities:

\paragraph{Case 1:} If there exists three overlapping points $x_{i_{1}},x_{i_{2}},x_{i_{3}} \in \{x_1,x_2,...,x_{P_{\epsilon Y}}\}$ which belongs to the same connected component of $\epsilon Y \cap \Omega(t_0)$ and $y$ belongs to the convex closure of the set, i.e. , $y \in \textrm{Conv} \{ x_{i_{1}},x_{i_{2}},x_{i_{3}} \}$, then we have:
\begin{equation}
 \tilde{f}_1^{\epsilon Y}(y) =  \tilde{f}_1^{\epsilon Y}(x_{i_1}) \phi_{x_{i_1}}(y) +\tilde{f}_1^{\epsilon Y}(x_{i_{2}}) \phi_{x_{i_2}}(y)
+\tilde{f}_1^{\epsilon Y}(x_{i_3}) \phi_{x_{i_3}}(y).
\label{eq:sourceTerm_case_1}
\end{equation}

\paragraph{Case 2:} If $y$ does not satisfies the conditions in Case 1, then we have
\begin{equation}
\tilde{f}_1^{\epsilon Y}(y) = 0.
\end{equation}

For the source term $f_2^{\epsilon Y}$, we use the same approximation method as above, except that there is only one case taken into consideration which is similar in equation \eqref{eq:sourceTerm_case_1} according to the definition of $f_2^{\epsilon Y}$. Now we could obtain the source terms $ \tilde{f}_1^{\epsilon Y} $ and $\tilde{f}_2^{\epsilon Y}$ on each microdomain $\epsilon Y$ with zero initial condition and Neumann boundary conditions and furthermore use the finite element method to solve the reaction-diffusion equations \eqref{equ:uPA} - \eqref{equ:plasmin} on $\epsilon Y$ over the time interval $[0,t_0 + \Delta t]$. Then, we use bilinear elements on a square mesh, the numerical scheme for the micro processes occurring on each $\epsilon Y$ is finally obtained by involving a trapezoidal predictor-corrector method for the time integration.

Then, for each microdomain we use the midpoint rule to compute the transitional probability described in \eqref{eq:probability}. For simplicity, now the numerical implementation of the multiscale model for cancer invasion proposed above is slightly simplified in the following way: provided that the transitional probability exceeds an associated threshold  $\omega_{\epsilon Y} \in (0, 1)$, the boundary mesh-point $x_{\epsilon Y}^*$ will move in direction $\eta_{\epsilon Y}$ to the macro-mesh point from $\partial \epsilon Y \backslash [\Omega(t_0)]_{\epsilon Y}$ that is closest (in Euclidean distance) to $x_{\epsilon Y}^*$ . If the threshold is not satisfied, then $x_{\epsilon Y}^*$ remains at the same spatial location. Therefore, the new boundary $\partial \Omega(t_0+\Delta t)$ is now obtained by the interpolation of the set of points given in \eqref{eq:newboundary} , and the computational process is continued on the new domain $\Omega(t_0+\Delta t)$ by using as a discretised version of \eqref{eq:newinitial} as a new initial condition at the macroscopic stage, i.e.,

\begin{align}
c(x_{i,j}, t_0+\Delta t) &= \left \{ \begin{array} {ll}
c_{i,j}^k,    & x_{i,j} \!\in \overline{\Omega(t_0)},\\
\frac{1}{4}(c_{i-1,j}^k\!+\!c_{i+1,j}^k\!+\!c_{i,j-1}^k\!+\!c_{i,j+1}^k), & x_{i,j}\! \in \overline{\mathbf{B}(\overline{\Omega(t_0)},h)} \backslash \overline{\Omega(t_0)},\\
0, & x_{i,j} \!\notin \overline{\mathbf{B}(\overline{\Omega(t_0)},h)},
\end{array} \right.
\label{eq:newinitial_approximation_1}
\end{align}
and,
\begin{align}
v(x_{i,j}, t_0+\Delta t) &= v_{i,j}^k, & u(x_{i,j}, t_0+\Delta t) &= u_{i,j}^k, \notag \\
p(x_{i,j}, t_0+\Delta t) &= p_{i,j}^k,  &m(x_{i,j}, t_0+\Delta t) &= m_{i,j}^k.
\label{eq:newinitial_approximation_2}
\end{align}
where $\{x_{i,j} \ i,j = 1,...,q \}$ is the macroscopic mesh in Y, $\overline{\Omega(t_0)}$ is the topological closure of $\Omega(t_0)$, and $\overline{\mathbf{B}(\overline{\Omega(t_0)},h)}$ represents the topological closure of the $h$-bundle of $\overline{\Omega(t_0)}$., i.e., $\overline{\mathbf{B}(\overline{\Omega(t_0)},h)} := \{ x \in Y | \exists z_x \in \overline{\Omega(t_0)} \; \textrm{such that}  \; || x-z_x  ||_2 \leq h            \}$.

\subsection{Overall algorithm steps}
\label{subsec:algorithm}
To sum up, the overall algorithm package of the macro-microscopic method consists of the following steps:
 
\paragraph{Step 1:} At the very begining time $t_0$, first of all, we discretise the macro-domain $[a,b] \times [c,d]$ by  
\begin{align*}
a &= x_0,\ldots,x_i = a+ i\Delta x ,\ldots, x_m = a+m\Delta x =b, \\
c  &= y_0,\ldots, y_j = c+ j\Delta y ,\ldots, y_n = c+n\Delta y =d. \\
\end{align*}
where $\Delta x = \Delta y = h$, $h:=\frac{\epsilon}{2}$ and let $a=c=0, c=d=4$. Also, we number each point on the macro-domain, record their coordinates all sorts of data of the domain that might be used later. 

\paragraph{Step 2:} Define initial conditions for cancer and ECM distribution on macro-domain:
\begin{align*}
c(x,t_0) &=: c_0(x), & x\in\Omega(t_0)\\
v(x,t_0) &=: v_0(x), & x\in\Omega(t_0)\\
u(x,t_0) &=: u_0(x), & x\in\Omega(t_0)\\
p(x,t_0) &=: p_0(x), & x\in\Omega(t_0)\\
m(x,t_0) &=: m_0(x), & x\in\Omega(t_0)
\end{align*}
where $c(x,t_0)$ is set as zero at the mesh points located outside the closure of the macroscopic domain $\Omega(t_0)$.

\paragraph{Step 3:} Start the main time loop (from time stage 1 to certain time stage $N$), and at the current time stage,
\begin{itemize}
\item[a)] Run the macro-solver, which applies the finite difference scheme mentioned above, to obtain the distribution of components in the system $c_{i,j}^{n+1}$  $v_{i,j}^{n+1}$, $u_{i,j}^{n+1}$, $p_{i,j}^{n+1}$, and $m_{i,j}^{n+1}$, where $i,j = 1,...,q$. 
\item[b)] Run the micro-solver, in which we loop over each points that was on the boundary of tumour at previous time, and at an arbitrary boundary points, 
         \begin{itemize}
         \item[i.] Define the micro-domain $\epsilon Y$ centring at the current point on the boundary,  which consists of nine points on macrodomain. For simplicity, we first construct the domain on $[0,\epsilon] \times [0,\epsilon]$ and on this domain  compute the source terms $f_1^{\epsilon Y}$ and $f_2^{\epsilon Y}$, and by interpolation, 
we uniformly decompose the domain into sixty-four square elements consists of eighty-one points in total, with the source term values and concentration values for uPA, PAI-1 and plasmin on a finer mesh (see Figure~\ref{fig:relocation}).

         \item[ii.] On the microdomain $\epsilon Y$, apply the finite element method to solve the microscopic dynamics system~\eqref{equ:uPA}-\eqref{equ:plasmin}, to obtain the spatial distribution of plasmin at the final micro-time $m(\cdot, \tau_f)$ (involving a proposed midpoint rule formula for the integral source terms $f_1^{\epsilon Y}$ and $f_2^{\epsilon Y}$, and for time integration a trapezoidal predictor-corrector), which will be used in the regulation functions of cancer cells' movement.
          
         \item[iii.] Translated the coordinates on this microdomain back to where the micro spatial position was before.

         \item[iv.] Using the trasitional probability $q^{*}$ defined in (\ref{eq:probability}), compute the invasion strength as $q^{*}(x_{\epsilon Y}^*):=q^{*}(\epsilon Y)$. 
         \item[v.] If and only if the microenvironment induced probability $q^*(x_{\epsilon Y}^*)$ is greater than some tissue threshold value $\omega_{\epsilon Y} \in (0,1)$, we further compute the direction $\eta_{\epsilon Y}$ and magnitude $\xi_{\epsilon Y}$ of the movement as described by \eqref{eq:direction} and \eqref{eq:magnitude}.
         
\end{itemize}
\item[c)] Once finishing both macro-solver and micro-solver at the current time stage, we obtained new macroscopic distribution for each components in the system $c_{i,j}^{n+1}$, $v_{i,j}^{n+1}$, $u_{i,j}^{n+1}$, $p_{i,j}^{n+1}$, and $m_{i,j}^{n+1}$; also for each midpoint $x_{\epsilon Y}^*$ on the tumour boundary, we have the possibility $q^*(x_{\epsilon Y}^*)$, direction $\eta_{\epsilon Y}$ and magnitude $\xi_{\epsilon Y}$ of their movement, therefore we could use all these information to determine the new position $\widetilde{x^*_{\epsilon Y}}$ and the points remain where they were on the cancer interface $\partial \Omega(t_0+\Delta t)$. This is schematically shown in Figure~\ref{fig:relocation} 

\begin{figure}[h!]
\begin{center}
\includegraphics[scale = 0.2]{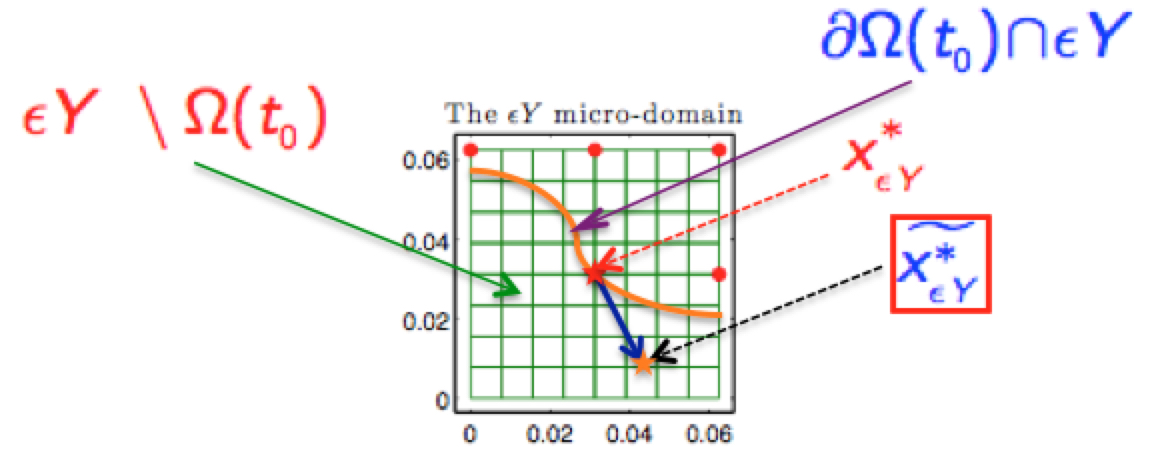}
\caption{Plot shows the relocation of one point on the boundary moves to a new spatial position in the microdomain $\epsilon Y$.  }
\label{fig:relocation}
\end{center}
\end{figure}
\revdt{where the red dots represent the discrete macro-mesh location for the where micro-scale source induced by the macro-scale was calculated via the integral formula (\ref{int_formula15March}).}\\

\item[d)] Finally, by using approximations shown in \eqref{eq:newinitial_approximation_1} and \eqref{eq:newinitial_approximation_2}, replace the initial values of cancer and ECM distribution in macroscopic dynamics  with the solution at the final time of the previous macro-step.

\end{itemize}

\paragraph{Step 4:} Using the new initial conditions for macroscopic dynamics, continue the invasion process by coupling the next-step macro-process given by the system \eqref{equ:Cancer} and \eqref{equ:plasmin_macro} on the expand domain $\Omega(t_0+\Delta t)$ with the corresponding micro processes \eqref{equ:uPA} - \eqref{equ:plasmin} occurring on its boundary, which means repeating the Step 3 above with new initial conditions for macroscopic dynamics and new boundary of cancer.
\revdt{
\section{Appendix: Table for the parameter set $\mathscr{P}$}\label{param_app}
In Table 1 we present a description of the parameters included in $\mathscr{P}$.
\begin{table}[h!]
\begin{tabular}{|lcc|}
  \hline
  && \vspace{-7pt}\\
  parameter & value & description \\ && \vspace{-7pt}\\
  \hline
  $D_c$ &  $4.3 \times 10^{-3}$ & diffusion of cancer cells \\ 
  $\chi_u$ & $3.05 \times 10^{-2}$ & chemotaxis to uPA \\
  $ \chi_p $ &$3.75 \times 10^{-2}$& chemotaxis to PAI-1\\
  $\chi_v$&$2.85 \times 10^{-2}$&haptotaxis to ECM(vitronectin)\\
  $\mu_1$&$0.25$& proliferation of cancer cells\\
  $ \delta$&$1.5$& degradation of ECM\\
  $\phi_{21}$&$0.75$&binding of uPA and PAI-1\\
  $\phi_{22}$&$0.55$&binding of PAI-1 and VN\\
  $\mu_2$&$0.01$& proliferation of ECM\\
  $D_u$&$2.5 \times 10^{-3}$& diffusion of uPA\\
  $\phi_{31}$&$0.75$&binding of uPA of PAI-1\\
  $\phi_{33}$&$0.3$& binding of uPA and uPAR\\
  $\alpha_{31}$&$0.215$& production of uPA\\
  $D_p$&$3.5 \times 10^{-3}$&diffusion of PAI-1\\
  $\phi_{41} $&$0.75$&binding of uPA and PAI-1\\
  $\phi_{42}$&$0.55$&binding of PAI-1 and VN\\
  $\alpha_{41}$&$0.5$& production of PAI-1\\
  $D_m $&$4.91 \times 10^{-3}$&diffusion of plasmin\\
  $\phi_{52}$&$0.11$&increase rate due to binding of PAI-1 and VN\\
  $\phi_{53}$&$0.75$& increase rate due to binding of  uPA and uPAR\\
  $\phi_{54}$&$0.5$&degradation of plasmin\\
        \hline
\end{tabular}
        \caption{\it The parameters in $\mathscr{P}$.}
        \label{tab:parameters}
\end{table}
}
\bibliographystyle{spmpsci}     
\bibliography{ReferenceForSummerDissertation_revised}  

\end{document}